\documentclass[10pt,twoside]{article}

\usepackage{amsmath, amssymb, amsthm, amsfonts, amsgen, stmaryrd}
\usepackage{indentfirst}
\usepackage{fancyhdr}

\numberwithin{equation}{section}

\makeatletter
\def\rightharpoonupfill@{\arrowfill@\relbar\relbar\rightharpoonup}
\newcommand{\xrightharpoonup}[2][]{\ext@arrow 0359\rightharpoonupfill@{#1}{#2}}
\makeatother

\newcommand{\ds}{\displaystyle}

\newcommand{\A}{{\mathcal{A}}}
\newcommand{\C}{{\mathcal{C}}}
\newcommand{\E}{{\mathcal{E}}}
\newcommand{\I}{{\mathcal{I}}}
\newcommand{\J}{{\mathcal{J}}}

\newcommand{\Nb}{{\mathbb{N}}}

\newcommand{\Rb}{{\mathbb{R}}}

\newcommand{\medi}{- \hskip-0.9em \int}
\newcommand{\med}{- \hskip-1.0em \int}

\let\e=\varepsilon
\let\a=\alpha

\let\g=\gamma

\let\O=\Omega
\let\o=\omega
\let\G=\Gamma

\newtheorem{thm}{Theorem}[section]
\newtheorem{defi}[thm]{Definition}
\newtheorem{rmk}[thm]{Remark}
\newtheorem{lemma}[thm]{Lemma}
\newtheorem{proposition}[thm]{Proposition}
\newtheorem{coro}[thm]{Corollary}

\begin{document}

\chead[{\sc Jean-Fran\c cois Babadjian}]{{\sc Quasistatic evolution of a
brittle thin film}} \thispagestyle{empty}

\title{{\bf Quasistatic evolution of a brittle thin film}}
\author{Jean-Fran\c cois Babadjian}
\date{}
\maketitle
\begin{abstract}
\noindent This paper deals with the quasistatic crack growth of a
homogeneous elastic brittle thin film. It is shown that the quasistatic
evolution of a three-dimensional cylinder converges, as its thickness
tends to zero, to a two-dimensional quasistatic evolution associated with
the relaxed model. Firstly, a $\Gamma$-convergence analysis is performed
with a surface energy density which does not provide weak compactness in
the space of Special Functions of Bounded Variation. Then, the asymptotic
analysis of the quasistatic crack evolution is presented in the case of
bounded solutions that is with the simplifying assumption that every minimizing
sequence is uniformly bounded in $L^\infty$.\\


\noindent {\bf MSC 2000 ($\mathcal{AMS}$): }74K30, 49J45, 74K30,
35R35, 49Q20.
\end{abstract}

\section{Introduction}\label{intro}
\noindent Following Griffith's theory of brittle fracture, the
variational model of quasistatic crack evolution proposed in
\cite{FM} is based on the competition between the elastic energy and
a surface energy which is necessary to produce a new crack or
extending a preexisting one. The classical model was plagued by a
few defects, being unable to initiate a crack, or to predict its
path during the propagation; the authors were able to overcome these
weaknesses, assuming neither a preexisting crack, nor a pre-defined
crack path. In their formulation, the time-continuous growth of the
cracks is seen as a limit of a discrete time evolution as the time
step tends to zero. The first precise mathematical justification of
this limit process in the scalar-valued case
 was  given in \cite{FL} in the framework of {\it  Special
Functions of Bounded Variation} ($SBV$). It was subsequently
generalized to the vector-valued case in \cite{DMFT1,DMFT2} (see
references therein) in the framework of {\it Generalized Special
Functions of Bounded Variation} ($GSBV$). To obtain compactness, the
authors had to either add some conservative body and surface
loadings with appropriate coerciveness as in \cite{DMFT1}, or to
impose an empirical $L^\infty$-boundness hypothesis on every
minimizing sequence as in \cite{DMFT2}. This latter assumption
permits to work in the space $SBV^p$  of  all $SBV$-functions with
$p$-integrable approximate gradient and whose jump set has finite
area, in lieu of $GSBV$ (see Definition \ref{sbvp} below). The
limit crack was defined through a new notion of convergence of
rectifiable sets, called $\sigma^p$-convergence, related to the
notion of jump sets of $SBV$ functions and based on the weak
convergence in $SBV^p$.

Before dealing with quasistatic evolutions, let us briefly describe
the static model. Let $U$ an open subset of $\Rb^n$ representing the
reference configuration of a homogeneous elastic material with
cracks, whose stored energy density is given by the function $W :
\Rb^{m \times n} \to \Rb$, where $\Rb^{m \times n}$ stands for the
set of real $m \times n$ matrices. According to Griffith's theory,
the total energy under the deformation $u : U \setminus K \to \Rb^m$
is given by
$$\E(u,K):=\int_{U \setminus K} W (\nabla u)\, d\mathcal L^n + \mathcal H^{n-1}(U \cap K),$$
where $K$ is an unknown crack. Throughout the paper, $\mathcal L^n$
and $\mathcal H^{n-1}$ denote the $n$-dimensional Lebesgue measure
and the $(n-1)$-dimensional Hausdorff measure in $\Rb^n$
respectively (in the sequel, $n$ will always be equal to $2$ or
$3$). It is often convenient to use the weak formulation of this
problem in the framework of Special Functions of Bounded Variation,
replacing the crack $K$ by the jump set $S(u)$ of the deformation $u
\in SBV(U;\Rb^m)$ and where $\nabla u$ is now the approximate
gradient of $u$. We refer to \cite{AFP,EG} for the definitions and
basic properties of Functions with Bounded Variation. It is also
usual to impose polynomial growth of order $1<p<\infty$ on $W$.
Thus, a natural space is the space $SBV^p(U;\Rb^m)$ of functions $u
\in SBV(U;\Rb^m)$ such that $\mathcal H^{n-1} (S(u) \cap U)
<+\infty$ and $\nabla u \in L^p(U;\Rb^{m \times n})$; for such
functions, the energy becomes
$$\E(u):=\E(u,S(u))=\int_U W(\nabla u)\, d\mathcal L^n + \mathcal H^{n-1}(S(u) \cap U).$$

Let us now return  to the quasistatic evolution. Adding some
appropriate boundary conditions and possibly some body and surface
loadings, at each time $t \in [0,T]$, we seek to minimize the total
energy  $(u,K) \mapsto \E(u,K)$ among all legal competitors.
Namely, $(u(t),K(t))$ is a minimum energy configuration provided that
$$\E(u(t),K(t)) \leq \E(u,K)$$
for every crack $K$ containing $K(t)$ and every deformation  $u$,
possibly discontinuous across $K$ and satisfying the boundary
conditions. More precisely, an irreversible quasistatic evolution
of minimum energy configuration is an application $[0,T] \ni t
\mapsto (u(t),K(t))$ satisfying the following conditions:
\begin{itemize}
\item[(i)]{\it Irreversibility: }$K(t)$ increases with $t$;
\item[(ii)]{\it Static equilibrium: }for every $t \in [0,T]$, the pair $(u(t),K(t))$
is a minimum energy configuration;
\item[(iii)]{\it Nondissipativity: }the function $t \mapsto \E(u(t),K(t))$ is absolutely continuous.
\end{itemize}
Condition (iii) express the conservation of the energy in the sense
that the derivative of the internal energy $\E(u(t),K(t))$ is equal
to the power of the applied forces.

Sometimes a natural small parameter, denoted by $\e>0$, is involved
in the model and one should look at the behavior of the energy and
of the quasistatic evolution $(u^\e(t),K^\e(t))$ when $\e$ tends to
zero. The notion of $\Gamma$-convergence -- see \cite{DM} for a
complete treatment on that subject -- has proved useful in
investigating the variational convergence for static problems; in
the present context of variational evolution, it remains to see
whether such a process  is compatible with the evolution. A problem
of this type have been studied in \cite{GP} where the authors proved
a stability result of variational models of quasistatic crack
evolution under $\Gamma$-convergence in the antiplanar case. Indeed,
they have shown that $(u^\e(t),K^\e(t))$ converges, in a certain
sense, to a quasistatic evolution of the relaxed model. They had to
define a new variational notion of convergence for rectifiable sets
that they called $\sigma$-convergence. The limit crack is thus seen
as the $\sigma$-limit of $K^\e(t)$ while the limit deformation is
nothing but the weak $SBV^p(U)$-limit of $u^\e(t)$. In the same
spirit, \cite{G1} investigates a notion of quasistatic evolution for
the elliptic approximation of the Mumford-Shah functional and proves
its convergence to a quasistatic growth of brittle fractures in
linearly elastic bodies.

In this paper we treat  cylindrical bodies whose thickness becomes
arbitrarily small and ask how the quasistatic evolution behaves
when the small parameter, the thickness of the cylinder, tends to
zero. It leads us to perform a 3D-2D asymptotic analysis of a thin
film, taking into account the possibility of fracture. An abundant
literature exists on dimensional reduction problems and we point the
reader to e.g. \cite{BhaFF,BraFo,BFF,LDR} and references therein.

Let $\O_\e:=\o \times (-\e,\e)$ be a three-dimensional cylinder of
basis  $\o$, a bounded open subset of $\Rb^2$ and of thickness
$2\e$, representing the reference configuration of a homogeneous
elastic body containing some cracks. The study of the static problem
relies on the computation of the $\G$-limit of the sequence of
functionals associated to the internal energy, as $\e$ tends to
zero. Namely, we want to look at the behavior of the functional
$$SBV^p(\O_\e;\Rb^3) \ni u \mapsto \int_{\O_\e} W(\nabla u)\, d\mathcal L^3 + \mathcal H^2(S(u) \cap \O_\e)$$
when $\e \to 0$, in the sense of $\G$-convergence. It is convenient
to rescale the problem on the unit cylinder $\O:=\O_1$. Denoting by
$x_\a:=(x_1,x_2)$ the in-plane variable, we set
$v(x_\a,x_3/\e)=u(x_\a,x_3)$. Replacing $u$ by $v$ in the previous
energy, changing  variables and dividing by $\e$, we are led to
studying the equivalent sequence of functionals
$$SBV^p(\O;\Rb^3) \ni v \mapsto \int_\O W\left( \nabla_\a v \Big|\frac{1}{\e} \nabla_3 v\right)d\mathcal L^3 + \int_{S(v) \cap \O} \left|\left( \left(\nu_v \right)_\a \Big|\frac{1}{\e} \left(\nu_v \right)_3 \right)\right| d\mathcal H^2,$$
where, from now on, $\nabla_\a$ (resp. $\nabla_3$) will stand for
the (approximate) differential with respect to $x_\a$ (resp. $x_3$),
$\xi=(\xi_\a|\xi_3)$ for some matrix $\xi \in \Rb^{3 \times 3}$ and
$z=(z_\a|z_3)$ for some vector $z \in \Rb^3$. Such studies have been
performed in \cite{BhaFF} where the authors imposed to the crack to
live on the mid-surface of the plate and also in \cite{BFLM,BraFo}.
Unlike the latest references, we are working with materials
satisfying Griffith's principle. This was not previously taken into
account because the authors forced the surface energy to grow
linearly with respect to the jump of the deformation. This
restriction was due to the application of Integral Representation
Theorems in $SBV(\O;\Rb^3)$ (see e.g. Theorem 2.4 in \cite{BCP} or
Theorem 1 in \cite{BFLM}) which hold for such surface energies.
Furthermore, it was necessary to obtain a bound in $BV(\O;\Rb^3)$ of
any minimizing sequence. We are convinced that such a behavior is
not very realistic because it does not permit to take into account
neither Griffith nor Barenblatt theory in which the surface energy
behaves asymptotically like a constant: the toughness. However, in
\cite{BraFo}, Example 2.10, the authors suggest a way to recover a
brittle elastic material obeying  Griffith's law by a singular
perturbation approach (see also \cite{BDV}, Section 8).

The originality of our work comes from the fact that we are directly
dealing with Griffith surface energy, which was proscribed in
\cite{BFLM,BraFo} as explained above. We propose here an alternative
proof of this result, deriving the $\G$-limit by an argument very
closed to that used in the proof of Theorem 2.1 in \cite{FF}.
Indeed, to overcome the lack of compactness in $BV(\O;\Rb^3)$, we
use a regular truncation function. It permits us to  show that for a
deformation belonging to $L^\infty(\O;\Rb^3)$, there is no loss of
generality in requiring that any minimizing sequence is uniformly
bounded in $L^\infty(\O;\Rb^3)$ and thus, weak
$SBV^p(\O;\Rb^3)$-compactness follows. The two-dimensional surface
energy remains of Griffith type whereas the bulk energy follows the
one obtained in \cite{BFLM,BraFo}: the stored energy density is
given by the function $\mathcal Q W_0 : \Rb^{3 \times 2} \to \Rb$
where $W_0 (\overline \xi):= \inf \left\{ W(\overline \xi |z) : z
\in \Rb^3 \right\}$ and $\mathcal QW_0$ is the
2D-quasiconvexification of $W_0$. Note that we obtain exactly the
same energy density as in \cite{LDR}, in which the authors were
treating healthy materials, that is without any crack.

Then, it remains to pass to the limit in the quasistatic evolution.
We need a notion of convergence of rectifiable sets, in the spirit
of $\sigma^p$-convergence introduced in \cite{DMFT1}, but better
adapted to the dimensional reduction problems. Intuitively, any
limit crack should be a one-dimensional set, that is a
two-dimensional set which is invariant by translation in the $x_3$
direction. It seems that $\sigma^p$-convergence do not permit to
obtain such limit cracks. This is why we are led to defining a new
notion of convergence (see Definition \ref{3d2dconv}) similar to
that of $\sigma^p$-convergence (see Definition 4.1 in \cite{DMFT1}),
but for the fact that we impose to any sequence of test functions to
have an approximate scaled gradient uniformly bounded in $L^p$.
Nevertheless it is not sufficient to pass to the limit in the
surface term because we need to have compactness for sequences of
cracks with bounded scaled surface energy. This difficulty is
overcome thanks to Proposition \ref{helly}. The last tool is a Jump
Transfer Theorem (see Theorem 2.1 in \cite{FL}) stated in a rescaled
version in Theorem \ref{JTT} and whose proof uses a De Giorgi type
slicing argument.

As in \cite{DMFT2}, we will assume that any minimizing sequence
involved in the quasistatic evolution, remains bounded in
$L^\infty(\O;\Rb^3)$ uniformly in $\e>0$ and in time $t \in [0,T]$
(see (\ref{Hyp}) and the proof of Lemma \ref{eto0}). We do not
attempt to justify this hypothesis which appears naturally in the
scalar case by a truncation argument, whenever the prescribed
boundary deformation is also bounded in $L^\infty(\O \times
[0,T];\Rb)$, uniformly in $\e>0$. It will yield weak
$SBV^p(\O;\Rb^3)$-compactness and thus allow to define a limit
deformation field. Note that we could also have  taken the path to
add appropriate body and surface loadings, as in \cite{DMFT1}, with
the right order of magnitude, in which case we should be working in
a suitable subspace of $GSBV(\O;\Rb^3)$ instead of
$SBV^p(\O;\Rb^3)$. We insist on the fact that this boundness
hypothesis is completely empirical but has the advantage of avoiding
technical difficulties connected with the
space $GSBV(\O;\Rb^3)$.\\

As mentioned before, we adopt the following
\begin{defi}\label{sbvp} For any open set $U \subset \Rb^n$,
$$SBV^p(U;\Rb^m) := \{w \in SBV(U;\Rb^m): \nabla w \in L^p(U;\Rb^{m \times n}) \text{ and }
\mathcal H^{n-1}(S(w) \cap U) <\infty \}.$$
\end{defi}
Throughout the text, letters as $C$ stand for generic constants
which may vary from line to line. We will denote by $\widetilde
\subset$ and $\widetilde =$ inclusion and equality of sets up to a
set of zero $\mathcal H^{n-1}$-measure respectively. As for
notation, $\mathcal M_b(U)$ stands for the space of signed Radon
measures with finite total variation. If $\mu \in \mathcal M_b(U)$
and $E$ is a Borel subset of $\Rb^n$, we will denote by
$\mu_{\lfloor E}$ the restriction of the measure $\mu$ to $E$ that
is, for every Borel subset $B$ of $\Rb^N$, $\mu_{\lfloor E}(B)=\mu(E
\cap B)$. For any $\mu$-measurable set $A \subset \Rb^n$ and any
$\mu$-measurable function $f:\Rb^n \to \Rb$, we write $\medi_A f\,
d\mu:= \mu(A)^{-1}\int_A f\, d\mu$ for the average of $f$ over $A$.
For the remainder of the paper, strong convergence will always be
denoted by $\to$, whereas weak (resp. weak-$*$) convergence will be
denoted by $\rightharpoonup$ (resp. $\xrightharpoonup[]{*}$). We
recall that a sequence $\{u_k\} \subset SBV^p(U;\Rb^m)$ is said to
converge weakly to some $u \in SBV^p(U;\Rb^m)$, and we write $u_k
\rightharpoonup u$ in $SBV^p(U;\Rb^m)$, if
$$\left\{
\begin{array}{l}
u_k \to u \text{ in } L^1(U;\Rb^m),\\
\\
u_k \xrightharpoonup[]{*} u \text{ in } L^\infty(U;\Rb^m),\\
\\
\nabla u_k  \rightharpoonup \nabla u \text{ in }L^p(U;\Rb^{m \times n}),\\
\\
(u_k^+-u_k^-) \otimes \nu_{u_k} \mathcal H^{n-1}_{\lfloor S(u_k)}
\xrightharpoonup[]{*} (u^+-u^-) \otimes \nu_{u} \mathcal
H^{n-1}_{\lfloor S(u)} \text{ in }\mathcal M_b(U;\Rb^{m \times n}).
\end{array}
\right.$$


\section{Formulation of the problem}\label{model}

\subsection{The physical configuration}\label{physical}

\noindent {\bf Reference configuration: }Let $\o$ be a bounded open
subset of $\Rb^2$ with Lipschitz boundary and define $\O_\e:= \o
\times (-\e,\e)$. The cylinder $\overline \O_\e$ denotes the
reference configuration of a homogeneous elastic body with cracks,
whose stored energy density is given by the function $W : \Rb^{3
\times 3} \to \Rb$. We assume $W$ to be a $\mathcal C^1$ and
quasiconvex function satisfying standard $p$-growth and
$p$-coercivity conditions ($1<p<\infty$): there exists $ \beta >0$
such that
\begin{equation}\label{pg1}
\frac{1}{\beta}|\xi|^p -\beta \leq W(\xi) \leq \beta (1+|\xi|^p), \text{ for all } \xi \in \Rb^{3 \times 3}.
\end{equation}
Let $p'=(p-1)/p$ be the conjugate exponent of $p$. In particular
(see \cite{D}), the derivative of $W$, denoted by $\partial W$,
satisfies some $(p-1)$-growth condition, namely
\begin{equation}\label{pg2}
|\partial W(\xi)| \leq \beta (1+|\xi|^{p-1}), \text{ for all } \xi \in \Rb^{3 \times 3}.
\end{equation}

\noindent{\bf Boundary conditions: }Let $\o' \subset \Rb^2$ be such
that $\overline \o \subset \o'$ and $\O'_\e:= \o' \times (-\e,\e)$
be the enlarged cylinder. We submit the body to a ``smooth'' given
deformation $$G^\e \in W^{1,1}([0,T];W^{1,p}(\O'_\e;\Rb^3))$$ on the
lateral boundary $\partial \o \times (-\e,\e)$ of the cylinder
$\O_\e$. For  homogeneity, we assume that $G^\e$ satisfies
\begin{equation}\label{CLe}
\|G^\e \|_{W^{1,1}([0,T];W^{1,p}(\O'_\e;\Rb^3))}\leq C \e^{1/p}.
\end{equation}
As  in \cite{DMFT2,FL}, we express the lateral boundary condition at
time $t$ by requiring that $U=G^\e(t)$ a.e. on $[\o' \setminus
\overline \o] \times (-\e,\e)$ whenever $U \in SBV^p(\O'_\e;\Rb^3)$
is a kinematically admissible deformation field.  Since the
Dirichlet boundary condition is only prescribed on the lateral
boundary, we do not need, as in \cite{DMFT2,FL}, to extend the whole
cylinder but only its lateral part, which can be trivially done by
extending the base domain $\o$ to $\o'$. As a consequence, any
admissible crack will necessarily be contained in $\overline \o
\times (-\e,\e)$.  We will assume that the remainder of the
boundary, $\o \times \{-\e,\e\}$, is traction free so that the
prescribed boundary deformation is the only driving mechanism. Since
the admissible cracks never run into $\o \times \{-\e,\e\}$, there
is no need, in contrast to \cite{DMFT2,FL}, to remove the part $\o
\times  \{-\e,\e\}$ from the surface energy of the crack, so that
the surface energy of a crack associated to a test function
$V \in SBV^p(\O'_\e;\Rb^3)$ will be exactly $\mathcal H^2(S(V) )$.\\

\noindent{\bf Initial conditions: }As in \cite{FL}, we consider a
body without any preexisting crack (see Remark \ref{crack}). If
$U^\e_0 \in SBV^p(\O'_\e;\Rb^3)$ is a given initial deformation
satisfying $U^\e_0=G^\e(0)$ a.e. on $[\o' \setminus \overline \o]
\times (-\e,\e)$,  we suppose that the Griffith equilibrium
condition is satisfied, that is $U^\e_0$ minimizes
$$V \mapsto \int_{\O_\e}W(\nabla V)\, dx + \mathcal H^2(S(V) )$$
among $\{V \in SBV^p(\O'_\e;\Rb^3) : \quad V=G^\e(0)\text{ a.e. on }[\o' \setminus \overline \o] \times (-\e,\e)\}$.\\

\noindent{\bf Quasistatic evolution: }If we follow word for word the
arguments developed   in \cite{DMFT1,DMFT2},  applying them to
$\o'\times (-\e,\e)$ in place of  what is used in those references,
namely $\Omega''_\e$, a Lipschitz extension of $\o\times (-\e,\e)$
with $\overline \o\times [-\e,\e]\subset \Omega''_\e$, we get   the
existence of a crack $K^\e(t) \subset \overline \o \times (-\e,\e)$
increasing in time and a deformation field $U^\e(t) \in
SBV^p(\O'_\e;\Rb^3)$ such that
\begin{itemize}
\item $U^\e(0)=U^\e_0$, $K^\e(0) \; \widetilde =\; S(U^\e_0)$ and $U^\e(0)$ minimizes
$$V \mapsto \int_{\O_\e}W(\nabla V)\, dx + \mathcal H^2(S(V))$$
among $\{V \in SBV^p(\O'_\e;\Rb^3) : \quad V=G^\e(0)\text{ a.e. on }[\o' \setminus \overline \o] \times (-\e,\e)\}$;

\item for any $t \in (0,T]$, $S(U^\e(t)) \; \widetilde{\subset}\;  K^\e(t)$ and $U^\e(t)$ minimizes
$$V \mapsto \int_{\O_\e}W(\nabla V)\, dx + \mathcal H^2(S(V) \setminus K^\e(t) )$$
among $\{ V\in SBV^p(\O'_\e;\Rb^3) : \quad V=G^\e(t) \text{ a.e. on }[\o' \setminus \overline \o] \times (-\e,\e)\}$;

\item the total energy
$$E_\e(t):=\int_{\O_\e}W(\nabla U^\e(t))\, dx + \mathcal H^2(K^\e(t))$$
is absolutely continuous with respect to the time $t$ and
$$E_\e(t)=E_\e(0)+\int_0^t \int_{\O_\e}\partial W(\nabla U^\e(\tau))\cdot \nabla \dot G^\e(\tau)\, dx\, d\tau.$$
\end{itemize}


\subsection{The rescaled configuration}\label{scaled}

\noindent As usual in dimension reduction, we perform a scaling so
as to study an equivalent problem stated on a cylinder of unit
thickness $\O:= \o \times I$,
where $I:=(-1,1)$. Let $\O':=\o' \times I$ be the enlarged rescaled cylinder.\\

\noindent{\bf Boundary conditions (rescaled): }Let
$g^\e(t,x_\a,x_3)=G^\e(t,x_\a,\e x_3)$ be the rescaled boundary
deformation. Changing variables, (\ref{CLe}) implies that
\begin{equation}\label{CL}
\|g^\e\|_{W^{1,1}([0,T];L^p(\O';\Rb^3))} + \left\| \left( \nabla_\a
g^\e \Big| \frac{1}{\e} \nabla_3 g^\e \right)
\right\|_{W^{1,1}([0,T];L^p(\O';\Rb^{3 \times 3}))} \leq C.
\end{equation}
In fact, we will further restrict $g^\e$ in (\ref{CLcvforte}) by
requiring the strong convergence of both $g^\e$ and its scaled
gradient together with a uniform bound in $L^\infty(\O' \times
[0,T];\Rb^3)$. In particular, the limit $g \in
W^{1,1}([0,T];W^{1,p}(\O';\Rb^3))$ satisfies $\nabla_3 g=0$. In the
sequel, we will identify the space $SBV^p(\o;\Rb^3)$ (resp.
$W^{1,p}(\o;\Rb^3)$, $L^p(\o;\Rb^3)$) with those functions $v \in
SBV^p(\O;\Rb^3)$ (resp. $W^{1,p}(\O;\Rb^3)$, $L^p(\O;\Rb^3)$) such
that $D_3 v =0$ in the sense of distributions.
As a consequence $g \in W^{1,1}([0,T];W^{1,p}(\o';\Rb^3))$.\\

\noindent{\bf Initial conditions (rescaled): }Let $u^\e_0(x_\a,x_3)=U^\e_0(x_\a,\e x_3)$
the rescaled initial deformation, then $u_0^\e \in
SBV^p(\O';\Rb^3)$, $u_0^\e = g^\e(0)$ a.e. on $[\o' \setminus
\overline \o] \times I$ and $u_0^\e$ minimizes
$$v \mapsto \int_{\O} W\left( \nabla_\a v \Big| \frac {1}{\e} \nabla_3 v \right) \, dx + \int_{S(v)} \left| \left( \left(\nu_v \right)_\a \Big| \frac{1}{\e} \left(\nu_v \right)_3 \right) \right| \, d \mathcal H^2,$$
among $\{ v \in SBV^p(\O';\Rb^3) : \quad v=g^\e(0) \text{ a.e. on }[\o' \setminus \overline \o] \times I\}$.\\

\noindent{\bf Quasistatic evolution (rescaled): }We set
$$u^\e(t,x_\a,x_3)=U^\e(t,x_\a,\e x_3) \text{ and }\G^\e(t)= \big\{ x \in \overline \o \times I : \quad (x_\a,\e x_3) \in K^\e(t) \big\}.$$
Then, $u^\e(t) \in SBV^p(\O';\Rb^3)$ and $\G^\e(t) \subset \overline \o \times I$ is increasing in time. Moreover,

\begin{itemize}
\item $u^\e(0)=u^\e_0$, $\G^\e(0) \; \widetilde = \; S(u_0^\e);$
\item for all $t \in (0,T]$, $S(u^\e(t)) \; \widetilde{\subset}\;  \G^\e(t)$ and  $u^\e(t)$ minimizes
\begin{equation}\label{min_t_e}
v \mapsto \int_{\O} W\left( \nabla_\a v \Big| \frac {1}{\e} \nabla_3 v \right) \, dx + \int_{S(v) \setminus \G^\e(t) } \left| \left( \left(\nu_v \right)_\a \Big| \frac{1}{\e} \left(\nu_v \right)_3 \right) \right| \, d \mathcal H^2,
\end{equation}
among $\{ v \in SBV^p(\O';\Rb^3) : \quad v=g^\e(t) \text{ a.e. on }[\o' \setminus \overline \o] \times I\}$;

\item the total energy
\begin{eqnarray*}
\E_\e(t):=\int_{\O} W\left( \nabla_\a u^\e(t) \Big| \frac {1}{\e}
\nabla_3 u^\e(t) \right)\, dx+ \int_{\G^\e(t)} \left| \left( \left(
\nu_{\G^\e(t)} \right)_\a \Big| \frac{1}{\e} \left( \nu_{\G^\e(t)}
\right)_3 \right) \right| \, d \mathcal H^2 \end{eqnarray*} is
absolutely continuous in time and
\begin{eqnarray}\label{Etot}
 \E_\e(t)=\E_\e(0) + \int_0^t \! \int_{\O} \partial W\left( \nabla_\a u^\e(\tau)
\Big| \frac{1}{\e} \nabla_3 u^\e(\tau) \right)  \cdot  \left(
\nabla_\a \dot g^\e(\tau) \Big| \frac{1}{\e} \nabla_3 \dot
g^\e(\tau) \right) dx\, d\tau.
\end{eqnarray}
\end{itemize}

We would like to perform an asymptotic analysis of this quasistatic
evolution when the thickness $\e$ tends to zero. To this end, we
start by stating a $\G$-convergence result in order to guess how the
energy is behaving through the dimensional reduction. In fact, we
will prove in Section \ref{gammaconvergence} that the functional
$$SBV^p(\O;\Rb^3) \ni v \mapsto \int_{\O} W\left( \nabla_\a v \Big| \frac {1}{\e} \nabla_3 v \right) \, dx + \int_{S(v)} \left| \left( \left(\nu_v \right)_\a \Big| \frac{1}{\e} \left(\nu_v \right)_3 \right) \right| d \mathcal
H^2$$ $\G$-converges for the strong $L^1(\O;\Rb^3)$-topology to
$$SBV^p(\o;\Rb^3) \ni v \mapsto 2\int_\o \mathcal Q W_0(\nabla_\a
v)\, dx_\a +2 \mathcal H^1(S(v)),$$ where $W_0 : \Rb^{3 \times 2}
\to \Rb$ is defined by $W_0(\overline \xi):= \inf \{W(\overline
\xi|z) : z \in \Rb^3\}$,
$$\mathcal QW_0(\overline \xi):= \inf_{\varphi \in W^{1,\infty}_0(Q';\Rb^3)} \int_{Q'}W_0(\overline \xi +
\nabla_\a \varphi(x_\a))\, dx_\a$$ is the 2D-quasiconvexification of
$W_0$ and $Q':=(0,1)^2$ is the unit square of $\Rb^2$. A
two-dimensional quasistatic evolution relative to the boundary data
$g(t)$ for the relaxed model is an application $[0,T] \ni t \mapsto
(u(t),\g(t))$ such that $u(t) \in SBV^p(\o';\Rb^3)$, $u(t)=g(t)$
a.e. on $\o' \setminus \overline \o$, $\g(t) \subset \overline \o$
and the three following properties hold:
\begin{itemize}
\item[(i)] {\it Irreversibility}: $\g(t_1) \;\widetilde \subset \; \g(t_2)$, for every $0 \leq t_1 \leq t_2 \leq T$;
\item[(ii)] {\it Minimality}: $S(u(0)) \; \widetilde = \; \g(0)$, $u(0)$ minimizes
$$v \mapsto 2\int_\o \mathcal Q W_0(\nabla_\a v)\, dx_\a+ 2 \mathcal H^1(S(v)),$$
among $\{ v \in SBV^p(\o';\Rb^3) : \quad v=g(0) \text{ a.e. on }\o'
\setminus \overline \o\}$ and, for every $t \in (0,T]$, $S(u(t)) \;
\widetilde \subset \; \g(t)$ and $u(t)$ minimizes
$$v \mapsto 2\int_\o \mathcal Q W_0(\nabla_\a v)\, dx_\a+ 2 \mathcal H^1(S(v)\setminus \g(t)),$$
among $\{ v \in SBV^p(\o';\Rb^3) : \quad v=g(t) \text{ a.e. on }\o'
\setminus \overline \o\}$;
\item[(iii)] {\it Nondissipativity}: The total energy
$$\E(t):=2\int_\o \mathcal Q W_0(\nabla_\a u(t))\, dx_\a +2\mathcal H^1(\g(t))$$
is absolutely continuous in time and
$$\E(t)=\E(0)+2\int_0^t \int_\o \partial(\mathcal Q W_0)(\nabla_\a u(\tau)) \cdot \nabla \dot g(\tau)\, dx_\a\, d\tau.$$
\end{itemize}

Note that the previous equality makes sense because we will prove in
Proposition \ref{classC1} that $\mathcal QW_0$ is of class $\C^1$
provided $W$ is also of class $\C^1$. The main result of this paper
is that the three-dimensional quasistatic evolution converges, in
the sense detailed below, towards a two-dimensional quasistatic
evolution associated with the $\G$-limit model. This is formally
expressed as the following

\begin{thm}\label{qse}
For all $\e>0$, let $[0,T] \ni t \mapsto (u^\e(t),\G^\e(t))$ be
a three-dimensional (rescaled) quasistatic evolution relative to the
boundary data $g^\e(t)$. Assume that
\begin{equation}\label{CLcvforte}
\left\{
\begin{array}{l}
\ds \sup_{\e >0}\|g^\e\|_{L^\infty(\O' \times [0,T];\Rb^3))} <+\infty,\\
\ds g^\e \to g \text{ in }W^{1,1}([0,T];W^{1,p}(\O';\Rb^3)), \quad
\frac{1}{\e}\nabla_3 g^\e \to H\text{ in
}W^{1,1}([0,T];L^p(\O';\Rb^3))
\end{array}
\right.
\end{equation}
and that
\begin{equation}\label{Hyp}
\sup_{\e>0 } \|u^\e(t)\|_{L^\infty(\O';\Rb^3)} <+\infty
\end{equation}
uniformly with respect to  $t \in [0,T]$. Then, there
exists a two-dimensional quasistatic evolution $[0,T] \ni t \mapsto
(u(t),\g(t))$ relative to the boundary data $g(t)$ for the relaxed
model and a sequence $\{\e_n\} \searrow 0^+$ such that for every $t
\in [0,T]$,
\begin{itemize}
\item  $\G^{\e_n}(t)$ converges to $\g(t)$ in the sense of Definition \ref{3d2dconv} and $u^{\e_n}(0) \rightharpoonup u(0)$ in $SBV^p(\O';\Rb^3)$;
\item $u^{\e_{n_t}}(t) \rightharpoonup u(t) \text{ in }SBV^p(\O';\Rb^3)$, for some $t$-dependent subsequence $\{\e_{n_t}\} \subset \{\e_n\}$;
\item the total energy $\E_{\e_n}(t)$ converges to $\E(t)$ and more precisely,
$$\left\{\begin{array}{l}
\ds \int_{\O} W\left(\nabla_\a u^{\e_{n}}(t)\Big| \frac{1}{\e_{n}} \nabla_3 u^{\e_{n}}(t) \right) \, dx \to 2\int_{\o} \mathcal Q W_0(\nabla_\a u(t))\, dx_\a,\\
\\
\ds \int_{\G^{\e_{n}}(t)} \left| \left( \left( \nu_{\G^{\e_{n}}(t)} \right)_\a \Big| \frac{1}{\e_{n}} \left( \nu_{\G^{\e_{n}}(t)} \right)_3 \right) \right| \, d \mathcal H^2 \to 2\mathcal H^1(\g(t)).
\end{array}\right.$$
\end{itemize}
\end{thm}

As it has been discussed in the introduction, it seems that the
right functional setting for this kind of problems is that of {\it
Generalized Special Functions of Bounded Variation} (see
\cite{DMFT1}) adding appropriate conservative surface and body
forces. Here, as in \cite{DMFT2}, we will only deal with {\it
Special Functions of Bounded Variation}, imposing, without any
justification, that the minimizing fields are bounded in
$L^\infty(\O;\Rb^3)$, uniformly in time (see (\ref{Hyp}) and the
proof of Lemma \ref{eto0}). Note that this assumption is
automatically satisfied in the scalar case by a truncation argument,
provided that $g^\e$ is also bounded in $L^\infty(\O' \times
[0,T];\Rb)$. This will let us get a bound in $BV(\O;\Rb^3)$ for
any sequence with bounded energy and thus  apply Ambrosio's
Compactness Theorem in $SBV$ (Theorem 4.8 in \cite{AFP}). In a
$GSBV$ context,  we could state a  result similar to that in
\cite{DMFT1} upon adding some appropriate conservative body and
surface loadings with the right order of magnitude.   Adequate
coerciveness assumptions would allow us to get rid of the empirical
$L^\infty$-bound (\ref{Hyp}) because we would have natural
compactness in a suitable subspace of
$GSBV(\O;\Rb^3)$ and we would then obtain a membrane limit model.\\

Let us first state a compactness result in $SBV^p(\O;\Rb^3)$ which
ensures that any limit deformation field does not depend on the
$x_3$ variable.

\begin{lemma}\label{comp}
Let $\{\e_n\} \searrow 0^+$ and $\{u_n\} \subset SBV^p(\O;\Rb^3)$ such that
\begin{eqnarray}\label{bound}
\sup_{n \in \Nb} \left\{ \|u_n\|_{L^\infty(\O;\Rb^3)} + \int_{\O}
\left|\left( \nabla_\a u_n \Big| \frac{1}{\e_n} \nabla_3 u_n \right)
\right|^p dx +\!\! \int_{S(u_n)} \left|\left( \left(\nu_{ u_n}
\right)_\a \Big| \frac{1}{\e_n} \left(\nu_{u_n} \right)_3 \right)
\right| d\mathcal H^2 \right\} <+\infty.
\end{eqnarray}
Then, there exists a subsequence $\{\e_{n_k}\} \subset \{\e_n\}$ and a function $u \in SBV^p(\o;\Rb^3)$ such that $u_{n_k} \rightharpoonup u$ in $SBV^p(\O;\Rb^3)$.
\end{lemma}

\noindent{\it Proof. }In view of (\ref{bound}), we have in particular
$$\sup_{n \in \Nb} \left\{ \|u_n\|_{L^\infty(\O;\Rb^3)} + \| \nabla u_n \|_{L^p(\O;\Rb^{3 \times 3})} + \mathcal H^2( S(u_n) ) \right\} <+\infty.$$
Thus, according to Ambrosio's Compactness Theorem (Theorem 4.8 in \cite{AFP}), we can find a subsequence $\{\e_{n_k}\} \subset \{\e_n\}$ and a function $u \in SBV^p(\O;\Rb^3)$ such that $u_{n_k} \rightharpoonup u$ in $SBV^p(\O;\Rb^3)$. Let us show that $u$ does not depend on $x_3$. Indeed, (\ref{bound}) implies that
$$\frac{1}{\e_{n_k}}\left[ \int_\O |\nabla_3 u_{n_k}|^p\, dx + \int_{S(u_{n_k})} \left|\left( \nu_{u_{n_k}} \right)_3 \right|\, d\mathcal H^2 \right] \leq C.$$
By lower semicontinuity of
$$v \mapsto \int_\O |\nabla_3 v|^p \, dx + \int_{S(v)} \left|\left( \nu_v \right)_3\right| \, d\mathcal H^2$$
with respect to the $SBV^p(\O;\Rb^3)$-convergence (see
\cite{BFLM,BraFo}), we deduce that $\nabla_3 u=0$ $\mathcal
L^3$-a.e. in $\O$ and that $\big(\nu_u \big)_3=0$ $\mathcal
H^2$-a.e. in $S(u)$. Since $u \in SBV^p(\O;\Rb^3)$, it implies that
$D_3 u=0$ in the sense of Radon measures and thus $u \in
SBV^p(\o;\Rb^3)$.
\hfill$\Box$\\

Note that the $L^\infty$-bound in (\ref{bound}) will  follow
from assumption (\ref{Hyp}), an assumption that we do not attempt to justify, while
the two other bounds will appear naturally in the energy estimates.


\section{A $\G$-convergence result}\label{gammaconvergence}

\noindent This section is devoted to the study of the static
problem. Let us define $\I_\e : BV(\O;\Rb^3) \to \overline \Rb$ by
$$\I_\e(u):=\left\{
\begin{array}{cl}
\ds \int_{\O} W\left(\nabla_\a u\Big|\frac{1}{\e}\nabla_3 u\right)dx
+ \int_{S(u)} \left| \left( \left(\nu_u\right)_\a \Big|\frac{1}{\e}
\left( \nu_u \right)_3 \right)\right| d\mathcal
H^2 & \text{if }u \in SBV^p(\O;\Rb^3),\\
&\\
\ds +\infty & \text{otherwise}
\end{array}
\right.$$

\noindent Then, the following $\G$-convergence result holds:

\begin{thm}\label{gammaconv}
Let $\o$ be a bounded open subset of $\Rb^2$ and $W : \Rb^{3 \times
3} \to \Rb$ be a continuous function satisfying (\ref{pg1}). Then
the functional $\I_\e$ $\G$-converges for the strong
$L^1(\O;\Rb^3)$-topology towards $\J : BV(\O;\Rb^3) \to \overline
\Rb$ defined by
$$\J(u):=\left\{
\begin{array}{cl}
\ds 2 \int_\o \mathcal Q W_0( \nabla_\a u)\, dx_\a +2 \mathcal H^1(S(u)) & \text{if }  u \in SBV^p(\o;\Rb^3),\\
&\\
\ds  +\infty & \text{otherwise}.
\end{array}
\right.$$
\end{thm}

\begin{rmk}\label{qcvx}{\rm
It has been noted in \cite{LDR} p. 556 (see also \cite{BraFo} p.
306) that the function $\Rb^{3 \times 3} \ni (\overline \xi|\xi_3)
\mapsto \mathcal Q W_0(\overline \xi)$ is quasiconvex: for any
$\varphi \in W^{1,\infty}_0(Q;\Rb^3)$,
$$\mathcal Q W_0(\overline \xi) \leq \int_{Q} \mathcal Q W_0(\overline \xi + \nabla_\a \varphi(y))\, dy,$$
where $Q:=(0,1)^3$ is the unit square of $\Rb^3$.
}
\end{rmk}

We first localize our functionals on $\A(\o)$, the family of open
subsets of $\o$. Define $\I_\e$ and $\J : BV(\O;\Rb^3) \times \A(\o)
\to \overline \Rb$ by
\begin{equation}\label{IeA}
\I_\e(u;A):=\left\{
\begin{array}{cl}
\begin{array}{l} \ds \int_{A \times I}  \ds W\left(\nabla_\a u\Big|\frac{1}{\e}\nabla_3 u\right)dx  \\
\hspace{2.0cm} \ds + \int_{S(u) \cap [A \times I]} \left| \left(
\left(\nu_u\right)_\a \Big|\frac{1}{\e} \left( \nu_u \right)_3
\right)\right| d\mathcal H^2
\end{array} & \text{if }u \in SBV^p(A \times I;\Rb^3),\\
&\\
\ds +\infty & \text{otherwise}
\end{array}
\right.
\end{equation}
and
\begin{equation}\label{JA}
\J(u;A):=\left\{
\begin{array}{cl}
\ds 2 \int_A \mathcal Q W_0( \nabla_\a u)\, dx_\a +2 \mathcal H^1(S(u) \cap A) & \text{if }  u \in SBV^p(A;\Rb^3),\\
&\\
\ds  +\infty & \text{otherwise}.
\end{array}
\right.
\end{equation}
For every sequence $\{\e_j\} \searrow 0^+$ and all $(u;A) \in
BV(\O;\Rb^3) \times \A(\o)$, we define the $\G$-lower limit by
\begin{equation}\label{IA}
\I (u;A):= \inf_{\{u_j\}} \left\{ \liminf_{j \to +\infty} \; \I_{\e_j}(u_j;A) : u_j \to u \text{ in }L^1(A \times I;\Rb^3)\right\}.
\end{equation}
Let $\mathcal R(\o)$ be the countable subfamily of $\A(\o)$ obtained
by taking every  finite union of open cubes of $\o$, centered at
rational points and with rational edge length. Theorem 8.5 in
\cite{DM} together with a diagonalization argument imply  the
existence of a subsequence $\{\e_n\} \equiv \{\e_{j_n}\}$ such that,
for any $C \in \mathcal R(\o)$ (or $C=\o$), $\I(\cdot;C)$ is the
$\G$-limit of $\I_{\e_n}(\cdot;C)$ for the strong $L^1(C \times
I;\Rb^3)$-topology. To prove Theorem \ref{gammaconv}, it is enough
to show that $\I(u;\o)=\J(u;\o)$.


\subsection{A truncation argument}\label{truncation}

\noindent The main problem with definition (\ref{IA}) of $\I$ is
that any minimizing sequence is not necessarily bounded in
$BV(\O;\Rb^3)$ and thus, not necessarily weakly convergent in this
space. What is missing is either a bound on the jump part of the
derivative (this bound was appearing naturally in \cite{BraFo,BFLM}
because the authors were considering a surface energy density with
linear growth with respect to the jump of the deformation), or an
$L^\infty$-bound on the minimizing sequences. In the spirit of
\cite{FF}, we define for all $(u;A) \in BV(\O;\Rb^3) \times \A(\o)$
$$\I_\infty (u;A):= \inf_{\{u_n\}} \left\{ \liminf_{n \to +\infty} \; \I_{\e_n}(u_n;A) : u_n \to u \text{ in }L^1(A \times I;\Rb^3),\; \sup_{n \in \Nb} \|u_n\|_{L^\infty(A \times I;\Rb^3)} <+\infty \right\}.$$
Obviously, we have that $\I(u;A) \leq \I_\infty (u;A)$. In fact, we will prove that both functionals coincide if the deformation belongs to $BV(\O;\Rb^3) \cap L^\infty(\O;\Rb^3)$ in Lemma \ref{bd} below. Thus for a deformation $u \in BV(\O;\Rb^3) \cap L^\infty(\O;\Rb^3)$, strong $L^1(\O;\Rb^3)$-convergence and weak $BV(\O;\Rb^3)$-convergence are, in a sense, equivalent for the computation of the $\G$-limit. This will permit us to prove in Lemma \ref{measure} that for such deformations $\I_\infty(u;\cdot)$ is a Radon measure on $\o$, absolutely continuous with respect to sum of the Lebesgue measure and of the restriction of the Hausdorff measure to $S(u)$. By Lebesgue's Decomposition Theorem, it will be enough to identify the Radon-Nikodym derivatives of $\I_\infty(u;\cdot)$ with respect to $\mathcal L^2$ and $\mathcal H^1_{\lfloor S(u)}$. This will be done for the proof of the upper bound in Lemma \ref{gs1} by a blow up argument. The general case will be treated in Lemma \ref{gs2} thanks to a truncation argument. The lower bound will be proved in Lemma \ref{gi} using the quasiconvexity properties of $\mathcal QW_0$ (see Remark \ref{qcvx}) and a lower semicontinuity result in $SBV^p(\O;\Rb^3)$.

\begin{lemma}\label{bd}
For any $C \in \mathcal R(\o)$ (or $C=\o$) and all $u \in BV(\O;\Rb^3) \cap L^\infty(\O;\Rb^3)$, $\I(u;C) = \I_\infty (u;C)$.
\end{lemma}

\noindent{\it Proof. }It is enough to show that $\I(u;C) \geq \I_\infty (u;C)$ for all $u \in BV(\O;\Rb^3) \cap L^\infty(\O;\Rb^3)$. If $\I(u;C)=+\infty$, the result is obvious, thus there is no loss of generality in assuming that $\I(u;C)<+\infty$. By the very definition of the $\G$-limit, there exists a sequence $u_n \to u$ in $L^1(C \times I;\Rb^3)$ such that
\begin{equation}\label{gammalim}
\I(u;C) = \lim_{n \to +\infty}\; \I_{\e_n}(u_n;C).
\end{equation}
Since $\I(u;C)<+\infty$, in view of (\ref{IeA}), we deduce that, for $n$ large enough, $u_{n} \in SBV^p(C \times I;\Rb^3)$. We consider only those $n$'s.

Let us define a smooth truncation function $\varphi_i \in \mathcal C^1_c(\Rb^3 ; \Rb^3)$ satisfying
\begin{equation}\label{fii}
\varphi_i(z)=\left\{
\begin{array}{rcl}
z & \text{if} & |z|<e^i,\\
0 & \text{if} & |z| \geq e^{i+1}
\end{array}\right. \quad \text{ and }\quad |\nabla \varphi_i(z)|\leq 1.
\end{equation}
Let $w_n^i:=\varphi_i(u_{n})$, thanks to the Chain Rule formula,
Theorem 3.96 in \cite{AFP}, $w_n^i \in SBV^p(C \times I;\Rb^3)$ and
\begin{equation}\label{wik}
\left\{\begin{array}{l}
\|w_n^i\|_{L^\infty(C \times I;\Rb^3)} \leq e^i,\\
\\
S(w_n^i) \subset S(u_{n}),\\
\\
\nabla w_n^i(x) = \nabla \varphi_i(u_n(x)) \circ \nabla u_{n}(x) \quad \mathcal L^3\text{-a.e. on }C \times I.
\end{array}\right.
\end{equation}
Since $u \in L^\infty(\O;\Rb^3)$, we can choose $i$ large enough so that $u=\varphi_i(u)$, thus according to (\ref{fii})
\begin{equation}\label{L1}
\|w_n^i-u\|_{L^1(C \times I;\Rb^3)}  =   \|\varphi_i(u_{n})-\varphi_i(u)\|_{L^1(C \times I;\Rb^3)} \leq  \|u_{n}-u\|_{L^1(C \times I;\Rb^3)}.
\end{equation}
The growth condition (\ref{pg1}), (\ref{fii}) and (\ref{wik}) imply that
\begin{eqnarray*}
&&\int_{C \times I}W\left( \nabla_\a w_n^i \Big|\frac{1}{\e_{n}} \nabla_3 w_n^i\right)dx\\
&&\hspace{2.0cm}= \int_{\{u_{n} < e^i\} }W\left( \nabla_\a u_{n} \Big|\frac{1}{\e_{n}} \nabla_3 u_{n} \right)dx \\
&&\hspace{2.5cm} +\int_{\{e^i \leq  u_{n} < e^{i+1}\} }W\left( \nabla \varphi_i(u_n) \circ \nabla_\a u_{n} \Big|\frac{1}{\e_{n}} \nabla \varphi_i(u_n) \circ \nabla_3 u_{n} \right)dx\\
&&\hspace{3.0cm} + W(0) \mathcal L^3(\{u_{n} \geq e^{i+1}\})\\
&&\hspace{2.0cm}\leq \int_{C \times I}W\left( \nabla_\a u_{n} \Big|\frac{1}{\e_{n}} \nabla_3 u_{n} \right)dx + \beta \mathcal L^3(\{u_{n} \geq e^{i}\})\\
&&\hspace{2.5cm} + \beta \int_{\{e^i \leq u_{n} < e^{i+1}\} }\left|\left( \nabla_\a u_{n} \Big|\frac{1}{\e_{n}}  \nabla_3 u_{n} \right)\right|^p dx\\
&&\hspace{2.0cm} \leq \int_{C \times I}W\left( \nabla_\a u_{n} \Big|\frac{1}{\e_{n}} \nabla_3 u_{n} \right)dx + \frac{\beta}{e^i}\|u_{n}\|_{L^1(C \times I;\Rb^3)}\\
&&\hspace{2.5cm} + \beta \int_{\{e^i \leq u_{n} < e^{i+1}\} }\left|\left( \nabla_\a u_{n} \Big|\frac{1}{\e_{n}}  \nabla_3 u_{n} \right)\right|^p dx,\\
\end{eqnarray*}
where we have used Chebyshev's inequality. Let $M \in \Nb$, a summation for $i=1$ to $M$  implies using (\ref{wik}) and the fact that $\nu_{w_n^i}(x)=\pm \nu_{u_{n}}(x)$ for $\mathcal H^2$-a.e. $x \in S(w_n^i)$,
\begin{eqnarray*}
&&\frac{1}{M} \sum_{i=1}^M \left[ \int_{C \times I}W\left( \nabla_\a w_n^i \Big|\frac{1}{\e_{n}} \nabla_3 w_n^i\right)dx + \int_{S(w_n^i) \cap [C \times I] } \left| \left( \left(\nu_{w_n^i}\right)_\a \Big| \frac{1}{\e_{n}} \left(\nu_{w_n^i}\right)_3 \right)\right| d\mathcal H^2 \right]\\
&&\hspace{0.5cm} \leq \int_{C \times I} W\left( \nabla_\a u_{n} \Big|\frac{1}{\e_{n}} \nabla_3 u_{n} \right)dx + \int_{S(u_{n}) \cap [C \times I]} \left| \left( \left(\nu_{u_{n}}\right)_\a \Big| \frac{1}{\e_{n}} \left(\nu_{u_{n}}\right)_3 \right)\right| d\mathcal H^2+ \frac{c}{M},
\end{eqnarray*}
where $$c=\beta\sup_{n \in \Nb}\|u_{n}\|_{L^1(C \times I;\Rb^3)}
\sum_{i \geq 1}\frac{1}{e^i} + \beta \sup_{n \in \Nb}\left\|\left(
\nabla_\a u_{n} \Big|\frac{1}{\e_{n}}  \nabla_3 u_{n} \right)
\right\|^p_{L^p(C \times I;\Rb^{3 \times 3})} < +\infty.$$ We may
find some $i_n \in \{1,\ldots,M\}$ such that
\begin{eqnarray}\label{last}
&&\int_{C \times I}W\left( \nabla_\a w_n^ {i_n} \Big|\frac{1}{\e_{n}} \nabla_3 w_n^{i_n}\right)dx + \int_{S(w_n^{i_n}) \cap [C \times I]} \left| \left( \left(\nu_{w_n^{i_n}}\right)_\a \Big| \frac{1}{\e_{n}} \left(\nu_{w_n^{i_n}}\right)_3 \right)\right| d\mathcal H^2 \nonumber\\
&& \leq \int_{C \times I}W\left( \nabla_\a u_{n} \Big|\frac{1}{\e_{n}} \nabla_3 u_{n} \right)dx + \int_{S(u_{n}) \cap [C \times I]} \left| \left( \left(\nu_{u_{n}}\right)_\a \Big| \frac{1}{\e_{n}} \left(\nu_{u_{n}}\right)_3 \right)\right| d\mathcal H^2+ \frac{c}{M}.
\end{eqnarray}
Set $v_n:=w_n^{i_n}$, thus in view of (\ref{L1}), $v_n \to u$ in $L^1(C \times I;\Rb^3)$. Moreover, (\ref{wik}) implies that $$\|v_n\|_{L^\infty(C \times I;\Rb^3)} \leq e^{i_n} \leq e^M.$$
Finally, by virtue of (\ref{gammalim}) and (\ref{last}),
$$\I(u;C) +\frac{c}{M} \geq \liminf_{n \to +\infty} \I_{\e_{n}}(v_n;C) \geq \I_\infty (u;C).$$
The proof is achieved upon letting $M$ tend to $+\infty$
\hfill$\Box$\\

To prove the upper bound in Lemma \ref{gs1} below, we need a little
bit more than  the only continuity condition imposed on $W$, namely
a $p$-Lipschitz condition. If $W$ was quasiconvex, this property
would be immediate. Since we do not want to restrict too much the
stored energy density, we will show that there is no loss of
generality in assuming $W$ to be quasiconvex. Let $\mathcal Q W$ be
the 3D-quasiconvexification of $W$ defined by
$$\mathcal QW(\xi):= \inf_{\varphi \in W^{1,\infty}_0(Q;\Rb^3)} \int_Q W(\xi + \nabla \varphi(x))\, dx \quad \text{ for all }\xi \in \Rb^{3 \times 3}.$$

\begin{lemma}\label{QW}
For all $u \in SBV^p(\o;\Rb^3) \cap L^\infty(\o;\Rb^3)$, the value of $\I_\infty(u;\o)$ does not change if we replace $W$ by $\mathcal QW$ in (\ref{IeA}).
\end{lemma}

\noindent {\it Proof. }We denote $\I_\e^{\mathcal Q}$ (resp.
$\I^{\mathcal Q}$, $\I^{\mathcal Q}_\infty$), the value of $\I_\e$
(resp. $\I$, $\I_\infty$) with $\mathcal Q W$ instead of $W$ in
(\ref{IeA}). By the same arguments as above, we may assume that
$\{\e_n\}$ is a subsequence of $\{\e_j\}$ such that $\I^{\mathcal
Q}(u;\o)=\I^{\mathcal Q}_\infty(u;\o)$ is the $\G$-limit of
$\I_{\e_n}^{\mathcal Q}(u;\o)$, for every $u \in BV(\O;\Rb^3) \cap
L^\infty(\O;\Rb^3)$.

Let $u \in SBV^p(\o;\Rb^3) \cap L^\infty(\o;\Rb^3)$, since $W \geq \mathcal QW$, we obviously have $\I_\infty(u;\o) \geq \I^{\mathcal Q}_\infty(u;\o)$. Let us prove the converse inequality. By the definition of the $\G$-limit, we may find a sequence $\{u_n\} \subset SBV^p(\O;\Rb^3)$ such that $u_n \to u$ in $L^1(\O;\Rb^3)$, $\sup_{n \in \Nb}\| u_n \|_{L^\infty(\O;\Rb^3)} < +\infty$ and
$$\I^{\mathcal Q}_\infty(u;\o) = \lim_{n \to +\infty} \left\{ \int_\O \mathcal Q W\left( \nabla_\a u_n \Big|\frac{1}{\e_n} \nabla_3 u_n \right)dx  + \int_{S(u_n)} \left|\left( \big( \nu_{u_n} \big)_\a \Big|\frac{1}{\e_n} \big( \nu_{u_n} \big)_3 \right) \right| d \mathcal H^2 \right\}.$$
We undo the scaling by letting $v_n(x_\a,x_3):=u_n(x_\a,x_3/\e_n)$.
Then $v_n \in SBV^p(\O_{\e_n};\Rb^3)$,
$$\frac{1}{\e_n} \int_{\O_{\e_n}} |v_n - u|\, dx \to 0, \quad \sup_{n \in \Nb}\| v_n \|_{L^\infty(\O_{\e_n};\Rb^3)} < +\infty$$
and
\begin{equation}\label{scaling}
\I^{\mathcal Q}_\infty(u;\o) = \lim_{n \to +\infty} \frac{1}{\e_n} \left\{ \int_{\O_{\e_n}} \mathcal Q W ( \nabla v_n )\, dx  +  \mathcal H^2(S(v_n)) \right\}.
\end{equation}
For all $n \in \Nb$, Theorem 8.1 in \cite{BCP} and Proposition 2.8
in \cite{FF} yield the existence of a sequence $\{v_{n,k}\}_{k \in
\Nb} \subset SBV^p(\O_{\e_n};\Rb^3)$ satisfying $v_{n,k} \to v_n$ in
$L^1(\O_{\e_n};\Rb^3)$ as $k \to +\infty$,
\begin{equation}\label{Braides}
\int_{\O_{\e_n}} \mathcal Q W ( \nabla v_n )\, dx  +  \mathcal
H^2(S(v_n))= \lim_{k \to +\infty} \left\{ \int_{\O_{\e_n}} W (
\nabla v_{n,k} )\, dx  +  \mathcal H^2(S(v_{n,k})) \right\}
\end{equation}
and $\sup_{k \in \Nb}\| v_{n,k} \|_{L^\infty(\O_{\e_n};\Rb^3)} <
+\infty$. Since the previous bound is of the form $e^M$ (see
\cite{FF} p. 417), for some constant $M>0$ independent of $n$, this
last relation holds uniformly with respect to $n \in \Nb$. Gathering
(\ref{scaling}) and (\ref{Braides}), we get
$$\I^{\mathcal Q}_\infty(u;\o) = \lim_{n \to +\infty} \lim_{k \to +\infty} \frac{1}{\e_n} \left\{ \int_{\O_{\e_n}} W ( \nabla v_{n,k} )\, dx  +  \mathcal H^2(S(v_{n,k})) \right\}.$$
Let $u_{n,k}(x_\a,x_3):=v_{n,k}(x_\a,\e_n x_3)$, then $\{u_{n,k}\} \subset SBV^p(\O;\Rb^3)$,
$$\lim_{n \to +\infty} \lim_{k \to +\infty} \int_{\O} |u_{n,k} - u|\, dx = 0, \quad \sup_{n,k \in \Nb}\| u_{n,k} \|_{L^\infty(\O;\Rb^3)} < +\infty$$
and \begin{eqnarray*} &&\I^{\mathcal Q}_\infty(u;\o) = \lim_{n \to
+\infty}\lim_{k \to +\infty} \left\{ \int_\O W\left( \nabla_\a
u_{n,k} \Big|\frac{1}{\e_n} \nabla_3 u_{n,k}
\right)dx\right.\\
&&\hspace{6.0cm}\left. + \int_{S(u_{n,k})} \left|\left( \big(
\nu_{u_{n,k}} \big)_\a \Big|\frac{1}{\e_n} \big( \nu_{u_{n,k}}
\big)_3 \right) \right| d \mathcal H^2 \right\}. \end{eqnarray*} By
a diagonalization argument, we can find a sequence $k_n \nearrow
+\infty$ such that upon denoting $w_n := u_{n,k_n}$, then $w_n \in
SBV^p(\O;\Rb^3)$, $w_n \to u$ in $L^1(\O;\Rb^3)$, $\sup_{n \in
\Nb}\| w_n \|_{L^\infty(\O;\Rb^3)} < +\infty$ and
\begin{eqnarray*}
\I^{\mathcal Q}_\infty(u;\o) & = & \lim_{n \to +\infty} \left\{
\int_\O W\left( \nabla_\a w_n \Big|\frac{1}{\e_n} \nabla_3 w_n
\right)dx  + \int_{S(w_n)} \left|\left( \big( \nu_{w_n} \big)_\a
\Big|\frac{1}{\e_n} \big( \nu_{w_n} \big)_3 \right) \right| d
\mathcal H^2 \right\}\\
& \geq & \I_{\infty}(u;\o).
\end{eqnarray*}
\hfill$\Box$

\begin{rmk}\label{WQW}
{\rm
From Comments on Theorem 8. (iii) p. 560 in \cite{LDR}, we always have $\mathcal Q W_0=\mathcal Q ((\mathcal Q W)_0)$. As a consequence, by Lemma \ref{QW}, we may assume without loss of generality, upon replacing $W$ by $\mathcal Q W$, that $W$ is quasiconvex. In particular (see \cite{D}), the following $p$-Lipschitz condition holds,
\begin{equation}\label{pg1bis}
|W(\xi_1)- W(\xi_2)| \leq \beta (1+|\xi_1|^{p-1}+|\xi_2|^{p-1})|\xi_1-\xi_2|, \text{ for all } \xi_1,\; \xi_2 \in \Rb^{3 \times 3}.
\end{equation}
}
\end{rmk}


\subsection{Integral Representation of the $\Gamma$-limit}\label{interep}

\noindent Lemma \ref{bd} is essential for the proof of the following
result because it allows us to replace strong
$L^1(\O;\Rb^3)$-convergence of any minimizing sequence by strong
$L^p(\O;\Rb^3)$-convergence, where $1<p<\infty$ is the same exponent
as in (\ref{pg1}).

\begin{lemma}\label{measure}
For all $u \in SBV^p(\o;\Rb^3) \cap L^\infty(\o;\Rb^3)$, $\I_\infty(u;\cdot)$ is the restriction to $\A(\o)$ of a Radon measure absolutely continuous with respect to $\mathcal L^2 + \mathcal H^1_{\lfloor S(u)}$.
\end{lemma}

\noindent {\it Proof. }Let $u \in SBV^p(\o;\Rb^3) \cap
L^\infty(\o;\Rb^3)$. The $p$-growth condition (\ref{pg1}) implies
that
\begin{equation}\label{harraps}
\I_\infty(u;A) \leq 2\beta \int_A (1+|\nabla_\a u|^p)\, dx_\a +
2\mathcal H^1(S(u) \cap A).
\end{equation}
Thus, thanks to e.g. Lemma 7.3 in \cite{BFF}, it is enough to show
the existence of a Radon measure $\hat \mu$ on $\Rb^2$ such that for
every $A$, $B$ and $C \in \A(\o)$,
\begin{itemize}
\item[(i)] $\I_\infty(u;A) \leq \I_\infty(u;A \setminus \overline C) +\I_\infty(u;B)$ if $\overline C \subset B \subset A$;
\item[(ii)] for any $\delta>0$, there exists $C_\delta \in \A(\o)$ such that $\overline C_\delta \subset A$ and $\I_\infty(u;A \setminus \overline C_\delta) \leq \delta$;
\item[(iii)] $\I_\infty(u;\o) \geq \hat \mu(\Rb^2)$;
\item[(iv)] $\I_\infty(u;A) \leq \hat \mu(\overline A)$.
\end{itemize}

Since  $u \in SBV^p(\o;\Rb^3)$, $\mathcal H^1(S(u))<+\infty$, thus
$\mathcal H^1_{\lfloor S(u)}$ is a Radon measure. Then, for any
$\delta>0$, there exists $C_\delta \in \A(\o)$ such that $\overline
C_\delta \subset A$ and
$$ 2 \beta \int_{A \setminus \overline C_\delta}(1+|\nabla_\a u|^p)\, dx_\a + 2 \mathcal H^1(S(u) \cap [A \setminus \overline C_\delta]) \leq \delta.$$
Thus, thanks to the growth condition (\ref{pg1}), we have
$\I_\infty(u;A \setminus \overline C_\delta) \leq \delta$ and item
(ii) holds true. Furthermore, by Lemma \ref{bd} and the definition
of the $\G$-limit, there exists a sequence $\{u_n\} \subset
SBV^p(\O;\Rb^3)$ such that $u_n \to u$ in $L^1(\O;\Rb^3)$ and
$\I_{\e_n}(u_n;\o) \to \I_\infty(u;\o)$. Denoting by
$$\mu_n := W\left( \nabla_\a u_n \Big|\frac{1}{\e_n} \nabla_3 u_n \right) \mathcal L^3_{\lfloor \O} + \left|\left( \big( \nu_{u_n} \big)_\a \Big|\frac{1}{\e_n} \big( \nu_{u_n} \big)_3 \right) \right| \mathcal H^2_{\lfloor S(u_n) \cap \O},$$
for a subsequence of $\{\e_n\}$ (not relabeled), there exists a
Radon measure $\mu$ such that $\mu_n \xrightharpoonup[]{*} \mu$ in
$\mathcal M_b(\Rb^3)$. Let $\hat \mu (B):=\mu (B \times  [-1,1] )$,
for every Borel set $B \subset \Rb^2$. Thus, $\hat \mu (\Rb^2) \leq
\I_\infty(u;\o)$ and item (iii) follows. Moreover, for every $A \in
\A(\o)$,
$$\I_\infty(u;A)  \leq  \liminf_{n \to +\infty} \I_{\e_n}(u_n;A)  \leq  \limsup_{n \to +\infty} \mu_n( \overline A \times [-1,1]) \leq \mu (\overline A \times [-1,1]) = \hat \mu(\overline A)$$
which establishes item (iv). We now show the subadditivity condition
expressed in item (i). For any $\eta >0$, we can find a sequence
$\{v_n\} \subset SBV^p([A \setminus \overline C] \times I;\Rb^3)$
such that $v_n \to u$ in $L^1([A \setminus \overline C] \times
I;\Rb^3)$,
\begin{equation}\label{sup1}
\sup_{n \in \Nb}\|v_n\|_{L^\infty([A \setminus \overline C] \times I;\Rb^3)} <+\infty
\end{equation}
and
$$\liminf_{n \to +\infty}\I_{\e_n}(v_n;A \setminus \overline C) \leq \I_\infty(u;A \setminus \overline C) + \eta.$$
In particular, $v_n \to u$ in $L^p([A \setminus \overline C] \times I;\Rb^3)$ and we may extract a subsequence $\{\e_{n_k}\} \subset \{\e_n\}$ for which
\begin{equation}\label{l1}
\lim_{k \to +\infty}\I_{\e_{n_k}}(v_{n_k};A \setminus \overline C) \leq \I_\infty(u;A \setminus \overline C) + \eta.
\end{equation}
Let $R_0 \in \mathcal R(\o)$ satisfy $C \subset \subset R_0
\subset\subset B$, thus, since $\I(u;R_0)$ is the $\G$-limit of
$\I_{\e_{n_k}}(u;R_0)$, thanks to Lemma \ref{bd}, there exists a
sequence $\{u_k\} \subset SBV^p(R_0 \times I;\Rb^3)$ such that $u_k
\to u$ in $L^1(R_0 \times I;\Rb^3)$,
\begin{equation}\label{sup2}
\sup_{k \in \Nb}\|u_k\|_{L^\infty(R_0 \times I;\Rb^3)} <+\infty
\end{equation}
and
\begin{equation}\label{l2}
\I_{\e_{n_k}}(u_k;R_0) \to \I_\infty(u;R_0).
\end{equation}
In particular, we have $u_k \to u$ in $L^p(R_0 \times I;\Rb^3)$.
According to the $p$-coercivity condition (\ref{pg1}) the following
sequence of  Radon measures
\begin{eqnarray*}
\lambda_k := \left( 1+\left|\left( \nabla_\a v_{n_k} \Big| \frac{1}{\e_{n_k}} \nabla_3 v_{n_k} \right) \right|^p \right)
 \mathcal L^3_{\lfloor (R_0 \setminus \overline C) \times I} +  \left| \left( \big(\nu_{v_{n_k}}\big)_\a \Big|
 \frac{1}{\e_{n_k}} \big(\nu_{v_{n_k}}\big)_3 \right) \right| \mathcal H^2_{\lfloor S(v_{n_k}) \cap [(R_0 \setminus \overline C) \times I]} \\
+\left( 1+ \left|\left( \nabla_\a u_{k} \Big| \frac{1}{\e_{n_k}}
\nabla_3 u_{k} \right) \right|^p \right) \mathcal L^3_{\lfloor (R_0
\setminus \overline C) \times I} +  \left| \left( \big(\nu_{u_{k}}\big)_\a
\Big|\frac{1}{\e_{n_k}} \big(\nu_{u_{k}}\big)_3 \right) \right| \mathcal
H^2_{\lfloor S(u_{k}) \cap [(R_0 \setminus \overline C) \times I]}
\end{eqnarray*}
is uniformly bounded, and thus, for a subsequence that will  not be
relabeled, there exists a positive Radon measure $\lambda$ such that
$\lambda _k \xrightharpoonup[]{*} \lambda$
in $\mathcal M_b(R_0 \setminus \overline C)$.\\

Let $t>0$, define $R_t:=\{ x_\a \in R_0 : {\rm dist}(x_\a,\partial
R_0)>t\}$ and for any $0< \delta < \eta$, $L_\delta:=R_{\eta -
\delta} \setminus \overline R_{\eta + \delta}$. Since we are
localizing in $\Rb^2$, we consider a cut-off function
$\varphi_\delta \in \C^\infty_c(R_{\eta - \delta};[0,1])$ depending
only on $x_\a$ and satisfying $\varphi_\delta = 1$ on $R_{\eta}$ and
$\|\varphi_\delta\|_{L^\infty(R_{\eta-\delta})} \leq C/\delta$.
Define
$$w_k(x):=u_k(x) \varphi_\delta(x_\a) + v_{n_k}(x)(1-\varphi_\delta(x_\a)).$$
Then, $w_k \in SBV^p(A \times I;\Rb^3)$, $w_k \to u$ in $L^p(A
\times I;\Rb^3)$ and in view of (\ref{sup1}) and (\ref{sup2})
$$\sup_{k \in \Nb}\|w_k\|_{L^\infty(A \times I;\Rb^3)} \leq \sup_{k \in \Nb} \|u_k\|_{L^\infty(R_0 \times I;\Rb^3)} + \sup_{k \in \Nb} \|v_{n_k}\|_{L^\infty([A \setminus \overline C] \times I;\Rb^3)} < +\infty.$$
From (\ref{IA}), (\ref{l1}) and (\ref{l2}), we deduce that
\begin{eqnarray*}
\I_\infty(u;A) & \leq & \liminf_{k \to +\infty} \I_{\e_{n_k}} (w_k;A)\\
&\leq  & \liminf_{k \to +\infty} \left\{ \I_{\e_{n_k}}(u_k;R_{\eta+\delta}) + \I_{\e_{n_k}}(v_{n_k};A \setminus \overline R_{\eta - \delta})\right.\\
 & & \hspace{1.5cm} \left. + C\lambda_k(L_\delta) + \frac{C}{\delta^p} \int_{L_\delta \times I}|u_k - v_{n_k}|^p\, dx.  \right\}\\
 & \leq & \I_\infty(u;R_0) + \I_\infty(u;A \setminus \overline C) + \eta + C\lambda(\overline L_\delta)\\
 & \leq & \I_\infty(u;B) + \I_\infty(u;A \setminus \overline C) + \eta + C\lambda(\overline L_\delta),
\end{eqnarray*}
where we have used the fact that $\I_\infty(u;\cdot)$ is an
increasing set function. Note that the previous computation would
not hold if we had considered $\I$ instead of $\I_\infty$ because
the minimizing sequences would only converge in $L^1$. Letting
$\delta$ tend to zero, we obtain,
$$\I_\infty(u;A) \leq  \I_\infty(u;B) + \I_\infty(u;A \setminus \overline C) + \eta + C\lambda(\partial R_\eta).$$
Now choose a sequence $\eta_h \to 0$ such that $\lambda(\partial R_{\eta_h})=0$. Letting $h \nearrow +\infty$ yields
$$\I_\infty(u;A) \leq  \I_\infty(u;B) + \I_\infty(u;A \setminus \overline C)$$
which completes the proof of item (i). Thus, according to
(\ref{harraps}) and Lemma 7.3 in \cite{BFF}, $\I_\infty(u;\cdot)$ is
the restriction to $\A(\o)$ of the Radon measure $\hat \mu$ which is
absolutely continuous with respect to $\mathcal L^2 + \mathcal
H^1_{\lfloor S(u)}$.
\hfill$\Box$\\

As a consequence of Lemma \ref{measure} and Lebesgue's Decomposition Theorem, there exists a $\mathcal L^2$-measurable function $h$ and a $\mathcal H^1_{\lfloor S(u)}$-measurable function $g$ such that for every $A \in \A(\o)$,
\begin{equation}\label{intrep}
\I_\infty(u;A)= \int_A h\, d\mathcal L^2 + \int_{A \cap S(u)} g\, d\mathcal H^1.
\end{equation}
We denote by $Q'(x_0,\rho)$ the open cube of $\Rb^2$ centered at $x_0 \in \o$ and of side length $\rho>0$, where $\rho$ is small enough so that $Q'(x_0,\rho) \in \A(\o)$. Since the measures $\mathcal L^2$ and $\mathcal H^1_{\lfloor S(u)}$ are mutually singular, $h$ is the Radon-Nikodym derivative of $\I_\infty(u;\cdot)$ with respect to $\mathcal L^2$,
$$h(x_0)=\lim_{\rho \to 0}\frac{\I_\infty(u;Q'(x_0,\rho))}{\rho^2}, \quad \text{for } \mathcal L^2 \text{-a.e. } x_0 \in \o$$
and $g$ is the Radon-Nikodym derivative of $\I_\infty(u;\cdot)$ with respect to $\mathcal H^1_{\lfloor S(u)}$,
$$g(x_0)=\lim_{\rho \to 0}\frac{\I_\infty(u;Q'(x_0,\rho))}{\mathcal H^1(S(u) \cap Q'(x_0,\rho))}, \quad \text{for } \mathcal H^1 \text{-a.e. } x_0 \in S(u).$$

Now we would like to identify both densities $g$ and $h$. Note that we cannot use classical Integral Representation Theorems in $SBV$ (see e.g. Theorem 2.4 in \cite{BCP} or Theorem 1 in \cite{BFLM}) because the term of surface energy does not grow linearly in  the deformation jump.


\subsection{Characterization of the $\Gamma$-limit}\label{caraterisation}

\subsubsection{The upper bound}\label{upper}

\noindent We will proceed in two steps to prove the upper bound.
Firstly, we will show that the inequality holds for deformations
belonging to $L^\infty(\O;\Rb^3)$ (see Lemma \ref{gs1} below).
Indeed, for those, we will use the integral representation proved
above. In fact, we will show, with the help of a blow up argument,
that the inequality holds separately for the surface and the bulk
terms. Then, we will prove the inequality in its full generality in
Lemma \ref{gs2}, using a truncation argument as in the proof of
Lemma \ref{bd}.

\begin{lemma}\label{gs1}
For all $u \in BV(\O;\Rb^3) \cap L^\infty(\O;\Rb^3)$, $\I_\infty(u;\o) \leq \J(u;\o)$.
\end{lemma}

\noindent {\it Proof. }It is enough to consider the case where $\J(u;\o)<+\infty$ and
thus $u \in SBV^p(\o;\Rb^3)$. Let $u \in L^\infty(\o;\Rb^3) \cap SBV^p(\o;\Rb^3)$,
according to (\ref{intrep}) and (\ref{JA}), we must show that
 $h(x_0) \leq 2 \mathcal Q W_0 (\nabla_\a u(x_0))$ for $\mathcal L^2$-a.e. $x_0 \in \o$ and $g(x_0) \leq 2$ for $\mathcal H^1$-a.e. $x_0 \in S(u)$.\\

Let us first treat the surface term. We have for $\mathcal H^1$-a.e. $x_0 \in S(u)$,
\begin{eqnarray*}
g(x_0) & = & \lim_{\rho \to 0}\frac{\I_\infty(u;Q'(x_0,\rho))}{\mathcal H^1(S(u) \cap Q'(x_0,\rho))}\\
 & \leq & \limsup_{\rho \to 0} \frac{1}{\mathcal H^1(S(u) \cap Q'(x_0,\rho))}\left[2\int_{Q'(x_0,\rho)} W(\nabla_\a u|0) dx_\a + 2 \mathcal H^1(S(u) \cap Q'(x_0,\rho)) \right]\\
 & =  & \limsup_{\rho \to 0} \frac{\mu(Q'(x_0,\rho))}{\mathcal H^1(S(u) \cap Q'(x_0,\rho))} + 2,
\end{eqnarray*}
where we set $\mu:= 2W(\nabla_\a u|0) \mathcal L^2$. But since $\mu$ and $\mathcal H^1_{\lfloor S(u)}$ are mutually singular, we have for $\mathcal H^1$-a.e. $x_0 \in S(u)$
$$\lim_{\rho \to 0}\frac{\mu(Q'(x_0,\rho))}{\mathcal H^1(S(u) \cap Q'(x_0,\rho))} = 0.$$
This shows that $g(x_0) \leq 2$ for $\mathcal H^1$-a.e. $x_0 \in S(u)$.\\

Concerning the bulk term, choose $x_0 \in \o$ such that
\begin{equation}\label{x0}
\lim_{\rho \to 0} \med_{Q'(x_0,\rho)}| \nabla_\a u(x_\a)-\nabla_\a
u(x_0)|^p\, dx =0.
\end{equation}
and
\begin{equation}\label{sing}
\lim_{\rho \to 0} \frac{\mathcal H^1(S(u) \cap Q'(x_0,\rho))}{\rho^2}=0.
\end{equation}
Since $\nabla_\a u \in L^p(\o;\Rb^{3 \times 2})$ and the measures
$\mathcal L^2$ and $\mathcal H^1_{\lfloor S(u)}$ are mutually
singular, it follows that $\mathcal L^2$-a.e. $x_0 \in \o$ satisfies
(\ref{x0}) and (\ref{sing}). For every $\rho>0$, Theorem 2 in
\cite{LDR} implies the existence of a sequence $\{v^\rho_n\} \subset
W^{1,p}(Q'(x_0,\rho) \times I ;\Rb^3)$ such that $v^\rho_n \to
\nabla_\a u(x_0) \cdot x_\a$ in $L^p(Q'(x_0,\rho) \times I;\Rb^3)$
(thus a fortiori in $L^1(Q'(x_0,\rho) \times I;\Rb^3)$) and
\begin{equation}\label{ledret}
\int_{Q'(x_0,\rho) \times I} W\left( \nabla_\a v^\rho_n
\Big|\frac{1}{\e_{n}} \nabla_3 v^\rho_n \right) dx \to 2 \rho^2 \mathcal
QW_0(\nabla_\a u(x_0)).
\end{equation}
Moreover, by the construction of this recovery sequence (see
\cite{LDR}), there is no loss of generality in assuming that
$\{v^\rho_n\}$  further satisfies
$$\sup_{\rho>0, \, n \in \Nb} \|v^\rho_n\|_{L^\infty(Q'(x_0,\rho) \times I;\Rb^3)}
<+\infty.$$ From the coercivity condition (\ref{pg1}), we get
\begin{equation}\label{coerc}
M:=\sup_{\rho>0,\, n \in \Nb}\med_{Q'(x_0,\rho) \times I} \left|
\left( \nabla_\a v^\rho_n \Big|\frac{1}{\e_{n}} \nabla_3 v^\rho_n
\right)\right|^p dx < +\infty.
\end{equation}
Define $u_n^\rho(x):= u(x_\a) + v^\rho_n(x_\a,x_3) - \nabla_\a
u(x_0) \cdot x_\a$. Then,
$$u_n^\rho \to u \text{ in } L^1(Q'(x_0,\rho) \times I;\Rb^3) \text{ as } n \to +\infty,
 \quad \sup_{\rho>0, \, n \in \Nb} \|u_n^\rho\|_{L^\infty(Q'(x_0,\rho) \times I;\Rb^3)} <+\infty$$
and $S(u_n^\rho) \cap [Q'(x_0,\rho) \times I] = [S(u) \cap
Q'(x_0,\rho)] \times I$. Thus,
\begin{eqnarray*}
\frac{\I_\infty(u;Q'(x_0,\rho))}{\rho^2} & \leq & \liminf_{n \to +\infty}\frac{1}{\rho^2}
\left\{ \int_{Q'(x_0,\rho) \times I} W\left( \nabla_\a u_n^\rho \Big|\frac{1}{\e_n} \nabla_3 u_n^\rho \right) dx \right.\\
&&\hspace{2.0cm} \left. + \int_{S(u_n^\rho) \cap [Q'(x_0,\rho) \times I]}
\left| \left( \left(\nu_{u_n^\rho}\right)_\a \Big|\frac{1}{\e_n} \left(\nu_{u_n^\rho}\right)_3 \right)\right| d\mathcal H^2 \right\}\\
 & \leq & \liminf_{n \to +\infty} \frac{1}{\rho^2}\int_{Q'(x_0,\rho) \times I}
 W\left( \nabla_\a u(x_\a) - \nabla_\a u(x_0) + \nabla_\a v^\rho_n (x) \Big|\frac{1}{\e_{n}} \nabla_3 v^\rho_n(x) \right) dx \\
&&\hspace{2.0cm} + 2 \frac{ \mathcal H^1(S(u) \cap Q'(x_0,\rho))}{\rho^2}.
\end{eqnarray*}
Thus from (\ref{sing}), we obtain
$$h(x_0) \leq \liminf_{\rho \to 0}  \liminf_{n \to +\infty} \frac{1}{\rho^2} \int_{Q'(x_0,\rho) \times I}
W\left( \nabla_\a u(x_\a) - \nabla_\a u(x_0) + \nabla_\a v^\rho_n
(x) \Big|\frac{1}{\e_{n}} \nabla_3 v^\rho_n(x) \right) dx.$$ Relations (\ref{pg1bis}), (\ref{ledret}), (\ref{coerc})
and H\"older's inequality yield
\begin{eqnarray*}
h(x_0) & \leq & \liminf_{\rho \to 0}  \liminf_{n \to +\infty}
\frac{1}{\rho^2} \Bigg\{  \int_{Q'(x_0,\rho) \times I}
W\left(  \nabla_\a v^\rho_n  \Big|\frac{1}{\e_n} \nabla_3 v^\rho_n \right) dx\\
&&\hspace{0.75cm} + C \int_{Q'(x_0,\rho) \times I}
\bigg(1+|\nabla_\a u(x_\a) - \nabla_\a u(x_0)|^{p-1} \\
&&\hspace{1.5cm}   + \Big| \Big(  \nabla_\a v^\rho_n(x)
\Big|\frac{1}{\e_{n}} \nabla_3 v^\rho_n(x) \Big)\Big|^{p-1} \bigg)
|\nabla_\a u(x_\a) - \nabla_\a u(x_0)|\, dx \Bigg\}\\
 & \leq & 2\mathcal Q W_0(\nabla_\a u(x_0)) + C\,  \limsup_{\rho \to 0}  \left( \med_{Q'(x_0,\rho)}
 |\nabla_\a u(x_\a) - \nabla_\a u(x_0)|^p\, dx_\a \right) \\
&&\hspace{1.0cm} +C \big(1+M^{(p-1)/p}\, \big)\,  \limsup_{\rho \to
0}\left(\med_{Q'(x_0,\rho)}
 |\nabla_\a u(x_\a) - \nabla_\a u(x_0)|^p\, dx_\a \right)^{1/p}.
\end{eqnarray*}
Thanks to (\ref{x0}), we conclude  that $h(x_0) \leq 2\mathcal Q
W_0(\nabla_\a u(x_0))$ for $\mathcal L^2$-a.e. $x_0 \in \o$.
\hfill$\Box$\\

Let us now turn back to the general case.

\begin{lemma}\label{gs2}
For all $u \in BV(\O;\Rb^3)$, $\I(u;\o) \leq \J(u;\o)$.
\end{lemma}

\noindent{\it Proof. }As in the proof of Lemma \ref{gs1}, we can
assume without loss of generality that $\J(u;\o) < +\infty$ and thus
that  $u \in SBV^p(\o;\Rb^3)$. In particular, it implies that
$\I(u;\o) < +\infty$. Let $\varphi_i \in \C^1_c(\Rb^3;\Rb^3)$ be the
truncation function introduced in Lemma \ref{bd} and defined by
(\ref{fii}). The Chain Rule formula, Theorem 3.96 in \cite{AFP},
implies that $\varphi_i(u) \in SBV^p(\o;\Rb^3) \cap
L^\infty(\o;\Rb^3)$ and
\begin{equation}\label{ui}
\left\{\begin{array}{l}
\|\varphi_i(u)\|_{L^\infty(\o;\Rb^3)} \leq e^i,\\
\\
S(\varphi_i(u))  \subset S(u) ,\\
\\
\nabla_\a \big( \varphi_i(u(x_\a)) \big) = \nabla \varphi_i(u(x_\a)) \circ \nabla_\a u(x_\a) \quad \mathcal L^2\text{-a.e. in }\o
\end{array}\right.
\end{equation}
and $\varphi_i(u) \to u$ in $L^1(\o;\Rb^3)$ as $i \to +\infty$.
Since $\varphi_i(u) \in L^\infty(\o;\Rb^3) \cap SBV^p(\o;\Rb^3)$, it
follows from Lemmas \ref{bd} and \ref{gs1} that
$$\I(\varphi_i(u);\o) = \I_\infty(\varphi_i(u);\o) \leq \J(\varphi_i(u);\o)$$
and by lower semicontinuity of $\I(\cdot;\o)$ with respect to the strong $L^1(\o;\Rb^3)$-convergence, we have
\begin{equation}\label{1735}
\I(u;\o) \leq \liminf_{i \to +\infty}\I(\varphi_i(u);\o) \leq \limsup_{i \to +\infty}\J(\varphi_i(u);\o).
\end{equation}
But, in view of (\ref{ui}), $\mathcal H^1(S(\varphi_i(u))) \leq \mathcal H^1(S(u))$ and, thanks to (\ref{fii}),
$$\int_\o \mathcal Q W_0(\nabla_\a (\varphi_i(u)))\, dx_\a \leq \int_{\{|u| < e^i\}} \mathcal Q W_0(\nabla_\a u)\, dx_\a + \beta\int_{\{|u| \geq e^i\}}(1+|\nabla_\a u|^p)\, dx_\a.$$
Thus,
\begin{eqnarray}\label{1736}
&&\limsup_{i \to +\infty} \J(\varphi_i(u);\o)\nonumber \\
&&\hspace{1.0cm} \leq \limsup_{i \to +\infty} \left\{ 2\int_\o \mathcal Q W_0(\nabla_\a u)\, dx_\a +2 \beta \int_{\{|u| \geq e^i\}}(1+|\nabla_\a u|^p)\, dx_\a +2 \mathcal H^1(S(u)) \right\} \nonumber\\
&&\hspace{1.0cm} \leq  2\int_\o \mathcal Q W_0(\nabla_\a u)\, dx_\a  + 2\mathcal H^1(S(u))\nonumber\\
&&\hspace{1.0cm} = \J(u;\o),
\end{eqnarray}
where we have used the fact that, by Chebyshev's inequality, $\mathcal L^2(\{|u| \geq e^i\}) \leq \|u\|_{L^1(\o;\Rb^3)}/e^i \to 0$ as $i \to +\infty$. Gathering (\ref{1735}) and (\ref{1736}), we deduce that $\I(u;\o) \leq \J(u;\o)$ and this completes the proof of the Lemma.
\hfill$\Box$

\subsubsection{The lower bound}\label{lower}

\noindent Let us now prove the lower bound. The proof is essentially
based on a lower semicontinuity result in $SBV^p(\O;\Rb^3)$. The
main difficulty remains to show that any deformation $u \in
BV(\O;\Rb^3)$ satisfying $\I(u;\o)<+\infty$ belongs in fact to
$SBV^p(\o;\Rb^3)$.

\begin{lemma}\label{gi}
For all $u \in BV(\O;\Rb^3)$, $\I(u;\o) \geq \J(u;\o)$.
\end{lemma}

\noindent{\it Proof. }It is not restrictive to assume that $\I(u;\o)<+\infty$. By $\G$-convergence, there exists a sequence $\{u_n\} \subset SBV^p(\O;\Rb^3)$ such that $u_n \to u$ in $L^1(\O;\Rb^3)$ and
\begin{equation}\label{1805}
\lim_{n \to +\infty}\left[ \int_{\O} W\left( \nabla_\a u_n \Big|\frac{1}{\e_n}\nabla_3 u_n \right)dx + \int_{S(u_n)}\left| \left( \big( \nu_{u_n} \big)_\a \Big|\frac{1}{\e_n} \big( \nu_{u_n} \big)_3 \right) \right|d\mathcal H^2 \right] = \I(u;\o).
\end{equation}
Let us show that $u \in SBV(\o;\Rb^3)$. We will use the truncation
function $\varphi_i \in \C^1_c(\Rb^3;\Rb^3)$  defined in
(\ref{fii}). The Chain Rule formula, Theorem 3.96 in \cite{AFP},
implies that $\varphi_i(u_n) \in SBV^p(\O;\Rb^3)$ and

$$\left\{\begin{array}{l}
\|\varphi_i(u_n)\|_{L^\infty(\O;\Rb^3)} \leq e^i,\\
\\
S(\varphi_i(u_n))  \subset S(u_{n}),\\
\\
\nabla \big( \varphi_i(u_n)(x) \big)= \nabla \varphi_i(u_n(x)) \circ \nabla u_n(x) \quad \mathcal L^3\text{-a.e. in }\O.
\end{array}\right.$$
As $\nu_{u_n}(x)=\pm \nu_{\varphi_i(u_n)}(x)$ for $\mathcal H^2$-a.e. $x \in S(\varphi_i(u_n))$, we get
\begin{eqnarray*}
\sup_{n \in \Nb} \left[ \int_{\O} \left|\left( \nabla_\a \big(\varphi_i(u_n)\big) \Big|\frac{1}{\e_n} \nabla_3 \big(\varphi_i(u_n)\big) \right)\right|^p dx + \int_{S(\varphi_i(u_n))}\left| \left( \big( \nu_{\varphi_i(u_n)} \big)_\a \Big|\frac{1}{\e_n} \big( \nu_{\varphi_i(u_n)} \big)_3 \right) \right| d\mathcal H^2 \right]\\
\leq \sup_{n \in \Nb} \left[ \int_{\O} \left|\left( \nabla_\a u_n \Big|\frac{1}{\e_n}\nabla_3 u_n \right)\right|^p dx + \int_{S(u_n)}\left| \left( \big( \nu_{u_n} \big)_\a \Big|\frac{1}{\e_n} \big( \nu_{u_n} \big)_3 \right) \right|d\mathcal H^2 \right]
<+\infty,
\end{eqnarray*}
where we used (\ref{1805}) together with the coercivity condition
(\ref{pg1}). Lemma \ref{comp} and a diagonalization argument yield
the existence of a subsequence (still denoted by $\{\e_n\}$), and a
function $v_i \in SBV^p(\o;\Rb^3)$ such that $\varphi_i(u_n)
\rightharpoonup v_i$ in $SBV^p(\O;\Rb^3)$ as $n \to +\infty$. But,
since $u_n \to u$ and $\varphi_i(u_n) \to \varphi_i(u)$ in
$L^1(\O;\Rb^3)$ as $n \to +\infty$, we deduce that $v_i=\varphi_i(u)
\in SBV^p(\o;\Rb^3)$ for every $i \in \Nb$. By virtue of Theorem
3.96 in \cite{AFP},
$$0=D^c v_i =\nabla \varphi_i(\tilde u) \circ D^c u \quad \text{ in } \o \setminus S(u),$$
where $\tilde u$ denotes the approximate limit of $u$ at $x_\a \in
\o \setminus S(u)$. Define $E_i:= \{x_\a \in \o \setminus S(u) :
|\tilde u(x_\a)| < e^i \}$, since $\tilde u$ is a Borel function and
$S(u)$ is a Borel set (see Proposition 3.64 (a) in \cite{AFP}),
$E_i$ is a Borel set. Moreover, as $\{E_i\}$ is an increasing
sequence of sets whose union is $\o \setminus S(u)$, we get
$$|D^c u|(\o)= |D^c u|(\o \setminus S(u)) = \lim_{i \to +\infty}|D^c u|(E_i)= \lim_{i \to +\infty} |\nabla \varphi_i(\tilde u) \circ D^c u|(E_i)=\lim_{i \to +\infty} |D^c v_i|(E_i)=0,$$
where we have used the fact that $\nabla \varphi_i(\tilde
u(x_\a))={\rm Id}$ for all $x_\a \in E_i$. Thus $u \in
SBV(\o;\Rb^3)$ and by (\ref{1805}), Remark \ref{qcvx}, Theorem 5.29
in \cite{AFP} and Theorem 3.7 in \cite{A}
\begin{eqnarray*}
\I(u;\o) & \geq & \liminf_{n \to +\infty}\left[ \int_{\O} \mathcal QW_0 ( \nabla_\a u_n )\, dx + \mathcal H^2( S(u_n) ) \right]\\
 & \geq & 2\int_\o \mathcal QW_0 ( \nabla_\a u )\, dx_\a + 2 \mathcal H^1( S(u) ).
\end{eqnarray*}
 In particular, the $p$-coercivity of $\mathcal QW_0$ implies
that $u \in SBV^p(\o;\Rb^3)$ and thus, according to (\ref{JA}), that
$\I(u;\o)\geq \J(u,\o)$.
\hfill $\Box$\\

\noindent{\it Proof of Theorem \ref{gammaconv}. }We have shown that
for any sequence $\{\e_j\} \searrow 0^+$, there exists a further
subsequence $\{\e_{j_n}\}\equiv\{\e_n\}$ such that
$\I_{\e_n}(\cdot;\o)$ $\G$-converges to $\I(\cdot;\o)$ for the
strong $L^1(\O;\Rb^3)$-topology. By virtue of Lemmas \ref{gs2} and
\ref{gi}, we have $\I(\cdot;\o)=\J(\cdot;\o)$. Since the $\G$-limit
does not depend upon the extracted subsequence, we deduce, in the
light of Proposition 8.3 in \cite{DM}, that the whole sequence
$\I_\e(\cdot;\o)$ $\G$-converges to $\J(\cdot;\o)$. \hfill$\Box$


\subsection{Boundary conditions}\label{boundconditions}

\noindent Let us now deal with boundary condition constraints that
will be of use in Lemmas \ref{MIN0} and \ref{MIN} in order to prove
the minimality property of the limit quasistatic evolution. Indeed,
it will allow to extend functions on the enlarged cylinder $\O'$ by
the value of the boundary condition. The following result, very
close in spirit to Lemma 2.6 in \cite{BFF}, relies on De Giorgi's
slicing argument together with the fact that we can consider cut-off
functions depending only on $x_\a$ (see also the proof of Lemma
\ref{measure}). It is established that any recovery sequence can be
chosen so as to match the lateral boundary condition of its target.

We recall that $\o'$ is a bounded open subset of $\Rb^2$ containing
$\o$ and that $\O'= \o' \times I$. In all that follows, if $v \in
SBV^p(\O';\Rb^3)$, we will denote by $v^-$ (resp. $v^+$) the inner
(resp. outer) trace of $v$ on $\partial \o \times I$.

\begin{lemma}\label{lbc}
For every $u \in SBV^p(\o;\Rb^3) \cap  L^\infty(\o;\Rb^3)$, there
exists a sequence $\{\bar u_\e\} \subset SBV^p(\O;\Rb^3)$ such that
$\bar u_\e \to u$ in $L^p(\O;\Rb^3)$, $\bar u_\e = u$ in a
neighborhood of $\partial \o \times I$ and
$$\J(u) = \lim_{\e \to 0} \left[\int_\O W\left(\nabla_\a \bar u_\e \Big|\frac{1}{\e} \nabla_3 \bar u_\e \right)dx
+\int_{S(\bar u_\e)}\left|\left( \big(\nu_{\bar u_\e}\big)_\a \Big|\frac{1}{\e} \big( \nu_{\bar u_\e}\big)_3 \right) \right|
 d\mathcal H^2 \right].$$
\end{lemma}

\noindent {\it Proof. } According to Theorem  \ref{gammaconv} and
Lemma \ref{bd}, there exists a sequence $\{u_\e\} \subset SBV^p(\O;\Rb^3)$ strongly converging
to $u$ in $L^1(\O;\Rb^3)$, satisfying $\sup_{\e > 0}\|u_\e\|_{L^\infty(\O;\Rb^3)} < +\infty$ and
$$\J(u) = \lim_{\e \to 0} \left[\int_\O W\left(\nabla_\a u_\e \Big|\frac{1}{\e} \nabla_3 u_\e \right)dx
+\int_{S(u_\e)}\left|\left( \big(\nu_{u_\e}\big)_\a \Big|\frac{1}{\e} \big( \nu_{u_\e}\big)_3 \right) \right|
 d\mathcal H^2 \right].$$
In particular, $u_\e \to u$ in $L^p(\O;\Rb^3)$ and from the $p$-coercivity condition (\ref{pg1}), it follows that
\begin{equation}\label{1404}
C:=\sup_{\e > 0} \left[ \int_{\O} \left( 1+\left|\left( \nabla_\a u_\e \Big|\frac{1}{\e}\nabla_3 u_\e \right)\right|^p\right) dx
 + \int_{S(u_\e)}\left| \left( \big( \nu_{u_\e} \big)_\a \Big|\frac{1}{\e} \big( \nu_{u_\e} \big)_3 \right) \right|d\mathcal H^2 \right]
<+\infty.
\end{equation}
Set $$K_\e:= \left \llbracket
\frac{1}{\| u_\e -u \|_{L^p(\O;\Rb^3)}^{1/2}} \right \rrbracket,
 \quad M_\e:=\left \llbracket \sqrt{K_\e} \right
\rrbracket$$ and denote
$$\o(\e):= \left\{ x_\a \in \o: \; {\rm dist}(x_\a,\partial \o) <
\frac{M_\e}{K_\e}\right\} \quad \text{ and } \quad \o_i^\e:= \left\{
x_\a \in \o: \; {\rm dist}(x_\a,\partial \o) \in
\Bigg[\frac{i}{K_\e},\frac{i+1}{K_\e} \Bigg)\right\},$$ for all $i
\in \{0,\ldots ,M_\e-1\}$. From (\ref{1404}), we get the existence
of a $i(\e) \in \{0,\ldots ,M_\e-1\}$ such that
\begin{equation}\label{1405}
\int_{\o^\e_{i(\e)} \times I} \left( 1+ \left| \left(\nabla_\a u_\e
\Big|\frac{1}{\e} \nabla_3 u_\e\right)\right|^p \right) dx +
\int_{S(u_\e) \cap [\o^\e_{i(\e)} \times I]} \left| \left(
\big(\nu_{u_\e} \big)_\a \Big| \frac{1}{\e} \big(\nu_{u_\e} \big)_3
\right)\right| d\mathcal H^2 \leq \frac{C}{M_\e}.
\end{equation}
Let us now consider a cut-off function $\phi_\e \in
\C^\infty_c(\o;[0,1])$ independent of $x_3$ and satisfying
$$\phi_\e(x_\a)=\left\{
\begin{array}{rcl}
1 & \text{if} & \ds {\rm dist}(x_\a,\partial \o) >
\frac{i(\e)+1}{K_\e},\\
0 & \text{if} & \ds {\rm dist}(x_\a,\partial \o) \leq
\frac{i(\e)}{K_\e}
\end{array}
\right.  \quad \text{ and } \quad \|\nabla_\a \phi_\e
\|_{L^\infty(\o)} \leq 2 K_\e.$$ Define $\bar u_\e(x):=\phi_\e(x_\a)
u_\e (x) + (1-\phi_\e(x_\a)) u(x_\a)$; then $\bar u_\e \in
SBV^p(\O;\Rb^3)$, $\bar u_\e \to u$ in $L^p(\O;\Rb^3)$, $\bar u_\e =
u$ in a neighborhood of $\partial \o \times I$ and $S(\bar u_\e)
\subset S(u_\e) \cup \big( S(u) \times I \big)$. The $p$-growth
condition (\ref{pg1}) implies that
\begin{eqnarray*}
\J(u) & \geq & \limsup_{\e \to 0} \Bigg[\int_{\left\{x_\a \in \o : \,
{\rm dist}(x_\a,\partial \o) >\frac{i(\e)+1}{K_\e} \right\} \times
I} W \left(\nabla_\a \bar u_\e \Big|\frac{1}{\e} \nabla_3 \bar
u_\e\right) dx \\
&&\hspace{2cm}+ \int_{S(\bar u_\e) \cap \left[ \left\{x_\a \in \o :
\, {\rm dist}(x_\a,\partial \o) >\frac{i(\e)+1}{K_\e} \right\}
\times I \right]} \left| \left( \big(\nu_{\bar u_\e} \big)_\a
\Big|\frac{1}{\e}
\big(\nu_{\bar u_\e} \big)_3 \right)\right| d\mathcal H^2 \Bigg]\\
& \geq & \limsup_{\e \to 0} \Bigg[\int_\O W \left(\nabla_\a
\bar u_\e \Big|\frac{1}{\e} \nabla_3 \bar u_\e\right) dx +
\int_{S(\bar u_\e)} \left| \left( \big(\nu_{\bar u_\e} \big)_\a
\Big|\frac{1}{\e} \big(\nu_{\bar u_\e} \big)_3 \right)\right|
d\mathcal H^2\\
&& -2\beta \int_{\{x_\a \in \o : \, {\rm dist}(x_\a,\partial \o)
\leq \frac{i(\e)}{K_\e} \}}(1+|\nabla_\a u|^p) \, dx_\a \\
&&\hspace{5cm} -2 \mathcal H^1 \left (S(u)
 \cap \left\{ x_\a \in \o : \, {\rm
dist}(x_\a,\partial \o) \leq \frac{i(\e)}{K_\e} \right\} \right)\\
&&\hspace{1cm} -c_1 K_\e ^p \int_{\o_{i(\e)}^\e \times I} |u_\e -
u|^p \, dx - c_2\int_{\o_{i(\e)}^\e \times I}\left( 1 + |\nabla_\a
u|^p + \left| \left(\nabla_\a u_\e \Big|\frac{1}{\e} \nabla_3
u_\e\right)\right|^p\right)dx\\
&&\hspace{2cm} -c_3 \int_{S(u_\e) \cap [\o_{i(\e)}^\e \times I]}
\left| \left( \big(\nu_{u_\e} \big)_\a \Big| \frac{1}{\e}
\big(\nu_{u_\e} \big)_3 \right)\right| d\mathcal H^2 -c_4 \mathcal
H^1 \big(S(u) \cap \o_{i(\e)}^\e \big) \Bigg]\\
& \geq & \limsup_{\e \to 0} \Bigg[\int_\O W \left(\nabla_\a
\bar u_\e \Big|\frac{1}{\e} \nabla_3 \bar u_\e\right) dx +
\int_{S(\bar u_\e)} \left| \left( \big(\nu_{\bar u_\e} \big)_\a
\Big|\frac{1}{\e} \big(\nu_{\bar u_\e} \big)_3 \right)\right|
d\mathcal H^2\\
&& \hspace{2cm} -c \Bigg( \int_{\o(\e) }(1+|\nabla_\a
u|^p) \, dx_\a + \mathcal H^1 (S(u) \cap \o(\e) ) + \| u_\e
-u\|^{p/2}_{L^p(\O;\Rb^3)} + \frac{1}{M_\e} \Bigg) \Bigg],
\end{eqnarray*}
where we have used (\ref{1405}) in the last inequality. Thus, since
$M_\e \to +\infty$, $M_\e/K_\e \to 0$ and $\mathcal H^1(S(u)) <
+\infty$, we obtain from the previous relation and the $\G$-$\liminf$
inequality
$$\J(u) = \lim_{\e \to 0} \Bigg[\int_\O W \left(\nabla_\a
\bar u_\e \Big|\frac{1}{\e} \nabla_3 \bar u_\e\right) dx +
\int_{S(\bar u_\e)} \left| \left( \big(\nu_{\bar u_\e} \big)_\a
\Big|\frac{1}{\e} \big(\nu_{\bar u_\e} \big)_3 \right)\right|
d\mathcal H^2\Bigg].$$
\hfill$\Box$\\

Let us now state a $\G$-convergence result involving the boundary
conditions. Consider a sequence of boundary conditions $\{g_\e\}
\subset W^{1,p}(\O';\Rb^3)$, and let $g \in W^{1,p}(\o';\Rb^3)$ and
$H \in L^p(\O';\Rb^3)$ be such that
\begin{equation}\label{flyaway}
\left\{ \begin{array}{l}
\ds \sup_{\e >0} \|g_\e\|_{L^\infty(\O';\Rb^3)} < +\infty,\\
\ds g_\e \to g \text{  in } W^{1,p}(\O';\Rb^3), \quad
\frac{1}{\e}\nabla_3 g_\e \to H  \text{ in }L^p(\O';\Rb^3).
\end{array} \right.
\end{equation}
Then,

\begin{coro}\label{corol}
The functional $\I^{g_\e}_\e : BV(\O';\Rb^3) \to \overline \Rb$ defined by
$$\I^{g_\e}_\e(u):=\left\{
\begin{array}{cl}
\ds \int_{\O} W\left(\nabla_\a u\Big|\frac{1}{\e}\nabla_3 u\right)dx
+ \int_{S(u)} \left| \left( \left(\nu_u\right)_\a \Big|\frac{1}{\e}
\left( \nu_u \right)_3 \right)\right| d\mathcal
H^2 & \text{if } \left\{
\begin{array}{l}
u \in SBV^p(\O';\Rb^3),\\
u=g_\e \text{ on }[\o' \setminus \overline \o] \times I,
\end{array}
\right.\\
&\\
\ds +\infty & \text{otherwise}
\end{array}
\right.$$
$\G$-converges for the strong
$L^1(\O';\Rb^3)$-topology towards $\J^g : BV(\O';\Rb^3) \to \overline
\Rb$ defined by
$$\J^g(u):=\left\{
\begin{array}{cl}
\ds 2 \int_\o \mathcal Q W_0( \nabla_\a u)\, dx_\a +2 \mathcal H^1(S(u)) & \text{if }  \left\{
\begin{array}{l}
u \in SBV^p(\o';\Rb^3),\\
u=g \text{ on }\o' \setminus \overline \o,
\end{array}
\right.\\
&\\
\ds  +\infty & \text{otherwise}.
\end{array}
\right.$$
\end{coro}

\begin{rmk}{\rm
Note that in the statement of the previous Corollary, the bulk
integrals are still computed over $\O$ (resp. $\o$) as in Theorem
\ref{gammaconv}, however, since the jump set of the deformations can
now reach the lateral boundary $\partial \o \times I$, the surface
integrals are implicitly computed over $\overline \o \times I$
(resp. $\overline \o$) or equivalently $\O'$ (resp. $\o'$).}
\end{rmk}

\noindent {\it Proof. }Let us first prove the $\G$- $\liminf$ inequality. Consider
a sequence $\{u_\e\} \subset L^1(\O';\Rb^3)$ strongly converging to $u$
in $L^1(\O';\Rb^3)$. It is not restrictive to assume that
$$\liminf_{\e \to 0}\I^{g_\e}_\e(u_\e) < +\infty.$$
Then, for a (not relabeled) subsequence, $u_\e \in
SBV^p(\O';\Rb^3)$, $u_\e = g_\e$ on $[\o' \setminus \overline \o]
\times I$ and arguing as in the proof of Lemma \ref{gi}, we get that
$u \in SBV^p(\o';\Rb^3)$. Consequently, since $u=g$ on $\o'
\setminus \overline \o$, we get from Theorem \ref{gammaconv} and the
definition of $\J^g$ that
$$\J^g(u) \leq \liminf_{\e \to 0}\I^{g_\e}_\e(u_\e).$$

Let $u \in SBV^p(\o';\Rb^3)$ satisfying $u=g$ on $\o' \setminus
\overline \o$. It remains to construct a recovery sequence. We first
assume that $u \in L^\infty(\o';\Rb^3)$. Then, by virtue of Lemma
\ref{lbc}, there exists a sequence $\{\bar u_\e\} \subset
SBV^p(\O;\Rb^3)$ satisfying $\bar u_\e \to u$ in $L^1(\O;\Rb^3)$,
$\bar u_\e=u$ in a neighborhood of $\partial \o \times I$ and
\begin{eqnarray}\label{trux}
&&2\int_\o \mathcal QW_0(\nabla_\a u)\, dx_\a + 2 \mathcal H^1(S(u) \cap \o)\nonumber\\
&&\hspace{2cm}= \lim_{\e \to 0} \Bigg[\int_\O W \left(\nabla_\a
\bar u_\e \Big|\frac{1}{\e} \nabla_3 \bar u_\e\right) dx +
\int_{S(\bar u_\e) \cap \O} \left| \left( \big(\nu_{\bar u_\e} \big)_\a
\Big|\frac{1}{\e} \big(\nu_{\bar u_\e} \big)_3 \right)\right|
d\mathcal H^2\Bigg].
\end{eqnarray}
Since $g \in W^{1,p}(\o' \setminus \overline \o;\Rb^3)$, by Corollary 3.89 in \cite{AFP},
the function
$$v_\e:= \bar u_\e \chi_{\O} + g \chi_{[\o' \setminus \overline \o] \times I} \in BV(\O';\Rb^3)$$
and, viewing $D \bar u_\e$ (resp. $Dg$) as measures on all $\Rb^3$
and concentrated on $\O$ (resp. $[\o' \setminus \overline \o] \times
I$), we get as $v_\e^-=\bar u_\e ^-=u^-$ and $v_\e ^+=g^+$ on
$\partial \o \times I$
$$D v_\e = D \bar u_\e + (g^+ - \bar u^-) \otimes \nu_{\partial \o \times I} + D g.$$
In particular, we observe that $v_\e \in SBV^p(\O';\Rb^3)$ and $v_\e
\to u$ in $L^1(\O';\Rb^3)$ but we may have created some additional
jump set on $\partial \o \times I$. However, $S(v_\e) \cap [\partial
\o \times I] \; \widetilde = \; [S(u) \cap
\partial \o] \times I$, and since $\nu_{\partial \o \times I} =
\nu_{v_\e}$ $\mathcal H^2$-a.e. in $S(v_\e) \cap [\partial \o \times
I]$ and $\big(\nu_{\partial \o \times I} \big)_3=0$,
$$\int_{S(v_\e) \cap [\partial \o \times I]} \left| \left( \big(\nu_{v_\e} \big)_\a
\Big|\frac{1}{\e} \big(\nu_{v_\e} \big)_3 \right)\right| d\mathcal
H^2= \mathcal H^2\big( S(v_\e) \cap [\partial \o \times I] \big) =2
\mathcal H^1 \big( S(u) \cap \partial \o \big).$$ Replacing in
(\ref{trux}), it yields
\begin{eqnarray}\label{killers1}
\J^g(u) & = & 2\int_\o \mathcal QW_0(\nabla_\a u)\, dx_\a + 2 \mathcal H^1(S(u))\nonumber\\
& = & \lim_{\e \to 0} \Bigg[\int_\O W \left(\nabla_\a
v_\e \Big|\frac{1}{\e} \nabla_3 v_\e\right) dx +
\int_{S(v_\e)} \left| \left( \big(\nu_{v_\e} \big)_\a
\Big|\frac{1}{\e} \big(\nu_{v_\e} \big)_3 \right)\right|
d\mathcal H^2\Bigg].
\end{eqnarray}
Actually, $v_\e=g$ on $[\o' \setminus \overline \o] \times I$ so
that we need to modify $v_\e$ in order it to have the value $g_\e$
instead of $g$ on the enlarged part of the domain. Let $H_j \in
\C^\infty_c(\O;\Rb^3)$ be a sequence strongly converging to $H$ in
$L^p(\O;\Rb^3)$ and extended by the value zero on $[\o' \setminus
\overline \o] \times I$, and set
$b_j(x_\a,x_3):=\int_{-1}^{x_3}H_j(x_\a,s)\, ds$. We now define
$$u_\e^j(x):= v_\e(x) - g(x_\a) + g_\e(x) - \e b_j(x).$$
It follows that $u_\e^j \in SBV^p(\O';\Rb^3)$, $u_\e^j \to u$ in
$L^1(\O';\Rb^3)$ as $\e \to 0$, $u_\e^j= g_\e$ on $[\o' \setminus
\overline \o] \times I$. Furthermore, since $u_\e^j$ is a smooth
perturbation of $v_\e$ on the whole domain $\O'$, both sequences
have the same jump set, namely $S(u_\e^j) \; \widetilde = \;
S(u_\e)$, and consequently, the surface energy is not affected,
\begin{equation}\label{killers2}
\int_{S(u^j_\e)} \left| \left( \big(\nu_{u^j_\e} \big)_\a
\Big|\frac{1}{\e} \big(\nu_{u^j_\e} \big)_3 \right)\right|
d\mathcal H^2=\int_{S(v_\e)} \left| \left( \big(\nu_{v_\e} \big)_\a
\Big|\frac{1}{\e} \big(\nu_{v_\e} \big)_3 \right)\right|
d\mathcal H^2.
\end{equation}
Let us treat now the bulk energy. According to Remark \ref{WQW}, (\ref{pg1bis}) and H\"older's
Inequality, we have
\begin{eqnarray*}
&&\hspace{-0.7cm}\int_\O W \left(\nabla_\a u^j_\e \Big|\frac{1}{\e} \nabla_3 u^j_\e\right) dx\\
&&\hspace{-0.5cm} = \int_\O W \left(\nabla_\a v_\e -\nabla_\a g + \nabla_\a g_\e
-\e \nabla_\a b_j \Big|\frac{1}{\e} \nabla_3 v_\e + \frac{1}{\e} \nabla_3 g_\e -H_j\right) dx\\
&& \hspace{-0.5cm}\leq  \int_\O W \left(\nabla_\a v_\e \Big|\frac{1}{\e} \nabla_3 v_\e\right) dx +
\beta \int_\O \Bigg(1 + \left|\left(\nabla_\a v_\e \Big|\frac{1}{\e} \nabla_3 v_\e \right) \right|^{p-1}
+ \left|\left(\nabla_\a g_\e \Big|\frac{1}{\e} \nabla_3 g_\e \right) \right|^{p-1}\\
&& + |\nabla_\a g|^{p-1} + |(\e \nabla_\a b_j|H_j)|^{p-1} \Bigg) \left|\left(\nabla_\a g_\e-\nabla_\a g -\e
\nabla_\a b_j \Big|
\frac{1}{\e} \nabla_3 g_\e -H_j \right) \right| \, dx\\
&&\hspace{-0.5cm} \leq \int_\O W \left(\nabla_\a v_\e \Big|\frac{1}{\e} \nabla_3 v_\e\right) dx +
c\Bigg( 1 + \left\|\left(\nabla_\a v_\e \Big| \frac{1}{\e} \nabla_3 v_\e \right) \right\|_{L^p(\O;\Rb^{3\times 3})}^{p-1}
+  \left\|\left(\nabla_\a g_\e \Big| \frac{1}{\e} \nabla_3 g_\e \right) \right\|_{L^p(\O;\Rb^{3\times 3})}^{p-1}\\
&&  + \|\nabla_\a g\|^{p-1}_{L^p(\O;\Rb^{3\times 2})} + \|(\e \nabla_\a
b_j|H_j)\|_{L^p(\O;\Rb^{3\times 3})}^{p-1}\Bigg)
\left\|\left(\nabla_\a g_\e-\nabla_\a g -\e \nabla_\a b_j \Big|
\frac{1}{\e} \nabla_3 g_\e -H_j \right)
\right\|_{L^p(\O;\Rb^{3\times 3})}.
\end{eqnarray*}
Passing to the limit when $\e \to 0$, (\ref{flyaway}) yields
\begin{equation}\label{killers3}
\limsup_{\e \to 0}\int_\O W \left(\nabla_\a u^j_\e \Big|\frac{1}{\e} \nabla_3 u^j_\e\right) dx
\leq \limsup_{\e \to 0}\int_\O W \left(\nabla_\a v_\e \Big|\frac{1}{\e} \nabla_3 v_\e\right) dx
+ c'\|H -H_j\|_{L^p(\O;\Rb^3)}.
\end{equation}
Gathering (\ref{killers1}), (\ref{killers2}), (\ref{killers3}) and remembering that
$H_j \to H$ in $L^p(\O;\Rb^3)$, we get that
\begin{equation}\label{killers4}
\J^g(u)= \lim_{j \to +\infty}\lim_{\e \to 0} \Bigg[\int_\O W \left(\nabla_\a
u^j_\e \Big|\frac{1}{\e} \nabla_3 u^j_\e\right) dx +
\int_{S(u^j_\e)} \left| \left( \big(\nu_{u^j_\e} \big)_\a
\Big|\frac{1}{\e} \big(\nu_{u^j_\e} \big)_3 \right)\right|
d\mathcal H^2\Bigg]
\end{equation}
where we have also used the $\G$-$\liminf$ inequality. A standard diagonalization procedure
(see e.g. Lemma 7.1 in \cite{BFF}) implies the existence of a sequence $j_\e \nearrow
+\infty$ as $\e \to +\infty$ such that $u_\e:=u_\e^{j_\e} \in SBV^p(\O';\Rb^3)$, $u_\e \to u$ in $L^1(\O';\Rb^3)$,
$u_\e = g_\e$ on $[\o' \setminus \overline \o] \times I$ and
$$\J^g(u)= \lim_{\e \to 0} \Bigg[\int_\O W \left(\nabla_\a
u_\e \Big|\frac{1}{\e} \nabla_3 u_\e\right) dx +
\int_{S(u_\e)} \left| \left( \big(\nu_{u_\e} \big)_\a
\Big|\frac{1}{\e} \big(\nu_{u_\e} \big)_3 \right)\right|
d\mathcal H^2\Bigg].$$

If $u$ does not belong to $L^\infty(\o';\Rb^3)$, we can consider
$\varphi_i(u) \in SBV^p(\o';\Rb^3) \cap L^\infty(\o';\Rb^3)$ where
$\varphi_i \in \C_c^1(\Rb^3;\Rb^3)$ is the truncation function
defined in (\ref{fii}) and $i$ is large enough (independently of
$\e$) so that $e^i > \| g \|_{L^\infty(\o';\Rb^3)}$. In particular,
$\varphi_i(u)=\varphi_i(g)=g$ on $\o' \setminus \overline \o$ and we
can apply the previous case. It implies, for each $i \in \Nb$, the
existence of a sequence $\{u_\e^i\} \subset SBV^p(\O';\Rb^3)$
strongly converging to $\varphi_i(u)$ in $L^1(\O';\Rb^3)$ satisfying
$ u_\e^i =g_\e$ on $[\o' \setminus \overline \o] \times I$ and
$$\J(\varphi_i(u)) = \lim_{\e \to 0} \I^{g_\e}_\e(u_\e^i).$$
Since $\varphi_i(u) \to u$ in $L^1(\o';\Rb^3)$ we get that
\begin{equation}\label{alive}
\lim_{i \to +\infty}\lim_{\e \to 0}\|u_\e^i - u\|_{L^1(\O';\Rb^3)}=0.
\end{equation}
Furthermore, by (\ref{1736}) together with the lower semicontinuity of $\J$ with
respect to the strong $L^1(\o';\Rb^3)$-convergence, we obtain that
\begin{equation}\label{limie}
\J^g(u) = \lim_{i \to +\infty}\lim_{\e \to 0}\I^{g_\e}_\e(u_\e^i).
\end{equation}
A standard diagonalization argument (see e.g. Lemma 7.1 in \cite{BFF}) applied to (\ref{alive})
and (\ref{limie}) yields the existence of a sequence $i_\e \nearrow +\infty$ as $\e \to 0$ such that
$u_\e:= u_\e^{i_\e} \in SBV^p(\O',\Rb^3)$, $u_\e = g_\e$ on $[\o' \setminus \overline \o] \times I$ and
$$\J^g(u) = \lim_{\e \to 0} \I^{g_\e}_\e(u_\e).$$
\hfill$\Box$


\section{A few tools}\label{tools}

\subsection{Convergence of sets}\label{sigmapconv}

\noindent The notion of $\sigma^p$-convergence introduced in
\cite{DMFT1,DMFT2} does  not  seem to naturally provide a one
dimensional limit crack. Indeed, let $\G_n \subset \O'$ be a
sequence of $\mathcal H^2$-rectifiable sets; we denote by
$\nu_{\G_n}$ its generalized normal defined $\mathcal H^2$-a.e. on
$\G_n$. We assume that there is an a priori bound on the scaled
surface energy associated with this sequence of cracks i.e.
\begin{equation}\label{eq0}
\sup_{n \in \Nb} \int_{\G_n}\left| \left( \big(\nu_{\G_n} \big)_\a \Big| \frac{1}{\e_n} \big(\nu_{\G_n}\big)_3 \right) \right|d\mathcal H^2 < +\infty.
\end{equation}
Note that this bound will appear naturally in the energy estimates. Intuitively, we expect that any limit crack of $\G_n$ will be a subset of $\o'$ of Hausdorff dimension equal to one. But, the sequences of test functions taken in the definition of the $\sigma^p$-convergence do not contain enough information in order for this to be true. Indeed, (\ref{eq0}) implies in particular that $\mathcal H^2(\G_n) \leq C$, thus according to Lemma 4.7 in \cite{DMFT1}, we have (for a subsequence) that $\G_n$ $\sigma^p$-converges in $\O'$ to some $\mathcal H^2$-rectifiable set $\G \subset \O'$. We would like to be able to state that $\G$ is of the form $\g \times I$ for some $\mathcal H^1$-rectifiable set $\g \subset \o'$. By lower semicontinuity of
$$v \mapsto \int_{S(v)}\left| \big( \nu_v \big)_3 \right|\, d \mathcal H^2$$
for the weak $SBV^p(\O')$-convergence, we have, according to Lemma 4.3 in \cite{DMFT1} and (\ref{eq0}), that $\big(\nu_\G \big)_3(x)=0$ $\mathcal H^2$-a.e. $x \in \G$. But this does not tell us that $\G \; \widetilde = \; \g \times I$. We know, by the very definition of the $\sigma^p$-convergence, that there exists a function $u \in SBV^p(\O')$ and a sequence $u_n \rightharpoonup u$ in $SBV^p(\O')$ such that $S(u_n) \; \widetilde \subset \; \G_n$ and $\G \; \widetilde =\;  S(u)$. To prove that  $\G \; \widetilde = \; \g \times I$, it would be enough to show that $D_3 u=0$ in the sense of Radon measures. This would be immediate if the approximate scaled gradient of $u_n$ was bounded in $L^p(\O';\Rb^3)$. Since, in the sequel, we will only be interested in minimizing sequences satisfying this property, it prompts us to redefine the notion of $\sigma^p$-convergence in a 3D-2D dimensional reduction setting.

\begin{defi}\label{3d2dconv}{\rm
Let $\{\e_n\} \searrow 0^+$ and $\G_n \subset \O'$ be a sequence of $\mathcal H^2$-rectifiable sets. We say that $\G_n$ converges towards $\g$ in $\O'$ if $\g \subset \o'$, (\ref{eq0}) holds and
\begin{itemize}
\item[(a)]if $u_k \rightharpoonup u$ in $SBV^p(\O')$, $S(u_k) \; \widetilde \subset \; \G_{n_k}$ and
$$\sup_{k \in \Nb} \int_{\O'} \left| \left( \nabla_\a u_k \Big|\frac{1}{\e_{n_k}} \nabla_3 u_k \right)\right|^p dx < +\infty,$$
for a subsequence $\{\e_{n_k}\} \subset \{\e_n\}$, then $u \in SBV^p(\o')$ and $S(u)\;  \widetilde \subset\;  \g$;
\item[(b)]there exists a function $u \in SBV^p(\o')$ and a sequence $u_n \in SBV^p(\O')$ such that $u_n \rightharpoonup u$ in $SBV^p(\O')$, $S(u_n)\;  \widetilde \subset\; \G_n$,
$$\sup_{n \in \Nb} \int_{\O'} \left| \left( \nabla_\a u_n \Big|\frac{1}{\e_n} \nabla_3 u_n \right)\right|^p dx < +\infty$$
and $S(u)\; \widetilde = \;\g$.
\end{itemize}
}
\end{defi}

According to property (b) of Definition \ref{3d2dconv}, $\g$ is necessarily a $\mathcal H^1$-rectifiable set. In the following Remark, we state few properties of this kind of convergence as lower semicontinuity with respect to the Hausdorff measure and stability with respect to the inclusion.
\begin{rmk}\label{prop}
{\rm Let $\G_n \to \g$ in the sense of Definition \ref{3d2dconv}, then
\begin{enumerate}
\item for every Borel set $E \subset \o'$ such that $\mathcal H^1(E)<+\infty$ (or $E$ a compact set),
$$2 \mathcal H^1(\g \setminus E) \leq \liminf_{n \to +\infty}\mathcal H^2\big(\G_n \setminus (E \times [-1,1])\big);$$
\item if $\G_n \;\widetilde \subset \; \G_n'$ and $\G_n' \to \g'$ in the sense of Definition \ref{3d2dconv}, then $\g \; \widetilde \subset \; \g'$;
\item if $\G_n \xrightarrow[]{\sigma^p} \G$, then $\g \times I \; \widetilde \subset \; \G$
.
\end{enumerate}
}
\end{rmk}

Replacing every approximate gradients by approximate scaled gradients and using Lemma \ref{comp} instead of Ambrosio's Compactness Theorem, the exact analogues of the proofs of Lemma 4.5 and Proposition 4.6 in \cite{DMFT1} would demonstrate that any sequence of $\mathcal H^2$-rectifiable sets $\G_n \subset \O'$ satisfying (\ref{eq0}) admits a convergent subsequence in the sense of Definition (\ref{3d2dconv}). But this compactness result will not be sufficient because, in the proof of Theorem \ref{qse}, we will deal with sequence of $\mathcal H^2$-rectifiable sets which are increasing with respect to the time parameter $t$. The following Proposition, which is the analogue of Lemma 4.8 in \cite{DMFT1}, states a version of Helly's Theorem for a sequence of increasing $\mathcal H^2$-rectifiable sets.

\begin{proposition}\label{helly}
Let $[0,T] \ni t \mapsto  \G_n(t)$ a sequence of $\mathcal H^2$-rectifiable sets of $\O'$ that increases with  $t$,  i.e.
$$\G_n(s) \; \widetilde \subset \; \G_n(t) \subset \O', \quad \text{ for every }s,t \in [0,T] \text{ with }s<t.$$
Assume that
$$ \sup_{n \in \Nb} \int_{\G_n(t)}\left| \left( \big(\nu_{\G_n(t)}\big)_\a \Big| \frac{1}{\e_n} \big(\nu_{\G_n(t)}\big)_3 \right) \right|d\mathcal H^2 < +\infty,$$
uniformly in $t$. Then, there exists a subsequence $\G_{n_k}(t)$ and a $t$-increasing $\mathcal H^1$-rectifiable set $\g(t) \subset \o'$  such that for every $t \in [0,T]$, $\G_{n_k}(t)$ converges to $\g(t)$ in the sense of Definition \ref{3d2dconv}.
\end{proposition}


\subsection{Transfer of jump sets}\label{jumptt}

\noindent We now state a Jump Transfer theorem in a rescaled
version. It permits, under weak $SBV^p(\O';\Rb^3)$-convergence
assumptions of a sequence $\{u_n\}$ -- with associated bounded
scaled bulk energy -- toward its limit $u$, the transfer of the part
of the  jump set of a 2D admissible deformation that lies in the
jump set of $u$ onto that of the sequence $\{u_n\}$. The proof
relies on De Giorgi's slicing argument.

\begin{thm}[Jump Transfer]\label{JTT}
Let $\{u_n\} \subset SBV^p(\O';\Rb^3)$ and $u \in SBV^p(\o';\Rb^3)$ such that $S(u_n) \subset \overline \o \times I$, $u_n \to u$ in $L^1(\O';\Rb^3)$ and
$$M:=\sup_{n \in \Nb} \int_{\O'} \left| \left( \nabla_\a u_n \Big| \frac{1}{\e_n} \nabla_3 u_n \right) \right|^p\, dx <+\infty.$$
Then, for all $\phi \in SBV^p(\o';\Rb^3)$, there exists $\{\phi_n\} \subset SBV^p(\O';\Rb^3)$ such that
\begin{itemize}
\item $\phi_n=\phi$ a.e. on $[\o' \setminus \overline \o] \times I$,
\item $\phi_n \to \phi$ in $L^1(\O';\Rb^3)$,
\item $\ds \left( \nabla_\a \phi_n \Big| \frac{1}{\e_n} \nabla_3 \phi_n \right) \to (\nabla_\a \phi|0)$ in $L^p(\O';\Rb^{3\times 3})$,
\item $\ds \int_{[S(\phi_n) \setminus S(u_n)] \setminus [S(\phi) \setminus S(u)]} \left| \left( \big( \nu_{\phi_n} \big)_\a \Big| \frac{1}{\e_n} \big( \nu_{\phi_n} \big)_3 \right) \right| \, d\mathcal H^2 \to 0$.
\end{itemize}
\end{thm}

\noindent {\it Proof. }We first undo the scaling, coming back to the
cylinder of thickness $2\e_n$. Then, we extend periodically the
function in the transverse direction. Note that the periodic
extension may generate some additional jump at the interface of each
slice of thickness $2\e_n$. Despite this new discontinuities, we can
still apply the classical Jump Transfer Theorem (Theorem 2.1 in
\cite{FL}) and, by contradiction, we show that we can choose a slice
of thickness $2\e_n$ that satisfies good estimations. Finally, we
observe that, after translation and dilation, the restriction of the
function to this particular slice
satisfies the conclusion of Theorem \ref{JTT}.\\

{\bf Step 1. }We come back to the non rescaled cylinder $\O'_{\e_n}$ of thickness $2\e_n$. We set $v_n(x_\a,x_3):=u_n(x_\a,x_3/\e_n)$. Thus $v_n \in SBV^p(\O'_{\e_n};\Rb^3)$ and $S(v_n) \subset \overline \o \times (-\e_n,\e_n)$. Moreover,
\begin{equation}\label{v_n}
\left\{\begin{array}{l}
\ds \frac{1}{\e_n}\int_{\O'_{\e_n}}|v_n-u|\, dx=\int_{\O'}|u_n-u|\, dx,\\
\ds \frac{1}{\e_n}\int_{\O'_{\e_n}}|\nabla v_n|^p\, dx = \int_{\O'} \left| \left( \nabla_\a u_n \Big| \frac{1}{\e_n} \nabla_3 u_n \right) \right|^p\, dx.\\
\end{array}\right.
\end{equation}
We now to extend $v_n$ by periodicity in the $x_3$ direction. The
discontinuities of the resulting function will be those inherited
from the discontinuities of $v_n$ and from additional jumps that may
occur at the interface of each slice. Let
$$N_n:= \left\{
\begin{array}{ll}
\ds \frac{1}{2\e_n} - \frac{1}{2} & \text{if   } \ds \frac{1}{2\e_n}  + \frac{1}{2}\in \Nb,\\[5mm]
\ds \left \llbracket \frac{1}{2\e_n} + \frac{1}{2} \right \rrbracket
& \text{otherwise}.
\end{array}
\right.$$ For every $i\in \{-N_n, \ldots, N_n\}$, we set
$I_{i,n}:=\big( (2i-1)\e_n,(2i+1)\e_n\big)$ and $\O'_{i,n}:= \o'
\times I_{i,n}$. Note that $N_n$ is the smaller integer such that
$\O' \cap \O'_{i,n} \ne \emptyset$ for every $i\in \{-N_n, \ldots,
N_n\}$. We define the function $w_n$ on $\O'(n):= \o' \times
(-(2N_n+1)\e_n, (2N_n+1)\e_n)$ by extending $v_n$ by periodicity in
the $x_3$ direction on $\O'(n)$:
$$w_n(x_\a,x_3)=v_n(x_\a,x_3-2i\e_n) \text{ if } x_3 \in I_{i,n}.$$
Since $\O' \subset \O'(n)$, $w_n$ is a fortiori defined on $\O'$, $w_n \in SBV^p(\O';\Rb^3)$ and $S(w_n) \cap \O' \subset \overline \o \times I$.\\

{\bf Step 2. }We would like to apply the classical Jump Transfer
Theorem (Theorem 2.1 in \cite{FL}) to the function $w_n$. From
(\ref{v_n}), we have that
\begin{eqnarray*}
\int_{\O'}|w_n-u|\, dx & = & \sum_{i=-N_n}^{N_n}\int_{\O'_{i,n} \cap \O'} |v_n(x_\a,x_3-2i\e_n)-u(x_\a)|\, dx\\
 & \leq & (2N_n+1) \int_{\O'_{\e_n}} |v_n-u|\, dx\\
 & = & \e_n(2N_n+1) \int_{\O'}|u_n-u|\, dx\\
 & \leq & (1+2 \e_n) \int_{\O'}|u_n-u|\, dx \to 0
\end{eqnarray*}
and
\begin{eqnarray*}
\int_{\O'}|\nabla w_n|^p\, dx & = & \sum_{i=-N_n}^{N_n}\int_{\O'_{i,n} \cap \O'} |\nabla v_n(x_\a,x_3-2i\e_n)|^p\, dx\\
 & \leq & (2N_n+1) \int_{\O'_{\e_n} }|\nabla v_n|^p\, dx\nonumber\\
 & = & \e_n(2N_n+1) \int_{\O'}  \left| \left( \nabla_\a u_n \Big| \frac{1}{\e_n} \nabla_3 u_n \right) \right|^p\, dx\\
 & \leq & (1+2\e_n) M,
\end{eqnarray*}
which implies, thanks to De La Vall\'ee Poussin criterion (see
Proposition 1.27 in \cite{AFP}), that the sequence $\{|\nabla
w_n|\}$ is equi-integrable. We are now in position to apply Theorem
2.1 in \cite{FL} to the sequence $\{w_n\}$. Indeed, an inspection of
the proof of this result shows that the weak $L^1$-convergence
required by $\{|\nabla w_n|\}$ can be replaced, without passing to a
subsequence, by its equi-integrability (see p. 1477 in \cite{FL}).
Thus, for all $\phi \in SBV^p(\o';\Rb^3)$, we get the existence of a
sequence $\{\psi_n\} \subset SBV^p(\O';\Rb^3)$ such that
\begin{itemize}
\item $\psi_n=\phi$ a.e. on $[\o' \setminus \overline \o] \times I$,
\item $\psi_n \to \phi$ in $L^1(\O';\Rb^3)$,
\item $\nabla \psi_n \to (\nabla_\a \phi|0)$ in $L^p(\O';\Rb^{3\times 3})$,
\item $\mathcal H^2\big([S(\psi_n) \setminus S(w_n)]\setminus[S(\phi) \setminus S(u)]\big) \to 0$.
\end{itemize}

{\bf Step 3. }Since $\bigcup_{i=-N_n+1}^{N_n-1} \O'_{i,n} \subset
\O'$, we may find a $i_n \in \{ -N_n+1, \ldots, N_n-1\}$ such that
\begin{eqnarray*}
&&(2N_n -1) \Bigg\{ \int_{\O'_{i_n,n}}|\psi_n-\phi|\, dx + \int_{\O'_{i_n,n}}|\nabla \psi_n- (\nabla_\a \phi|0) |^p \, dx \\
&&\hspace{7.0cm} +\mathcal H^2_{\lfloor \O'_{i_n,n}}\big([S(\psi_n) \setminus S(w_n)]\setminus [S(\phi) \setminus S(u)]\big)\Bigg\} \\
&&\hspace{0.75cm}\leq \int_{\O'}|\psi_n-\phi|\, dx +
\int_{\O'}|\nabla \psi_n-(\nabla_\a \phi|0) |^p \, dx+ \mathcal
H^2\big([S(\psi_n) \setminus S(w_n)]\setminus[S(\phi) \setminus
S(u)]\big).
\end{eqnarray*}
Since $2N_n-1 \geq 1/\e_n-2$, we have
\begin{eqnarray*}
\frac{1}{\e_n} \left\{ \int_{\O'_{i_n,n}}|\psi_n-\phi|\, dx +\int_{\O'_{i_n,n}}|\nabla \psi_n- (\nabla_\a \phi|0) |^p \, dx +
\mathcal H^2_{\lfloor \O'_{i_n,n}}\big([S(\psi_n) \setminus S(w_n)]\setminus[S(\phi) \setminus S(u)]\big) \right\} \\
\leq 3\int_{\O'}|\psi_n-\phi|\, dx + 3\int_{\O'}|\nabla
\psi_n-(\nabla_\a \phi|0) |^p \, dx+ 3\mathcal H^2\big([S(\psi_n)
\setminus S(w_n)]\setminus[S(\phi) \setminus S(u)]\big).
\end{eqnarray*}

{\bf Step 4. }We will show that, after a translation and a dilation,
${\psi_{n}}_{\lfloor \O'_{i_n,n}}$ is the right candidate for
Theorem \ref{JTT}. Let us come back to the cylinder $\O'_{\e_n}=\o'
\times (-\e_n,\e_n)$; letting
$$\varphi_n(x_\a,x_3):={\psi_n}_{\lfloor \O'_{i_n,n}}(x_\a,x_3+2 i_n \e_n)\text{ if }x_3 \in (-\e_n,\e_n),$$
then $\varphi_n \in SBV^p(\O'_{\e_n};\Rb^3)$, $\varphi_n=\phi$ a.e.
on $[\o' \setminus \overline \o] \times (-\e_n,\e_n)$ and
\begin{eqnarray*}
\frac{1}{\e_n} \left\{ \int_{\O'_{\e_n}}|\varphi_n-\phi|\, dx + \int_{\O'_{\e_n}}|\nabla \varphi_n- (\nabla_\a \phi|0) |^p \, dx +  \mathcal H^2_{\lfloor \O'_{\e_n}}\big([S(\varphi_n) \setminus S(v_n)]\setminus[S(\phi) \setminus S(u)]\big)\right\}\\
 \leq 3\int_{\O'}|\psi_n-\phi|\, dx + 3\int_{\O'}|\nabla \psi_n-(\nabla_\a \phi|0) |^p \, dx + 3\mathcal H^2\big([S(\psi_n) \setminus S(w_n)]\setminus[S(\phi) \setminus S(u)]\big).
\end{eqnarray*}
Performing the scaling so as to come back to the unit cylinder, we get, upon setting $\phi_n(x_\a,x_3):=\varphi_n(x_\a,\e_n x_3)$, that $\phi_n \in SBV^p(\O';\Rb^3)$, $\phi_n=\phi$ a.e. on $[\o' \setminus \overline \o] \times I$ and
\begin{eqnarray*}
&&\int_{\O'}|\phi_n-\phi|\, dx + \int_{\O'}\left|\left( \nabla_\a \phi_n \Big| \frac{1}{\e_n}\nabla_3 \phi_n \right)- (\nabla_\a \phi|0) \right|^p \, dx\\
&&\hspace{4.0cm} +\int_{[S(\phi_n) \setminus S(u_n)] \setminus [S(\phi) \setminus S(u)]} \left| \left( \big( \nu_{\phi_n} \big)_\a \Big| \frac{1}{\e_n} \big( \nu_{\phi_n} \big)_3 \right) \right| \, d\mathcal H^2\\
&&\hspace{1.0cm} \leq 3\int_{\O'}|\psi_n-\phi|\, dx +
3 \int_{\O'}|\nabla \psi_n-(\nabla_\a \phi|0) |^p \,dx\\
&&\hspace{4.0cm}  + 3\mathcal H^2\big([S(\psi_n) \setminus
S(w_n)]\setminus[S(\phi) \setminus S(u)]\big) \to 0.
\end{eqnarray*}
\hfill$\Box$

\begin{rmk}\label{jtt++}{\rm
Since for $\mathcal H^2$-a.e. $x \in S(\phi_n) \cap S(\phi)$, $\nu_{\phi_n}(x)=\pm \nu_\phi(x)$, we have
\begin{eqnarray*}
&&\int_{[S(\phi_n) \setminus S(u_n)] \setminus [S(\phi) \setminus S(u)]} \left| \left( \big( \nu_{\phi_n} \big)_\a \Big| \frac{1}{\e_n} \big( \nu_{\phi_n} \big)_3 \right) \right| \, d\mathcal H^2 \\
&&\hspace{1.0cm} \geq  \int_{S(\phi_n) \setminus S(u_n)}\left| \left( \big( \nu_{\phi_n} \big)_\a \Big| \frac{1}{\e_n} \big( \nu_{\phi_n}\big)_3 \right) \right| \, d\mathcal H^2 - \int_{S(\phi) \setminus S(u)}\left| \left( \big( \nu_{\phi} \big)_\a | 0 \right) \right| \, d\mathcal H^2\\
&&\hspace{1.0cm} =  \int_{S(\phi_n) \setminus S(u_n)}\left| \left(
\big( \nu_{\phi_n} \big)_\a \Big| \frac{1}{\e_n} \big( \nu_{\phi_n}
\big)_3 \right) \right| \, d\mathcal H^2 - 2\mathcal H^1 (S(\phi)
\setminus S(u)),
\end{eqnarray*}
thus
$$\limsup_{n \to +\infty}\int_{S(\phi_n) \setminus S(u_n)}\left| \left( \big( \nu_{\phi_n} \big)_\a \Big| \frac{1}{\e_n} \big( \nu_{\phi_n} \big)_3 \right) \right| \, d\mathcal H^2 \leq 2\mathcal H^1 (S(\phi) \setminus S(u)).$$
}
\end{rmk}

The following Theorem establishes a link between the convergence in the sense of Definition \ref{3d2dconv} and the Jump Transfer Theorem. It will allow to pass to the limit in the surface energy.

\begin{thm}\label{JTT+conv}
Let $\G_n \subset \O'$ be a sequence of $\mathcal H^2$-rectifiable sets converging towards
$\g$ in the sense of Definition \ref{3d2dconv}. Then, for every $v \in SBV^p(\o';\Rb^3)$,
there exists $\{v_n\} \subset SBV^p(\O';\Rb^3)$ such that $v_n=v$ a.e. on $[\o' \setminus \overline \o] \times I$,
\begin{itemize}
\item $v_n \to v$ in $L^1(\O';\Rb^3)$,
\item $\ds \left(\nabla_\a v_n \Big|\frac{1}{\e_n} \nabla_3 v_n \right) \to (\nabla_\a v|0)$ in $L^p(\O';\Rb^{3 \times 3})$,
\item $\ds \limsup_{n \to +\infty} \int_{S(v_n) \setminus \G_n} \left|\left( \big( \nu_{v_n} \big)_\a \Big| \frac{1}{\e_n} \big( \nu_{v_n} \big)_3 \right) \right|d\mathcal H^2 \leq 2 \mathcal H^1(S(v) \setminus \g).$
\end{itemize}
\end{thm}

\noindent {\it Proof. }According to Definition \ref{3d2dconv} (b), there exists a function $u \in SBV^p(\o';\Rb^3)$ and a sequence $\{u_n\} \subset SBV^p(\O';\Rb^3)$ such that $u_n \rightharpoonup u$ in $SBV^p(\O';\Rb^3)$, $S(u_n) \; \widetilde \subset \; \G_n$,  $S(u)\; \widetilde = \;\g$ and
$$\sup_{n \in \Nb} \int_{\O'}\left| \left( \nabla_\a u_n \Big| \frac{1}{\e_n} \nabla_3 u_n\right) \right|^p dx < +\infty.$$
Theorem \ref{JTT} and Remark \ref{jtt++} yield, for any $v \in SBV^p(\o';\Rb^3)$, the existence of a sequence $\{v_n\} \subset SBV^p(\O';\Rb^3)$ such that $v_n=v$ a.e. on $[\o' \setminus \overline \o] \times I$,
\begin{itemize}
\item $v_n \to v$ in $L^1(\O';\Rb^3)$,
\item $\ds \left(\nabla_\a v_n \Big|\frac{1}{\e_n} \nabla_3 v_n \right) \to (\nabla_\a v|0)$ in $L^p(\O';\Rb^{3 \times 3})$,
\item $\ds \limsup_{n \to +\infty} \int_{S(v_n) \setminus S(u_n)} \left|\left( \big( \nu_{v_n} \big)_\a \Big| \frac{1}{\e_n}\big( \nu_{v_n} \big)_3 \right) \right|d\mathcal H^2 \leq 2 \mathcal H^1(S(v) \setminus S(u))$.
\end{itemize}
As $S(u_n) \;  \widetilde \subset \; \G_n$ and $S(u) \; \widetilde = \; \g$, we get
$$\limsup_{n \to +\infty} \int_{S(v_n) \setminus \G_n} \left|\left( \big( \nu_{v_n} \big)_\a \Big| \frac{1}{\e_n} \big( \nu_{v_n} \big)_3 \right) \right|d\mathcal H^2 \leq 2 \mathcal H^1(S(v) \setminus \g).$$
\hfill$\Box$


\subsection{Convergence of the stresses}\label{convstresses}

\noindent The energy conservation involves the derivative of the
stored energy density. Thus, we have to ensure that the $\C^1$
character of $W$ is preserved by passing to the $\G$-limit. The
following Proposition provides an answer to this question. For an
alternative proof of that result, we refer to \cite{babadjian}
Chapter 4.

\begin{proposition}\label{classC1}
Let $W : \Rb^{3\times 3} \to \Rb$ be a $\C^1$ function satisfying
(\ref{pg1}), then the function $\mathcal Q W_0 : \Rb^{3 \times 2}
\to \Rb$ is of class $\mathcal C^1$.
\end{proposition}

\noindent\textit{Proof. }According to \cite{LDR}, the function $W_0$
is continuous and satisfies
$$\frac{1}{\beta}|\overline \xi|^p - \beta \leq
W_0(\overline \xi) \leq \beta(1+|\overline \xi|^p)$$ for every
$\overline \xi \in \Rb^{3 \times 2}$. As a consequence, since $p-1>0$,
$$\liminf_{|\overline \xi| \to +\infty} \frac{W_0(\overline
\xi)}{|\overline \xi|^{p-1}} = +\infty \quad \text{ and } \quad
\limsup_{|\overline \xi| \to +\infty} \frac{W_0(\overline
\xi)}{|\overline \xi|^p}= \beta < +\infty.$$ Furthermore, for all
$\overline \xi \in \Rb^{3 \times 2}$ there exists $\xi_3 \in \Rb^3$
such that $W_0(\overline \xi)=W(\overline \xi|\xi_3)$. Since $W$
is differentiable,
$$\limsup_{|\overline \eta| \to 0}
\frac{W_0(\overline \xi + \overline \eta) - W_0(\overline \xi) - d
\cdot \overline \eta} {|\overline \eta|} \leq \limsup_{|\overline
\eta| \to 0} \frac{W((\overline \xi|\xi_3) + (\overline
\eta|0) ) - W(\overline \xi|\xi_3) -
\partial W(\overline \xi|\xi_3) \cdot (\overline \eta|0)}
{|(\overline \eta|0)|}=0,$$ where $d \in \Rb^{3 \times 2}$ is
defined by $d_{ij}:= \big(\partial W(\overline \xi|\xi_3) \big)_{ij}$ for all $i
\in \{1,2,3\}$ and all $j \in \{1,2\}$. It yields that $W_0$ is
upper semidifferentiable and the thesis follows from Theorem B in
\cite{BKK}.
\hfill$\Box$\\

From the previous Lemma, the function $\mathcal QW_0$ is of class
$\C^1$ and we denote by $\partial (\mathcal Q W_0)$ its
differential. The following result is the analogue of Lemma 4.11 in
\cite{DMFT1} in a 3D-2D dimensional reduction setting. It asserts,
under assumptions of weak $SBV^p(\O;\Rb^3)$-convergence of the
deformations together with the convergence of the bulk energy, the
weak $L^{p'}(\O;\Rb^{3 \times 3})$-convergence the stresses.

\begin{lemma}\label{stresses}
Let $\{u_n\} \subset SBV^p(\O;\Rb^3)$ and $u \in SBV^p(\o;\Rb^3)$
such that $u_n \rightharpoonup u$ in $SBV^p(\O;\Rb^3)$ and
\begin{equation}\label{hyp}
\int_\O W\left( \nabla_\a u_n \Big|\frac{1}{\e_n} \nabla_3 u_n \right)\, dx \to 2 \int_\o \mathcal Q W_0 (\nabla_\a u)\, dx_\a.
\end{equation}
Then,
$$\partial W \left(\nabla_\a u_n \Big|\frac{1}{\e_n} \nabla_3 u_n \right) \rightharpoonup \big(\partial (\mathcal Q W_0)(\nabla_\a u)|0 \big) \quad \text{ in } L^{p'}(\O;\Rb^{3 \times 3}).$$
\end{lemma}

\noindent {\it Proof. }Let $\Psi \in L^p(\O;\Rb^{3 \times 3})$, we
denote by $\overline \Psi \in  L^p(\O;\Rb^{3 \times 2})$ the
restriction of $\Psi$ to $\Rb^{3 \times 2}$ i.e. $\overline
\Psi_{ij}=\Psi_{ij}$ if $i \in \{1,2,3\}$ and $j \in \{1,2\}$. It is
enough to show that \begin{eqnarray*} \int_\O \partial (\mathcal
QW_0)(\nabla_\a u)\cdot \overline \Psi\, dx & = & \int_\O
\big(\partial (\mathcal QW_0)(\nabla_\a u) |0 \big)\cdot \Psi\, dx\\
& \leq & \liminf_{n \to +\infty} \int_\O \partial W \left(\nabla_\a
u_n \Big|\frac{1}{\e_n} \nabla_3 u_n \right) \cdot \Psi\, dx.
\end{eqnarray*} Let $h_k \searrow 0^+$, according to Remark \ref{qcvx},
Theorem 5.29 in \cite{AFP}  we have
\begin{eqnarray*}
\int_\O \mathcal Q W_0(\nabla_\a u +h_k \overline \Psi)\, dx & \leq
& \liminf_{n\to +\infty} \int_\O \mathcal Q W_0(\nabla_\a u_n+ h_k \overline \Psi) \, dx\\
& \leq & \liminf_{n\to +\infty} \int_\O W\left(\left(\nabla_\a
u_n\Big|\frac{1}{\e_n}\nabla_3u_n\right) + h_k \Psi \right) \, dx.
\end{eqnarray*}
As a consequence, from (\ref{hyp}) we get that
\begin{eqnarray*}
&&\int_\O\frac{\mathcal Q W_0(\nabla_\a u +h_k \overline \Psi) -\mathcal Q W_0(\nabla_\a u)}{h_k}\, dx\\
&&\hspace{2.0cm} \leq \liminf_{n\to +\infty} \int_\O \frac{1}{h_k} \left[ W\left(\left(\nabla_\a u_n\Big|\frac{1}{\e_n}\nabla_3u_n\right) + h_k \Psi \right) -W \left(\nabla_\a u_n \Big|\frac{1}{\e_n}\nabla_3 u_n\right)\right] \, dx.\\
\end{eqnarray*}
We may find a $n_k \in \Nb$ such that
\begin{eqnarray*}
&&\int_\O\frac{\mathcal Q W_0(\nabla_\a u +h_k \overline \Psi) -\mathcal Q W_0(\nabla_\a u)}{h_k}\, dx -\frac{1}{k}\\
&&\hspace{2.0cm} \leq \int_\O \frac{1}{h_k} \left[ W\left(\left(\nabla_\a u_n\Big|\frac{1}{\e_n}\nabla_3u_n\right) + h_k \Psi \right) -W \left(\nabla_\a u_n \Big|\frac{1}{\e_n}\nabla_3 u_n\right)\right] \, dx\\
\end{eqnarray*}
for all $n \geq n_k$. We define $\eta_n=h_k$ if $n_k \leq n \leq
n_{k+1}$ and pass to the limit when $n \to +\infty$. Since $\Phi
\mapsto \int_\O W(\Phi)\, dx$ is a $\C^1$-map from $L^p(\O;\Rb^{3
\times 3})$ to $\Rb$ with differential $\Psi \mapsto \int_\O
\partial W(\Phi) \cdot \Psi\, dx$, it follows that
\begin{eqnarray*}
&&\lim_{n \to +\infty}\int_\O\frac{\mathcal Q W_0(\nabla_\a u +\eta_n \overline \Psi) -\mathcal Q W_0(\nabla_\a u)}{\eta_n}\, dx\\
&&\hspace{2.0cm} \leq \liminf_{n\to +\infty} \int_\O \frac{1}{\eta_n} \left[ W\left(\left(\nabla_\a u_n\Big|\frac{1}{\e_n}\nabla_3u_n\right)
+ \eta_n \Psi \right) -W \left(\nabla_\a u_n \Big|\frac{1}{\e_n}\nabla_3 u_n\right)\right] \, dx\\
&&\hspace{2.0cm}  \leq \liminf_{n\to +\infty} \int_\O \partial
W\left(\left(\nabla_\a u_n\Big|\frac{1}{\e_n}\nabla_3 u_n\right) +
\tau_n \Psi \right) \cdot \Psi\, dx,
\end{eqnarray*}
for some $\tau_n \in [0,\eta_n]$. Lebesgue's Dominated Convergence
Theorem in the left hand side, together with Lemma 4.9 of
\cite{DMFT1} in the right hand side yield
\begin{eqnarray*}
\int_\O \partial (\mathcal Q W_0)(\nabla_\a u)\cdot \overline \Psi\, dx & = & \lim_{n \to +\infty} \int_\O\frac{\mathcal Q W_0(\nabla_\a u +\eta_n \overline \Psi) - \mathcal Q W_0(\nabla_\a u)}{\eta_n}\, dx\\
 & \leq & \liminf_{n \to +\infty} \int_\O \partial W\left( \left(\nabla_\a u_{n}  \Big| \frac{1}{\e_{n}}\nabla_3 u_{n}\right)+ \tau_n \Psi \right) \cdot \Psi\, dx\\
 & = & \liminf_{n \to +\infty}\int_\O \partial W\left(\nabla_\a u_{n}\Big|\frac{1}{\e_{n}}\nabla_3 u_{n} \right) \cdot \Psi\, dx.
\end{eqnarray*}
\hfill$\Box$


\section{Convergence of the quasistatic evolution}\label{convqse}

\noindent The first step of the analysis consists in defining a
limit deformation field and a crack. This is done in Lemma
\ref{eto0} by means of energy estimates which are possible,  thanks
to the $L^\infty$-boundness assumption (\ref{Hyp}) and to the bound
of the prescribed boundary deformation (\ref{CL}). The limit
deformation $u(t)$ turns out to be the weak $SBV^p(\O';\Rb^3)$-limit
of $u^\e(t)$ while the limit crack $\g(t)$ is obtained through the
convergence of $\G^\e(t)$ in the sense of Definition \ref{3d2dconv}.
Then, we derive a minimality property for $u(t)$. At time $t=0$ in
Lemma \ref{MIN0}, we use a $\G$-convergence argument. This is
possible because, in the absence of  preexisting cracks, the surface
term of the energy at  time $0$  is precisely that introduced in the
$\G$-limit analysis of Section \ref{gammaconvergence}. Nevertheless,
we cannot proceed in this way for the next times in Lemma \ref{MIN}
because of the presence of $\G^\e(t)$ in the surface term. We need
here to construct directly a sequence thanks to the Jump Transfer
Theorem, Theorem \ref{JTT}. Then, we show that $(u(t),\g(t))$ is a
quasistatic evolution for the relaxed model by proving that the
energy conservation holds. To do this, we use, on the one hand, the
approximation of the Lebesgue integral by Riemann sums  in Lemma
\ref{E1}, and, on the other hand, the convergence of the total
energy at the initial time (that can be proved directly) together
with the weak convergence of the stresses and the strong convergence
assumption (\ref{CLcvforte}) in Lemma \ref{E2}. Finally, in Lemma
\ref{convE}, we show the convergence
of the total energy at any time.\\

\subsection{Energy estimates and compactness}\label{NRJestim}

\begin{lemma}\label{eto0}
There exists a subsequence $\{\e_{n}\} \searrow 0^+$, a deformation
field  $u(t) \in SBV^p(\o';\Rb^3)$ satisfying $u(t)=g(t)$ $\mathcal
L^2$-a.e. on $\o' \setminus \overline \o$, and a time-increasing crack   $\g(t)
\subset \overline \o$ such that, for
every $t \in [0,T]$, $S(u(t)) \; \widetilde \subset \; \g(t)$.
Moreover, $\G^{\e_n}(t)$ converges to $\g(t)$ in the sense of
Definition \ref{3d2dconv}, $u^{\e_n}(0) \rightharpoonup u(0)$ in
$SBV^p(\O';\Rb^3)$ and, for every $t \in (0,T]$, there exists a
$t$-dependent subsequence $\{\e_{n_t}\} \subset \{\e_n\}$ such that
$u^{\e_{n_t}}(t) \rightharpoonup u(t)$ in $SBV^p(\O';\Rb^3)$.
\end{lemma}

\noindent {\it Proof. }Firstly, at time $t=0$, we test the
minimality of $u^\e(0)$ with $v=g^\e(0)$. Since $S(u^\e(0)) \;
\widetilde{=} \; \G^\e(0)$, we deduce by (\ref{CL}) and the growth
condition (\ref{pg1}) that $\E_\e(0) \leq C$.

Then, we take $v=g^\e(t)$ as test function in (\ref{min_t_e}) at
time $t$. As $S(u^\e(t)) \; \widetilde{\subset} \; \G^\e(t)$, it
follows from (\ref{CL}) together with the growth condition
(\ref{pg1}) that
\begin{eqnarray}\label{Ee0}
\int_{\O} W\left( \nabla_\a u^\e(t) \Big| \frac {1}{\e} \nabla_3 u^\e(t) \right)\, dx \leq \int_\O W\left( \nabla_\a g^\e(t) \Big| \frac {1}{\e} \nabla_3 g^\e(t) \right)\, dx \leq C.
\end{eqnarray}
Thus, H\"older's inequality, (\ref{pg2}), (\ref{Etot}) and
(\ref{Ee0}) imply the existence of a constant $C>0$, independent of
$t$ and $\e$ such that for every $t \in [0,T]$, $\E_\e(t) \leq C$.
Hence, by the coercivity condition (\ref{pg1}),
\begin{equation}\label{borne}
\int_{\O} \left|\left( \nabla_\a u^\e(t) \Big| \frac {1}{\e} \nabla_3 u^\e(t) \right)\right|^p\, dx + \int_{\G^\e(t)} \left| \left( \left( \nu_{\G^\e(t)} \right)_\a \Big| \frac{1}{\e} \left( \nu_{\G^\e(t)} \right)_3 \right) \right| \, d \mathcal H^2 \leq C.
\end{equation}
In view of Proposition \ref{helly}, we may find a subsequence
$\{\e_n\} \searrow 0^+$ and an $\mathcal H^1$-rectifiable set $\g(t)
\subset \o'$, increasing in time, such that $\G^{\e_n}(t)$ converges
to $\g(t)$ in the sense of Definition \ref{3d2dconv}. According to
Remark \ref{prop}-1, since $\overline \o$ is compact, $\g(t) \;
\widetilde \subset \; \overline \o$ and
\begin{equation}\label{sci}
2 \mathcal H^1(\g(t)) \leq \liminf_{n \to +\infty} \mathcal
H^2(\G^{\e_n}(t)) \leq \liminf_{n \to +\infty} \int_{\G^{\e_n}(t)}
\left| \left( \left( \nu_{\G^{\e_n}(t)} \right)_\a \Big|
\frac{1}{\e_n} \left( \nu_{\G^{\e_n}(t)} \right)_3 \right) \right|
\, d \mathcal H^2.
\end{equation}
As $S(u^{\e}(t)) \; \widetilde{\subset} \;  \G^{\e}(t)$ and
$u^{\e}(t)=g^{\e}(t)$ $\mathcal L^3$-a.e. on $[\o' \setminus
\overline \o] \times I$, we have by (\ref{Hyp}) and (\ref{borne})
\begin{eqnarray*}
&&\|u^{\e}(t)\|_{L^\infty(\O';\Rb^3)} + \int_{\O'} \left|\left(
\nabla_\a u^{\e}(t) \Big| \frac {1}{\e} \nabla_3 u^{\e}(t)
\right)\right|^p\, dx\nonumber\\
&&\hspace{5.0cm} + \int_{S(u^{\e}(t))} \left| \left( \left(
\nu_{u^{\e}(t)} \right)_\a \Big| \frac{1}{\e} \left( \nu_{u^{\e}(t)}
\right)_3 \right) \right| \, d \mathcal H^2 \leq C,
\end{eqnarray*}
for some constant $C>0$ independent of $\e>0$ and $t \in [0,T]$.
 We insist, once again, on the fact that we do not try to
justify the boundness assumption  on $u^\e(t)$. In view of Lemma
\ref{comp}, there exists  a further subsequence of $\{\e_n\}$ (still
denoted by $\{\e_n\}$) and $u(0) \in SBV^p(\o';\Rb^3)$ such that
$u^{\e_n}(0) \rightharpoonup u(0)$ in $SBV^p(\O';\Rb^3)$. Moreover,
as $u^{\e_n}(0)=g^{\e_n}(0)$ $\mathcal L^3$-a.e. on $[\o' \setminus
\overline \o] \times I$, from (\ref{CLcvforte}) we get $u(0)=g(0)$
$\mathcal L^2$-a.e. on $\o' \setminus \overline \o$. Thanks to
condition (a) of Definition \ref{3d2dconv}, we deduce that
$S(u(0)) \; \widetilde \subset \; \g(0)$.\\

We set for a.e. $t \in [0,T]$,
\begin{equation}\label{theta}
\left\{\begin{array}{rcl}
\theta_n(t) & := & \ds \int_\O \partial W \left( \nabla_\a u^{\e_n}(t) \Big| \frac{1}{\e_n} \nabla_3 u^{\e_n}(t) \right) \cdot \left( \nabla_\a \dot g^{\e_n}(t) \Big| \frac{1}{\e_n} \nabla_3 \dot g^{\e_n}(t) \right) dx,\\
&&\\
\theta(t) & := & \ds \limsup_{n \to +\infty}\; \theta_n(t).
\end{array}\right.
\end{equation}
From (\ref{pg2}), (\ref{CL}), (\ref{borne}), $\theta \in L^1(0,T)$
and by virtue of Fatou's Lemma
\begin{equation}\label{scs}
\limsup_{n \to +\infty} \int_0^t \theta_n(s)\, ds \leq \int_0^t \theta(s)\, ds.
\end{equation}
For a.e. $t \in [0,T]$, we extract a $t$-dependent subsequence $\{n_t\}$ such that
\begin{equation}\label{sst}
\theta(t) =  \lim_{n_t \to +\infty}\theta_{n_t}(t).
\end{equation}
Lemma \ref{comp} implies that for every $t \in (0,T]$, upon
extracting a further subsequence (not relabeled), $u^{\e_{n_t}}(t)
\rightharpoonup u(t)$ in $SBV^p(\O';\Rb^3)$ for some $u(t) \in
SBV^p(\o';\Rb^3)$. Moreover, as $u^{\e_{n_t}}(t)=g^{\e_{n_t}}(t)$
$\mathcal L^3$-a.e. on $[\o' \setminus \overline \o] \times I$, from
(\ref{CLcvforte}), we get $u(t)=g(t)$ $\mathcal L^2$-a.e. on $\o' \setminus
\overline \o$. By condition (a) of Definition \ref{3d2dconv}, we get
that $S(u(t)) \; \widetilde \subset \; \g(t)$. \hfill$\Box$


\subsection{Minimality property}\label{minprop}

\noindent For all $t \in [0,T]$, we define the limit energy by
\begin{equation}\label{E(t)}
\E(t):=2 \int_\o  \mathcal Q W_0(\nabla _\a u(t))\, dx_\a +2 \mathcal H^1(\g(t)).
\end{equation}

Our goal is to show that $u(t)$ satisfies some minimality property inherited from that of $u^\e(t)$. We will distinguish the initial time from the subsequent times. At time $t=0$, we will further show the convergence of the bulk and the surface energy to their two-dimensional counterpart respectively, for the subsequence $\{\e_n\}$. Concerning the next times, we will only be able to prove the convergence of the bulk energy toward its two-dimensional analogue, for the $t$-dependent subsequence $\{\e_{n_t}\}$. The convergence of the total energy, or equivalently of the surface energy, will be established later in Lemma \ref{convE} for a subsequence independent of the time.

\begin{lemma}\label{MIN0}
At time $t=0$, $u(0)$ minimizes
$$v \mapsto 2\int_\o \mathcal Q W_0(\nabla_\a v)\, dx_\a + 2 \mathcal H^1 \big(S(v)\big),$$
among $\{v \in SBV^p(\o';\Rb^3):\quad v=g(0) \text{ a.e. on } \o' \setminus \overline \o\}$. Moreover, $\g(0) \; \widetilde = \; S(u(0))$ and we have
$$\left\{\begin{array}{l}
\ds \int_\O W\left( \nabla_\a u^{\e_{n}}(0) \Big| \frac {1}{\e_{n}} \nabla_3 u^{\e_{n}}(0) \right)\, dx \to 2\int_\o \mathcal Q W_0(\nabla_\a u(0))\, dx_\a,\\
\ds  \int_{\G^{\e_n}(0)}\left|\left( \big( \nu_{\G^{\e_{n}}(0)} \big)_\a \Big|\frac{1}{\e_{n}} \big( \nu_{\G^{\e_{n}}(0)} \big)_3 \right)\right| d\mathcal H^2 \to 2\mathcal H^1(\g(0)).
\end{array}\right.$$
In particular, $\E_{\e_{n}}(0) \to \E(0)$ and
$$\partial W \left(\nabla_\a u^{\e_{n}}(0) \Big|\frac{1}{\e_{n}} \nabla_3 u^{\e_{n}}(0) \right) \rightharpoonup \big( \partial (\mathcal Q W_0)(\nabla_\a u(0)) |0 \big) \quad \text{ in } L^{p'}(\O;\Rb^{3 \times 3}).$$
\end{lemma}

\noindent {\it Proof. }Let $v \in  SBV^p(\o';\Rb^3)$ such that
$v=g(0)$ $\mathcal L^2$-a.e. on $\o' \setminus \overline \o$. By
virtue of Corollary \ref{corol}, there exists a sequence $\{w_n\}
\subset SBV^p(\O';\Rb^3)$ satisfying $w_n=g^{\e_n}(0)$ $\mathcal
L^3$-a.e. on $[\o' \setminus \overline \o] \times I$, $w_n \to v$ in
$L^1(\O';\Rb^3)$ and
\begin{eqnarray}\label{w_n}
&&2\int_\o \mathcal Q W_0(\nabla_\a v)\, dx_\a + 2 \mathcal H^1 \big(S(v)\big)\nonumber\\
&&\hspace{1.5cm} = \lim_{n \to +\infty}\left[ \int_\O W\left(\nabla_\a w_n \Big|\frac{1}{\e_n}\nabla_3 w_n \right)dx + \int_{S(w_n)}\left|\left( \big( \nu_{w_n} \big)_\a \Big|\frac{1}{\e_n} \big( \nu_{w_n} \big)_3 \right)\right| d\mathcal H^2 \right].
\end{eqnarray}
Taking $w_n$ as test function in  the minimality condition for $u^{\e_n}(0)$ we get
\begin{eqnarray}\label{mint0}
&& \int_\O W\left(\nabla_\a u^{\e_n}(0) \Big|\frac{1}{\e_n}\nabla_3 u^{\e_n}(0) \right)dx + \int_{S(u^{\e_n}(0))}\left|\left( \big( \nu_{u^{\e_n}(0)} \big)_\a \Big|\frac{1}{\e_n} \big( \nu_{u^{\e_n}(0)} \big)_3 \right)\right| d\mathcal H^2\nonumber\\
&&\hspace{1.5cm}\leq \int_\O W\left(\nabla_\a w_n \Big|\frac{1}{\e_n}\nabla_3 w_n \right)dx + \int_{S(w_n)}\left|\left( \big( \nu_{w_n} \big)_\a \Big|\frac{1}{\e_n} \big( \nu_{w_n} \big)_3 \right)\right| d\mathcal H^2.
\end{eqnarray}
Remark \ref{qcvx} and Theorem 5.29 in \cite{AFP} yield
\begin{eqnarray}\label{vli}
2\int_\o \mathcal QW_0(\nabla_\a u(0))\, dx_\a & \leq & \liminf_{n \to +\infty}  \int_\O \mathcal QW_0(\nabla_\a u^{\e_{n}}(0))\, dx\nonumber\\
& \leq & \liminf_{n \to +\infty} \int_\O W\left(\nabla_\a u^{\e_{n}}(0) \Big|\frac{1}{\e_{n}}\nabla_3 u^{\e_{n}}(0) \right)dx
\end{eqnarray}
and thanks to (\ref{sci}) together with the fact that $\G^{\e_n}(0) \; \widetilde = \; S(u^{\e_n}(0))$,
\begin{eqnarray}\label{sli}
2 \mathcal H^1(\g(0)) & \leq & \liminf_{n \to +\infty}  \mathcal H^2(\G^{\e_{n}}(0)) \nonumber\\
 & = & \liminf_{n \to +\infty}  \mathcal H^2(S(u^{\e_{n}}(0))) \nonumber\\
 & \leq & \liminf_{n \to +\infty}\int_{S(u^{\e_{n}}(0))}\left|\left( \big( \nu_{u^{\e_{n}}(0)} \big)_\a \Big|\frac{1}{\e_{n}} \big( \nu_{u^{\e_{n}}(0)} \big)_3 \right)\right| d\mathcal H^2.
\end{eqnarray}
Finally, from (\ref{w_n}), (\ref{mint0}), (\ref{vli}) and (\ref{sli}) we get by letting $n \to +\infty$,
$$2\int_\o \mathcal QW_0(\nabla_\a u(0))\, dx_\a + 2 \mathcal H^1(\g(0)) \leq 2\int_\o \mathcal Q W_0(\nabla_\a v)\, dx_\a + 2 \mathcal H^1 \big(S(v)\big).$$
Taking $v=u(0)$ in the previous inequality, we observe that $\mathcal H^1(\g(0)) \leq \mathcal H^1\big( S(u(0)) \big)$, which implies, as $S(u(0)) \; \widetilde \subset \; \g(0)$, that $S(u(0)) \; \widetilde = \; \g(0)$. It establishes the minimality property satisfied by $u(0)$. Taking still $v=u(0)$, (\ref{mint0}) and (\ref{sli}) give
$$\limsup_{n \to +\infty}  \int_\O W\left( \nabla_\a u^{\e_{n}}(0) \Big| \frac {1}{\e_{n}} \nabla_3 u^{\e_{n}}(0) \right)\, dx \leq 2\int_\o \mathcal QW_0(\nabla_\a u(0))\, dx_\a.$$
and this shows with (\ref{vli}) that
$$\int_\O W\left( \nabla_\a u^{\e_{n}}(0) \Big| \frac {1}{\e_n} \nabla_3 u^{\e_{n}}(0) \right)\, dx \to 2\int_\o \mathcal Q W_0(\nabla_\a u(0))\, dx_\a.$$
We report in (\ref{mint0}) and  obtain
$$\limsup_{n \to +\infty} \int_{S(u^{\e_{n}}(0))}\left|\left( \big( \nu_{u^{\e_{n}}(0)} \big)_\a \Big|\frac{1}{\e_{n}} \big( \nu_{u^{\e_{n}}(0)} \big)_3 \right)\right| d\mathcal H^2 \leq 2\mathcal H^1(S(u(0))),$$
which implies together with (\ref{sli}) that
$$ \int_{S(u^{\e_{n}}(0))}\left|\left( \big( \nu_{u^{\e_{n}}(0)} \big)_\a \Big|\frac{1}{\e_{n}} \big( \nu_{u^{\e_{n}}(0)} \big)_3 \right)\right| d\mathcal H^2 \to 2\mathcal H^1(S(u(0))).$$
This yields that $\E_{\e_{n}}(0) \to \E(0)$ and the convergence of the
stresses follows from Lemma \ref{stresses}.
\hfill$\Box$

\begin{rmk}\label{2005}
{\rm It is immediate from the previous lemma that $u(0)$ minimizes
$$v \mapsto 2\int_\o \mathcal Q W_0(\nabla_\a v)\, dx_\a + 2 \mathcal H^1 \big(S(v)\setminus \g(0)\big),$$
among $\{v \in SBV^p(\o';\Rb^3) : \quad v=g(0) \text{ a.e. on } \o'
\setminus \overline \o\}$.
}
\end{rmk}

\begin{rmk}\label{crack}{\rm
Note that the previous result holds because we did not allow the
body to contain a preexisting crack. Indeed, in this case, since the
energy we are minimizing at the initial time is the same as the
functional involved in the $\G$-limit analysis, we can take as
competitor in the minimization a recovery sequence. It permits us to
show the convergence of the total energy at time $t=0$; that is
essential if one is to prove that it still holds true at subsequent
times in Lemma \ref{E2}. If we had considered a body containing a
preexisting crack, we would be unable to obtain such a convergence,
but only a convergence of the bulk energy. Indeed, if $\G_0^\e
\subset \overline \o \times I$ denoted a preexisting (rescaled)
crack with bounded scaled surface energy, then, according to the
formulation in \cite{DMFT1}, $(u_0^\e,\G_0^\e)$ would have to
minimize
$$(v,\G) \mapsto \int_\O W \left(\nabla_\a v \Big|\frac{1}{\e} \nabla_3
v \right) dx + \int_{\G}\left| \left( \big(\nu_{\G}\big)_\a \Big|
\frac{1}{\e} \big(\nu_{\G}\big)_3 \right) \right| d\mathcal H^2$$
among every $\mathcal H^2$-rectifiable crack $\G \subset \overline
\o \times I$ with $\G_0^\e \; \widetilde \subset \; \G$, and every
deformation $v \in SBV^p(\O';\Rb^3)$ such that $v=g^\e(t)$ a.e. on
$[\o' \setminus \overline \o] \times I$ and $S(v) \; \widetilde
\subset \; \G$. In particular, setting $\G:= \G_0^\e \cup S(v)$ for
all $v \in SBV^p(\O';\Rb^3)$ satisfying $v=g^\e(0)$ a.e. on $[\o'
\setminus \overline \o] \times I$, we would get that $u_0^\e$ must
minimize
$$v\mapsto \int_\O W \left(\nabla_\a v \Big|\frac{1}{\e} \nabla_3
v \right) dx + \int_{S(v) \setminus \G_0^\e}\left| \left(
\big(\nu_{v}\big)_\a \Big| \frac{1}{\e} \big(\nu_{v}\big)_3 \right)
\right| d\mathcal H^2$$ among such $v$'s. Hence, since by Theorem
3.15 in \cite{DMFT1} (or Theorem 2.1 in \cite{DMFT2})
$(u^\e(0),\G^\e(0)) = (u_0^\e,\G_0^\e)$, the argument used in the
proof of Lemma \ref{MIN0} would not hold anymore. We would only be
able to state, as in the following Lemma \ref{MIN}, the convergence
of the bulk energy. Unfortunately, the convergence of the surface
energy would then remain an open question. }
\end{rmk}

We are now going to state a minimality property satisfied by $u(t)$
for $t \in (0,T]$. The following result ensures the convergence of
the three-dimensional bulk energy to its two-dimensional counterpart
for a $t$-dependent subsequence. But  the convergence of the total energy, or
equivalently of the surface energy, cannot be established at this stage in a manner similar to that
used in Lemma \ref{MIN0}  at the initial time.

\begin{lemma}\label{MIN}
For every $t \in (0,T]$, $u(t)$ minimizes
$$v \mapsto 2\int_\o \mathcal Q W_0(\nabla_\a v)\, dx_\a + 2 \mathcal H^1 \big(S(v)\setminus \g(t)\big),$$
among $\{v \in SBV^p(\o';\Rb^3) : \quad v=g(t) \text{ a.e. on } \o' \setminus \overline \o\}$. Moreover, we have
$$\int_\O W\left( \nabla_\a u^{\e_{n_t}}(t) \Big| \frac {1}{\e_{n_t}} \nabla_3 u^{\e_{n_t}}(t) \right)\, dx \to 2\int_\o \mathcal Q W_0(\nabla_\a u(t))\, dx_\a.$$
In particular,
$$\partial W \left(\nabla_\a u^{\e_{n_t}}(t) \Big|\frac{1}{\e_{n_t}} \nabla_3 u^{\e_{n_t}}(t) \right) \rightharpoonup \big( \partial (\mathcal Q W_0)(\nabla_\a u(t)|0) \big) \quad \text{ in } L^{p'}(\O;\Rb^{3 \times 3})$$
and thus, for a.e. $t \in [0,T]$,
\begin{equation}\label{theta0}
\theta(t)= 2 \int_\o \partial ( \mathcal Q W_0)(\nabla_\a u(t) ) \cdot \nabla_\a \dot g(t)\, dx_\a.
\end{equation}
\end{lemma}

\noindent {\it Proof. }We first prove the minimality property.
Unlike Lemma  \ref{MIN0}, we cannot use a $\G$-convergence argument
because of the presence of an $\e$-dependent crack in the surface
term. We will construct a minimizing sequence with the help of the
Jump Transfer Theorem.

Let $w \in SBV^p(\o';\Rb^3)$ such that $w=g(t)$ $\mathcal L^2$-a.e. on $\o' \setminus \overline \o$.
Since $\G^{\e_n}(t)$ converges to $\g(t)$ in the sense of Definition \ref{3d2dconv}, from Theorem
\ref{JTT+conv}, there exists a sequence $\{w_n\} \subset SBV^p(\O';\Rb^3)$ satisfying $w_n=w=g(t)$ a.e. on $[\o' \setminus \overline \o] \times I$, $w_n \to w$ in $L^1(\O';\Rb^3)$ and
\begin{equation}\label{minJTT}
\left\{\begin{array}{l}
\ds \left( \nabla_\a w_n \Big| \frac{1}{\e_{n}} \nabla_3 w_n \right) \to (\nabla_\a w|0)\text{ in }L^p(\O';\Rb^{3\times 3}),\\
\ds  \limsup_{n \to +\infty} \int_{S(w_n) \setminus \G^{\e_{n}}(t)} \left| \left( \big( \nu_{w_n} \big)_\a \Big| \frac{1}{\e_{n}} \big( \nu_{w_n} \big)_3 \right) \right| \, d\mathcal H^2 \leq 2 \mathcal H^1(S(w) \setminus \g(t)).
\end{array}\right.
\end{equation}
A measurable selection criterion (see e.g. \cite{ET}) together with
the coercivity condition (\ref{pg1}) imply the existence of $z \in
L^p(\o;\Rb^3)$ such that $W_0(\nabla_\a w)=W(\nabla_\a w|z)$
$\mathcal L^2$-a.e. in $\o$. By density, there exists a sequence
$z_j \in \mathcal C^\infty_c(\O;\Rb^3)$ such that $z_j \to H(t)-z$
in $L^p(\O;\Rb^3)$ where $H(t)$ is defined in (\ref{CLcvforte}).
Denoting by $b_j:= \int_{-1}^{x_3}z_j(\cdot,s)\, ds$, we take
$w_n+g^{\e_{n}}(t)-g(t)-\e_{n} b_j$ as test function in
(\ref{min_t_e}) and we get
\begin{eqnarray*}
&&\int_{\O} W\left( \nabla_\a u^{\e_n}(t) \Big| \frac {1}{\e_n} \nabla_3 u^{\e_n}(t) \right) \, dx\\
&&\hspace{1.0cm} \leq \int_{\O} W \left( \nabla_\a w_n + \nabla_\a g^{\e_n}(t) - \nabla_\a g(t) - \e_n \nabla _\a b_j \Big|\frac {1}{\e_n} \nabla_3 w_n + \frac {1}{\e_n} \nabla_3 g^{\e_n}(t) -z_j \right)\, dx\\
&&\hspace{5.0cm}  + \int_{S(w_n) \setminus \G^{\e_n}(t)} \left| \left( \big( \nu_{w_n} \big)_\a \Big| \frac{1}{\e_n} \big( \nu_{w_n} \big)_3 \right) \right| \, d \mathcal H^2.
\end{eqnarray*}
We replace $n$ by $n_t$ (see Lemma \ref{eto0}) and pass to the limit when $n_t$ tends to $+\infty$. In view of Theorem 5.29 in \cite{AFP} and of Remark \ref{qcvx},
\begin{eqnarray}\label{529}
2\int_\o \mathcal Q W_0(\nabla_\a u(t))\, dx_\a  & \leq  & \liminf_{n_t \to +\infty} \int_{\O} \mathcal Q W_0 ( \nabla_\a u^{\e_{n_t}}(t)) \, dx\nonumber\\
 & \leq & \liminf_{n_t \to +\infty} \int_{\O} W\left( \nabla_\a u^{\e_{n_t}}(t) \Big| \frac {1}{\e_{n_t}} \nabla_3 u^{\e_{n_t}}(t) \right) \, dx.
\end{eqnarray}
Thus, using (\ref{CLcvforte}) and (\ref{minJTT}) in the right hand side, we get
$$2\int_\o \mathcal Q W_0(\nabla_\a u(t))\, dx_\a \leq \int_\O W(\nabla_\a w|H(t)-z_j)\, dx +2 \mathcal H^1(S(w)\setminus \g(t)).$$
Passing to the limit when $j \to +\infty$ we obtain
\begin{eqnarray*}
2\int_\o \mathcal Q W_0(\nabla_\a u(t))\, dx_\a & \leq & \int_\O W(\nabla_\a w|z)\, dx +2 \mathcal H^1(S(w)\setminus \g(t) )\\
 & =& 2\int_\o W_0(\nabla_\a w)\, dx_\a +2 \mathcal H^1(S(w) \setminus \g(t) ).
\end{eqnarray*}
We would like to replace $W_0$ by its quasiconvexification in the
previous relation. To this end, we use a relaxation argument. First
of all, we approach $\g(t)$ from inside by a compact set, so as to
work on an open subset of $\o$. This is possible because, since
$\mathcal H^1(\g(t)) <+\infty$, then $\mathcal H^1_{\lfloor \g(t)}$
is a Radon measure. Thus, for any $\eta>0$, there exists a compact
set $K_t^\eta \subset \g(t)$ such that $\mathcal H^1(\g(t) \setminus
K_t^\eta) \leq \eta$. In particular,
\begin{equation}\label{relax}
2\int_\o \mathcal Q W_0(\nabla_\a u(t))\, dx_\a \leq 2\int_\o
W_0(\nabla_\a w)\, dx_\a +2 \mathcal H^1(S(w) \setminus K_t^\eta ).
\end{equation}
Let $v \in SBV^p(\o';\Rb^3)$ satisfying $v=g(t)$ $\mathcal L^2$-a.e.
on $\o' \setminus \overline \o$. In view of Theorem 8.1 together
with Remark 8.2 in \cite{BCP} and arguing as in the proof of Lemma
\ref{lbc} and Corollary \ref{corol}, it is easily deduced that there
exists a sequence $\{w_k\} \subset SBV^p(\o';\Rb^3)$ such that $w_k
\to v$ in $L^1(\o';\Rb^3)$, $w_k =g(t)$ $\mathcal L^2$-a.e. on $\o'
\setminus \overline \o$ and
$$\int_\o \mathcal Q W_0 (\nabla_\a v)\, dx_\a + \mathcal H^1(S(v) \setminus K_t^\eta) = \lim_{k \to +\infty} \left[ \int_\o W_0 (\nabla_\a w_k)\, dx_\a + \mathcal H^1(S(w_k) \setminus K_t^\eta) \right].$$
In (\ref{relax}), we replace $w$ by $w_k$ and we pass to the limit
when $k \to +\infty$; we get
\begin{eqnarray*}
2\int_\o \mathcal Q W_0(\nabla_\a u(t))\, dx_\a  & \leq & 2\int_\o \mathcal Q W_0(\nabla_\a v)\, dx_\a +2 \mathcal H^1(S(v) \setminus K_t^\eta )\\
 & \leq &  2\int_\o \mathcal Q W_0(\nabla_\a v)\, dx_\a +2 \mathcal H^1(S(v) \setminus \g(t) ) + 2\eta.
\end{eqnarray*}
The minimality property follows after letting $\eta \to 0$.\\

Concerning the convergence of the bulk energy, the previous calculation with $v=u(t)$ and the fact that $S(u(t)) \; \widetilde \subset \; \g(t)$ yield
\begin{equation}\label{lim+}
\limsup_{n \to +\infty}\int_\O W\left( \nabla_\a u^{\e_n}(t) \Big| \frac {1}{\e_n} \nabla_3 u^{\e_n}(t) \right)\, dx \leq 2 \int_\o \mathcal Q W_0(\nabla_\a u(t))\, dx_\a.
\end{equation}
Note that (\ref{lim+}) holds for the sequence $\{\e_n\}$ which is independent of the time. Thus, from (\ref{529}), we deduce that
$$\int_\O W\left( \nabla_\a u^{\e_{n_t}}(t) \Big| \frac {1}{\e_{n_t}} \nabla_3 u^{\e_{n_t}}(t) \right)\, dx \to 2\int_\o \mathcal Q W_0(\nabla_\a u(t))\, dx_\a.$$
In particular, Lemma \ref{stresses} implies the convergence of the
stresses and thanks to (\ref{CLcvforte}), (\ref{theta}) and
(\ref{sst}), we have for a.e. $t \in [0,T]$,
\begin{eqnarray*}
\theta(t) & = & \lim_{n_t \to +\infty} \int_\O \partial W \left(\nabla_\a u^{\e_{n_t}}(t) \Big|\frac{1}{\e_{n_t}} \nabla_3 u^{\e_{n_t}}(t) \right) \cdot \left(\nabla_\a \dot g^{\e_{n_t}}(t) \Big|\frac{1}{\e_{n_t}} \nabla_3 \dot g^{\e_{n_t}}(t) \right)dx\\
 & = & \int_\O \big( \partial (\mathcal Q W_0)(\nabla_\a u(t))|0 \big) \cdot \big( \nabla_\a \dot g(t)| \dot H(t) \big)\, dx\\
 & = & 2\int_\o \partial (\mathcal Q W_0)(\nabla_\a u(t)) \cdot \nabla_\a \dot g(t)\, dx_\a.
\end{eqnarray*}
\hfill$\Box$

\begin{rmk} {\rm
According to Remark \ref{2005} and Lemma \ref{MIN}, for every $t \in
[0,T]$, the function $u(t)$ minimizes
\begin{equation}\label{min1}
v \mapsto 2\int_\o \mathcal Q W_0(\nabla_\a v)\, dx_\a + 2 \mathcal
H^1 \big(S(v)\setminus \g(t)\big), \end{equation} among $\{v \in
SBV^p(\o';\Rb^3) : \quad v=g(t) \text{ a.e. on } \o' \setminus
\overline \o\}$. Equivalently, the pair $(u(t),\g(t))$ satisfies the
following unilateral minimality property:
\begin{equation}\label{min2}
2\int_\o \mathcal QW_0 (\nabla_\a u(t))\, dx_\a + 2 \mathcal H^1
(\g(t)) \leq 2\int_\o \mathcal QW_0 (\nabla_\a v)\, dx_\a + 2
\mathcal H^1 (\g'), \end{equation} for every $\mathcal
H^1$-rectifiable set $\g' \subset \overline \o$ such that $\g(t) \;
\widetilde \subset \; \g'$ and every $v \in SBV^p(\o';\Rb^3)$
satisfying $v=g(t)$ a.e. on $\o' \setminus \overline \o$ and $S(v)
\; \widetilde \subset \; \g'$. Indeed, for such pairs $(v,\g')$,
from (\ref{min1}) we get that
\begin{eqnarray*}
2\int_\o \mathcal QW_0 (\nabla_\a u(t))\, dx_\a  & \leq & 2\int_\o
\mathcal Q W_0(\nabla_\a v)\, dx_\a + 2 \mathcal H^1
\big(S(v)\setminus \g(t)\big)\\
 & \leq & 2\int_\o
\mathcal Q W_0(\nabla_\a v)\, dx_\a + 2 \mathcal H^1 \big(\g'
\setminus \g(t)\big)\\
 & = & 2\int_\o
\mathcal Q W_0(\nabla_\a v)\, dx_\a + 2 \mathcal H^1 (\g') - 2
\mathcal H^1 \big(\g(t)\big) \end{eqnarray*} where the second
inequality holds since $S(v) \; \widetilde \subset \; \g'$ and the
last equality because $\g(t) \; \widetilde \subset \; \g'$. On the
other hand, (\ref{min1}) follows from (\ref{min2}) by taking $\g' :=
S(v) \cup \g(t)$. }
\end{rmk}


\subsection{Energy conservation}\label{NRJconservation}

\noindent The last step in proving that $(u(t),\g(t))$ is a
quasistatic evolution relative to the boundary data $g(t)$ consists
in showing that the two-dimensional total energy $\E(t)$ defined in
(\ref{E(t)}) is absolutely continuous in time. This is the aim of
Lemmas \ref{E1} and \ref{E2} that follow.

\begin{lemma}\label{E1}
For every $t \in [0,T]$,
$$\E(t) \geq \E(0) +2 \int_0^t \int_\o \partial (\mathcal Q W_0)(\nabla_\a u(\tau)) \cdot \nabla_\a \dot g(\tau)\, dx_\a\, d\tau.$$
\end{lemma}

\noindent {\it Proof. }We proceed as in \cite{DMFT2} by approximation of the Lebesgue integral by Riemann sums. Let $s<t$, at time $s$ we test the minimality of $u(s)$ against $u(t)+g(s)-g(t)$. By Lemma \ref{MIN},
$$2\int_\o \mathcal Q W_0(\nabla_\a u(s)) \, dx_\a \leq 2\int_\o \mathcal Q W_0(\nabla_\a u(t)+\nabla_\a g(s)-\nabla_\a g(t)) \, dx_\a +2\mathcal H^1(S(u(t))\setminus \g(s)).$$
Thus, since $S(u(t)) \; \widetilde \subset \; \g(t)$ and $\g(s) \subset \g(t)$,
\begin{eqnarray*}
\E(s) & = & 2 \int_\o \mathcal Q W_0(\nabla _\a u(s))\, dx_\a +2 \mathcal H^1(\g(s))\\
 & \leq & 2\int_\o \mathcal Q W_0(\nabla_\a u(t)+\nabla_\a g(s)-\nabla_\a g(t)) \, dx_\a +2\mathcal H^1(\g(t))\\
 & = &  2\int_\o \mathcal Q W_0(\nabla_\a u(t)+\nabla_\a g(s)-\nabla_\a g(t)) \, dx_\a - 2\int_\o \mathcal Q W_0(\nabla_\a u(t)) \, dx_\a +\E(t).
\end{eqnarray*}
It implies that for some $\rho(s,t) \in [0,1]$,
\begin{equation}\label{abscont}
\E(t)-\E(s) \geq  2 \int_\o \left[ \partial (\mathcal Q W_0)
\left(\nabla_\a u(t) +\rho(s,t) \int_s^t \nabla_\a \dot g(\tau)\,
d\tau \right) \cdot \int_s^t \nabla_\a \dot g(\tau)\, d\tau
\right]\, dx.
\end{equation}
Fix $t \in [0,T]$, thanks to Lemma 4.12 in \cite{DMFT1}, there
exists a subdivision $0 \leq s_0^n \leq s_1^n \leq \ldots \leq
s_{k(n)}^n=t$ such that
$$\lim_{n \to +\infty}\sup_{1\leq i\leq k(n)} (s_i^n-s_{i-1}^n)=0$$ and
\begin{equation}\label{Riemann_sum}
\left\{\begin{array}{l}
\ds \lim_{n\to +\infty}\sum_{i=1}^{k(n)}\left\|(s_i^n-s_{i-1}^n)\nabla_\a \dot g(s_i^n)-\int_{s_{i-1}^n}^{s_i^n}\nabla_\a \dot g(\tau)\, d\tau \right\|_{L^p(\o;\Rb^{3 \times2})}=0\\
\ds  \lim_{n\to +\infty}\sum_{i=1}^{k(n)}\left|(s_i^n-s_{i-1}^n)\theta(s_i^n)-\int_{s_{i-1}^n}^{s_i^n} \theta(\tau)\, d\tau \right|=0.
\end{array}\right.
\end{equation}
For all $s \in (s_i^n,s_{i+1}^n]$, we define
$$u_n(s):=u(s_{i+1}^n), \quad \text{ and }\Psi_n(s):= \rho(s_i^n,s_{i+1}^n) \int_{s_i^n}^{s_{i+1}^n} \nabla_\a \dot g(\tau)\, d\tau.$$
As $\nabla_\a \dot g \in L^1(0,T;L^p(\o';\Rb^{3 \times 2}))$, we have
\begin{equation}\label{resid}
\|\Psi_n(s)\|_{L^p(\o';\Rb^{3 \times 2})} \to 0,
\end{equation}
uniformly with respect to $s \in [0,t]$. In (\ref{abscont}), we replace $s$ by $s_{i}^n$ and $t$ by $s_{i+1}^n$, then a summation for $i=0$ to $k(n)-1$ yields
$$\E(t)-\E(0) \geq 2 \int_0^t \int_\o \partial (\mathcal Q W_0)(\nabla_\a u_n(\tau)+\Psi_n(\tau)) \cdot \nabla_\a \dot g(\tau)\,  dx_\a \, d\tau.$$
From (\ref{resid}) and Lemma 4.9 in \cite{DMFT1}, we have for a.e. $\tau \in (0,t)$,
$$\left| \int_\o \partial (\mathcal Q W_0)(\nabla_\a u_n(\tau)+\Psi_n(\tau)) \cdot \nabla_\a \dot g(\tau)\,  dx_\a - \int_\o \partial (\mathcal Q W_0)(\nabla_\a u_n(\tau)) \cdot \nabla_\a \dot g(\tau)\,  dx_\a\right| \to 0.$$
Thus, according to (\ref{pg2}) together with Lebesgue's Dominated Convergence Theorem,
$$\int_0^t \left| \int_\o \partial (\mathcal Q W_0)(\nabla_\a u_n(\tau)+\Psi_n(\tau)) \cdot \nabla_\a \dot g(\tau)\,  dx_\a - \int_\o \partial (\mathcal Q W_0)(\nabla_\a u_n(\tau)) \cdot \nabla_\a \dot g(\tau)\,  dx_\a\right|\, d\tau \to 0.$$
Thus,
$$\E(t)-\E(0) \geq \limsup_{n \to +\infty} 2 \int_0^t \int_\o \partial (\mathcal Q W_0)(\nabla_\a u_n(\tau)) \cdot \nabla_\a \dot g(\tau)\,  dx_\a \, d\tau.$$
But in view of (\ref{Riemann_sum}), (\ref{pg2}) and H\"older's inequality,
\begin{eqnarray*}
&&\sum_{i=1}^{k(n)} \left| \int_\o \partial (\mathcal Q W_0)(\nabla_\a u(s_i^n)) \cdot \left( (s_i^n-s_{i-1}^n)\nabla_\a g(s_i^n)-\int_{s_{i-1}^n}^{s_i^n}\nabla_\a \dot g(\tau)\, d\tau \right)\, dx_\a \right|\\
& \leq & C\left( 1+\|\nabla_\a
u\|^{p-1}_{L^\infty(0,t;L^p(\o;\Rb^{3 \times 2}))}
\right)\sum_{i=1}^{k(n)} \left\|(s_i^n-s_{i-1}^n)\nabla_\a \dot
g(s_i^n)-\int_{s_{i-1}^n}^{s_i^n}\nabla_\a \dot g(\tau)\, d\tau
\right\|_{L^p(\o;\Rb^{3 \times2})}\\
&& \to  0,
\end{eqnarray*}
thus, using again (\ref{Riemann_sum}) and (\ref{theta0}),
\begin{eqnarray*}
\E(t)-\E(0) & \geq & 2 \limsup_{n \to +\infty}\sum_{i=1}^{k(n)}(s_i^n-s_{i-1}^n)\int_\o \partial (\mathcal Q W_0)(\nabla_\a u(s_i^n)) \cdot \nabla_\a \dot g(s_i^n)\, dx_\a \\
 & = & 2 \int_0^t \int_\o \partial (\mathcal Q W_0) (\nabla_\a u(\tau)) \cdot \nabla_\a \dot g(\tau)\, dx_\a \, d\tau.
\end{eqnarray*}
\hfill$\Box$\\

It now remains to show that the inequality proved in Lemma \ref{E1} is actually an equality. This is the object of the following Lemma.

\begin{lemma}\label{E2}
For every $t \in [0,T]$,
$$\E(t) \leq \E(0) +2 \int_0^t \int_\o \partial (\mathcal Q W_0)(\nabla_\a u(\tau)) \cdot \nabla_\a \dot g(\tau)\, dx_\a\, d\tau.$$
\end{lemma}

\noindent {\it Proof. }According to (\ref{sci}) and (\ref{529}),
\begin{eqnarray}\label{ineq1}
\liminf_{n_t \to +\infty} \E_{\e_{n_t}}(t) & = & \liminf_{n_t \to +\infty} \left[ \int_{\O} W\left( \nabla_\a u^{\e_{n_t}}(t) \Big| \frac {1}{\e_{n_t}} \nabla_3 u^{\e_{n_t}}(t) \right) \, dx \right.\nonumber\\
 & & \hspace{2.0cm}\left. + \int_{\G^{\e_{n_t}}(t)} \left| \left( \left( \nu_{\G^{\e_{n_t}}(t)}\right)_\a \Big| \frac{1}{\e_{n_t}} \left( \nu_{\G^{\e_{n_t}}(t)} \right)_3 \right) \right| \, d \mathcal H^2 \right]\nonumber\\
 & \geq & 2\int_\o \mathcal Q W_0(\nabla_\a u(t)) \, dx_\a + 2 \mathcal H^1 (\g(t))\nonumber\\
 & = & \E(t).
\end{eqnarray}
On the other hand, by Lemma \ref{MIN0}, (\ref{theta}), (\ref{scs})
and (\ref{theta0}) we have
\begin{eqnarray}\label{ineq2}
\limsup_{n_t \to +\infty} \E_{\e_{n_t}}(t) & \leq & \limsup_{n \to +\infty} \E_{\e_n}(t)\nonumber\\
& \leq &\lim_{n \to +\infty} \E_{\e_n}(0)+ \limsup_{n \to +\infty}
\int_0^t \theta_n(\tau)\, d\tau\nonumber\\
& = & \E(0) + \int_0^t \theta (\tau)\, d\tau \nonumber \\
& = &  \E(0) + \int_0^t \int_\o \partial (\mathcal Q W_0)(\nabla_\a
u(\tau))\cdot \nabla_\a \dot g(\tau)\, dx_\a \, d\tau.
\end{eqnarray}
Accordingly, relations (\ref{ineq1}) and (\ref{ineq2}) complete the
proof of the Lemma.
\hfill$\Box$\\

By virtue of Lemmas \ref{E1} and \ref{E2}, the two-dimensional total
energy $\E(t)$ is absolutely continuous with respect to the time $t$
and
$$\E(t) = \E(0) +2 \int_0^t \int_\o \partial (\mathcal Q W_0)(\nabla_\a u(\tau)) \cdot \nabla_\a \dot g(\tau)\, dx_\a\, d\tau$$
hence, $(u(t),\g(t))$ is a quasistatic evolution relative to the
boundary data $g(t)$. Let us show now that the three-dimensional
bulk and surface energies are converging towards the two-dimensional
bulk and surface energies respectively. Note that the following
convergence result holds for a subsequence $\{\e_n\}$ independent of
$t$ unlike in Lemma \ref{MIN} where we stated the convergence of the
volume energy for a $t$-dependent subsequence $\{\e_{n_t}\}$.

\begin{lemma}\label{convE}
For every $t \in [0,T]$,
$$\left\{\begin{array}{l}
\ds \int_\O W\left( \nabla_\a u^{\e_n}(t) \Big| \frac {1}{\e_n} \nabla_3 u^{\e_n}(t) \right)\, dx \to 2\int_\o \mathcal Q W_0(\nabla_\a u(t))\, dx_\a,\\
\\
\ds \int_{\G^{\e_n}(t)} \left| \left( \left( \nu_{\G^{\e_n}(t)} \right)_\a \Big| \frac{1}{\e_n} \left( \nu_{\G^{\e_n}(t)} \right)_3 \right) \right| \, d \mathcal H^2 \to 2 \mathcal H^1(\g(t)).
\end{array}\right.$$
In particular, $\E_{\e_n}(t) \to \E(t)$.
\end{lemma}

\noindent{\it Proof. }For $t=0$, the result is already proved in Lemma \ref{MIN0}. Assume now that $t \in (0,T]$ and let $\{n_j\}$ be a $t$-dependent subsequence such that
\begin{equation}\label{sssuite}
\liminf_{n \to +\infty}\int_\O W\left( \nabla_\a u^{\e_n}(t) \Big| \frac {1}{\e_n} \nabla_3 u^{\e_n}(t) \right)\, dx
=\lim_{j \to +\infty}\int_\O W\left( \nabla_\a u^{\e_{n_j}}(t) \Big| \frac {1}{\e_{n_j}} \nabla_3 u^{\e_{n_j}}(t) \right)\, dx.
\end{equation}
Arguing as in the proofs of Lemmas \ref{eto0} and \ref{MIN}, we can
suppose that, for a subsequence of $n_j$ (still denoted by $n_j$),
$u^{\e_{n_j}}(t) \rightharpoonup u^*(t)$ in $SBV^p(\O';\Rb^3)$ for
some $u^*(t) \in SBV^p(\o';\Rb^3)$ with $u^*(t)=g(t)$ a.e. on $\o'
\setminus \overline \o$, $S(u^*(t)) \; \widetilde \subset \;
\g(t)$, and which is also a minimizer of
$$v \mapsto 2\int_\o \mathcal QW_0(\nabla_\a v)\, dx_\a + 2\mathcal H^1(S(v) \setminus \g(t)),$$
among $\{v \in SBV^p(\o';\Rb^3): \quad v=g(t) \text{ a.e. on }\o' \setminus \overline \o\}$. Hence,
\begin{equation}\label{=}
\int_\o \mathcal QW_0(\nabla_\a u^*(t))\, dx_\a=\int_\o \mathcal QW_0(\nabla_\a u(t))\, dx_\a.
\end{equation}
According to Remark \ref{qcvx} and Theorem 5.29 in \cite{AFP},
\begin{eqnarray*}
2 \int_\o \mathcal QW_0(\nabla_\a u^*(t))\, dx_\a & \leq & \liminf_{j \to +\infty} \int_\O  \mathcal QW_0(\nabla_\a u^{\e_{n_j}}(t))\,dx\\
 & \leq & \lim_{j \to +\infty}\int_\O W\left( \nabla_\a u^{\e_{n_j}}(t) \Big| \frac {1}{\e_{n_j}} \nabla_3 u^{\e_{n_j}}(t) \right)\, dx.
\end{eqnarray*}
Thus, (\ref{sssuite}) and (\ref{=}) imply that
$$2 \int_\o \mathcal QW_0(\nabla_\a u(t))\, dx_\a \leq \liminf_{n \to +\infty}\int_\O W\left( \nabla_\a u^{\e_n}(t) \Big| \frac {1}{\e_n} \nabla_3 u^{\e_n}(t) \right)\, dx,$$
which ensure together with (\ref{lim+}) the convergence of the bulk energy i.e.
\begin{equation}\label{volume}
\int_\O W\left( \nabla_\a u^{\e_n}(t) \Big| \frac {1}{\e_n} \nabla_3 u^{\e_n}(t) \right)\, dx \to 2\int_\o \mathcal Q W_0(\nabla_\a u(t))\, dx_\a.
\end{equation}
But in view of (\ref{Etot}), (\ref{scs}) and Lemmas \ref{MIN0} and
\ref{MIN},
\begin{eqnarray}\label{NRJ}
\limsup_{n \to +\infty} \E_{\e_n}(t) & = & \limsup_{n \to +\infty} \E_{\e_n}(0) +\limsup_{n \to +\infty} \int_0^t \theta_n(\tau)\, d\tau\nonumber\\
 & \leq & \E(0) + \int_0^t \theta(\tau)\, d\tau\nonumber\\
 & = & \E(0) + 2 \int_0^t \int_\o \partial (\mathcal Q W_0)(\nabla_\a u(\tau)) \cdot \nabla_\a \dot g(\tau)\, dx_\a\, d\tau\nonumber\\
 & = & \E(t).
\end{eqnarray}
Thus (\ref{volume}) and (\ref{NRJ}) yield
$$\limsup_{n \to +\infty}\int_{\G^{\e_n}(t)} \left| \left( \left( \nu_{\G^{\e_n}(t)} \right)_\a \Big| \frac{1}{\e_n} \left( \nu_{\G^{\e_n}(t)} \right)_3 \right) \right| \, d \mathcal H^2 \leq 2 \mathcal H^1(\g(t))$$
which, together with (\ref{sci}), gives the convergence of the surface term
$$\int_{\G^{\e_n}(t)} \left| \left( \left( \nu_{\G^{\e_n}(t)} \right)_\a \Big| \frac{1}{\e_n} \left( \nu_{\G^{\e_n}(t)} \right)_3 \right) \right| \, d \mathcal H^2 \to 2 \mathcal H^1(\g(t)).$$
\hfill$\Box$\\

\noindent {\bf Acknowledgments. }The author is indebted to Gilles
Francfort for having proposed him this problem and for his fruitful
suggestions and comments. He also wishes to thank the referee for
his remarks and improvements.


\vspace{0.5cm}

\begin{center}
\begin{small}

Jean-Fran\c{c}ois Babadjian\\
 \textsc{L.P.M.T.M., Universit\'e Paris Nord, 93430, Villetaneuse, France}\\
\textit{E-mail address}: {\bf jfb@galilee.univ-paris13.fr}

\end{small}
\end{center}

\end{document}